\documentclass[11pt]{article}
\usepackage{amsfonts, epsfig, amsmath, amssymb, color, url}

 \usepackage[linkcolor=black,citecolor=black]{hyperref}
	
\usepackage[english]{babel}

\usepackage{paralist}
\usepackage{ulem}
\usepackage{enumitem} 
\usepackage[]{color}
\definecolor{DarkGreen}{rgb}{0.1,0.7,0.3}   

\definecolor{DarkGreen}{rgb}{0.1,0.7,0.3}   

\newcommand{\frack}[2]{\mbox{$\frac{#1}{#2}$}}
\newcommand{\spac}[2]{\hspace{#1}{#2}\hspace{#1}}
\newcommand{\eq}[0]{eqnarray}

\newcommand{\E}{\mathbf{E}}

\renewcommand{\P}{\mathbf{P}}

\renewcommand {\epsilon}{\varepsilon}

\newtheorem{theorem}{Theorem}
\newtheorem{defin}{Definition}[section]
\newtheorem{prop}{Proposition}[section]
\newtheorem{cor}{Corollary}[section]
\newtheorem{lemma}{Lemma}[section]

\newtheorem{rem}{Remark}[section]

\DeclareMathSymbol{\ophi}{\mathalpha}{letters}{"1E}

\newcommand{\refq}[1]{(\ref{#1})}

\newcommand{\e}{\varepsilon}

\newcommand{\be}{\begin{equation}}
\newcommand{\ee}{\end{equation}}
\newcommand{\ben}{\begin{equation*}}
\newcommand{\een}{\end{equation*}}

\newcommand{\ba}{\begin{equation}\begin{aligned}}
\newcommand{\ea}{\end{aligned}\end{equation}}
\newcommand{\ban}{\begin{equation}\begin{aligned*}}
\newcommand{\ean}{\end{aligned}\end{equation*}}

\newcommand{\mb}[2]{\hspace{#1}\mbox{#2}\hspace{#1}}

\DeclareMathOperator{\Law}{Law}
\DeclareMathOperator{\x}{x}
\DeclareMathOperator{\z}{z}

\newcommand{\cA}{\mathcal{A}}
\newcommand{\cB}{\mathcal{B}}

\newcommand{\cF}{\mathcal{F}}

\newcommand{\lime}{\lim_{\e\rightarrow0}}

\newcommand{\bR}{\mathbb{R}}
\newcommand{\bN}{\mathbb{N}}

\newfont{\cyrfnt}{wncyr10}
\def\J3{\cyrfnt{\rm \u{\cyrfnt I}}}
\def\j3{\cyrfnt{\rm \u{\cyrfnt i}}}

\allowdisplaybreaks[4]

\title{Large deviations for light-tailed L\'evy bridges on short time scales}

\author{Michael A. H\"ogele\\ 
Departamento de Matem\'aticas, Universidad de los Andes,
Bogot\'a, Colombia\\
\url{ma.hoegele@uniandes.edu.co}\\
\hfill\\Torsten \ Wetzel\\ 
Berlin, Germany\\
\url{w@zel.berlin}\\
}

\date{\today}

\begin{document}

\maketitle
\begin{abstract}
Let $L = (L(t))_{t\geq 0}$ be a multivariate L\'evy process with L\'evy measure $\nu(dy) = \exp(-f(|y|)) dy$ for a smoothly regularly varying function $f$ of index $\alpha>1$. The process $L$ is renormalized as $X^\e(t) = \e L(r_\e t)$, $t\in [0, T]$, for a scaling parameter $r_\e = o(\e^{-1})$, as $\e \to 0$. 
We study the behavior of the bridge $Y^{\e, \x}$ of the renormalized process $X^\e$ conditioned on the event $X^\e(T) = \x$ for a given end point $\x\neq 0$ and end time $T>0$ in the regime of small $\e$. 
Our main result is a sample path large deviations principle (LDP) for $Y^{\x, \e}$ with a specific speed function $S(\e)$ and an entropy-type rate function $I_{\x}$ on the Skorokhod space 
in the limit $\e \rightarrow 0+$. 
We show that the asymptotic energy minimizing path of $Y^{\e, \x}$ 
is the linear parametrization of the straight line between $0$ and $\x$, 
while all paths leaving this set are exponentially negligible. 
We also infer a LDP for the asymptotic number of jumps 
and establish asymptotic normality 
of the jump increments of $Y^{\e, \x}$. 
Since on these short time scales ($r_\e = o(\e^{-1})$) 
direct LDP methods cannot be adapted we use an alternative direct approach 
based on convolution density estimates of the marginals $X^{\e}(t)$, $t\in [0, T]$,
for which we solve a specific nonlinear functional equation. 
\end{abstract}

\section{Introduction}\setcounter{equation}{0}

Recall that a parameterized family of probability measures $(\P^\e)_{\e>0}$ on a topological space $\mathfrak{X}$ equipped with its Borel-$\sigma$-algebra $\mathcal{B}$ obeys a \textbf{large deviations principle (LDP) with speed function $S$ and rate function $I$}, 
if there is a function $S: (0, \infty) \to (0, \infty)$ satisfying $S(\e) \to 0$ as $\e \to 0+$ 
and an upper semicontinuous function $I: \mathfrak{X} \to [0, \infty]$ with compact sub-level sets 
such that for all $E\in \mathcal{B}$ 
\begin{align*}
-\inf_{z\in E^\circ} I(z) \leq \liminf_{\e \to 0+} S(\e) \ln \P^\e(E) \leq 
\limsup_{\epsilon \to 0+} S(\e) \ln \P^\e(E) \leq -\inf_{z\in \bar E} I(z).
\end{align*}
For precise definitions and references we refer to Subsection~\ref{sss:LDPs} below. 

This article establishes a LDP and asymptotic path properties for the bridges of a paradigmatic class of rescaled multidimensional L\'evy processes $X^\e_t:=\e L_{r_\e t}$ with light-tailed jump measures and some $r_\e=o(\e^{-1})$ as $\e\to 0+$. Those bridges connect the origin and a given end-point $\x\in \mathbb{R}^n$, $n\in \mathbb{N}$, $\x\neq 0$, similarly to the classical Brownian bridge (see e.g. \cite{KaratzasShreve}). The LDP includes a speed function $S$ which is given with the help of the solution of a particular functional equation, and a rate function $I_{\x}$ which strongly resembles the differential entropy of information theory. 

\noindent For convenience, we start the exposition of our results in a paradigmatic simplified setting. Consider a compound Poisson process $(L_t)_{t\geq 0}$ with values in $\mathbb{R}^n$ with the absolutely continuous, rotationally invariant, light-tailed jump measure given by 
\[
\nu(dz) = e^{-|z|^{\alpha}}dz, \qquad \alpha>1.
\] 
For any $\e>0$ and $\rho\in (-\infty, 1)$ we define the stochastic process $(X^\e_t)_{t\in [0, T]}$ by 
\[
X^\e_t = \e L_{t \cdot \e^{-\rho}}, \qquad t\in [0, 1].  
\]
Note that while $\P(X^\e_t = 0)>0$ for $t, \e>0$, the law of $X^\e_t$ outside the point mass in $0$ is absolutely continuous with respect to the Lebesgue measure in $\bR^n$ for any $t>0$ and $\e>0$. Denote by $\mu_t$ the density of $\Law(L_t)$ on $\bR^n\setminus\{0\}$.

\noindent Our main results (Theorem~\ref{th:1Rn}, ~\ref{Th:2Rn} and Theorem~\ref{Th:3Rn}$(i)$ and $(ii)$) are generalizations of the items 1.-4. of the following theorem in the respective order. 
\begin{theorem}\label{th:resumen}\hfill
\begin{enumerate}
 \item[1.] \textbf{Exponential density estimates of $\Law(L_t)$ by a nonlinear functional equation: } For all $n\in\bN$, $\alpha>1$, $\delta>0$ and $\rho < \gamma <1$ there is some $k>0$ such that for all $y\in \bR^n$ with $|y|>k$ and $t\in [|y|^\rho, |y|^\gamma]$ we have 
 \begin{align}
 -|y|\bigg(\alpha \Big(g\big(\frack {|y|}t\big)\Big)^{\alpha-1} -(1-\delta)g\big(\frack {|y|}t\big)^{-1}\bigg)
 \le\ln\mu_t(y) \le -|y|\bigg(\alpha \Big(g\big(\frack {|y|}t\big)\Big)^{\alpha-1}-(1+\delta)g\big(\frack {|y|}t\big)^{-1}\bigg),\label{e:fine}
 \end{align}
where for some $\Lambda_o>0$ the function $g:(\Lambda_o,\infty)\to\bR$ is given as the unique pointwise solution of the following nonlinear functional equation 
\begin{align}\label{e:funktional}
&\big(g(\Lambda)\big)^\alpha + C_1 \ln\big(g(\Lambda)\big) =  C_2 + C_3 \ln\big(\Lambda \big),\qquad \Lambda >\Lambda_0, 
\end{align}
\begin{align*}
&\mbox{where }\quad  C_1 = \frac{2-(\alpha-2) n}{2(\alpha-1)}, \quad C_2=\frac{\ln(\alpha-1)+n\ln(\alpha)-n\ln(2\pi)}{2(\alpha-1)},\qquad~~C_3 = \frac{1}{\alpha-1}.\nonumber
 \end{align*}
The function $g$ is slowly varying and $\Lambda\mapsto g(e^\Lambda)$ is a smoothly regularly varying function of order $\alpha^{-1}$.

 \item[2.] \textbf{Sample path LDP for the L\'evy bridges $Y^{\x, \e}$ with speed function $S$ and entropy-like rate function $I_{\x}$ : } For fixed $\x \in\mathbb{R}^n\backslash \{0\}$ we denote by $Y^{\e, \x} = (Y^{\e, \x}_{t})_{t\in [0,1]}$ the process $(X^\e_t)_{t\in [0,1]}$ being conditioned on the event $\{X^\e_1 = \x\}$. 
 
 Then the family $(\P^{\x, \e})_{\e>0}$, $\P^{\x, \e} = \Law(Y^{\x, \e})$ 
 satisfies a LDP on the space of c\`adl\`ag functions with values in $\bR^n$ 
 equipped with the uniform norm $(\mathcal{D}_{[0, 1],\mathbb{R}^n}, \|\cdot\|_\infty)$ 
 with speed function $S(\e)$ given by 
 \begin{equation}\label{e:simplerate}
 S(\e) := \e \cdot g(\e^{-(1-\rho)}), \quad \e>0,   
 \end{equation}
 where $g$ is defined by the solution of \eqref{e:funktional}, 
 and the good rate function 
 \[
 I_{\x}(\varphi) := 
 \begin{cases}  
\displaystyle \int_0^1 |\varphi|' ~ \ln |\varphi|' \,dt- |\x| \ln |\x|, & \mbox{if } |\varphi|(\cdot)\mbox{ is }
 \begin{cases}
 \mbox{continuous,}&\\
 \mbox{nondecreasing,}&\\
 \varphi([0,1])=[[0,\x]],
 \end{cases}\\
   \infty, & \mbox{ otherwise.}
  \end{cases}
 \]
Here we denote by $|\varphi|(t) = |\varphi(t)|$, $|\varphi|'(t) = \frac{d}{dt} |\varphi(t)|$ the total derivative, whenever it exists and set it equal to $0$ otherwise. We set $r \ln r = 0$, whenever $r = 0$, and $[[0, \x]] = \{\vartheta\x\in \bR^n~|~\vartheta\in [0, 1]\}$. 
 \item[3.] \textbf{LDP for the jump frequency of the L\'evy bridges $Y^{\x, \e}$: } For $g$ being defined by \refq{e:fine} we denote by $N^{\x, \e}$ the number of jumps of $Y^{\x, \e}$, and 
  \[
\bar N^{\x, \e}:= \sqrt{k_\e} (N^{\x, \e}- m_{\x, \e}) \qquad \mbox{ for }\qquad m_{\x,\e} := g\Big(\frac{|\x|}{\e^{1-\rho}}\Big)^{-1} \frac{|\x|}{\e}\quad \mbox{ and }\quad k_\e := \alpha  \frac{\e}{|\x|} g\Big(\frac{|\x|}{\e^{1-\rho}}\Big). 
 \]
Then, for $S$ being defined by \eqref{e:simplerate}, the family $(\mathbf{Q}^{\x,\e})_{\e>0}$, $\mathbf{Q}^{\x,\e}:= \Law(\sqrt{S(\e)} \bar N^{\x, \e})$, satisfies a LDP on $(\mathbb{R}, |\cdot|)$ with speed function $S(\e)$ and good rate function
\[
J(y):= \frac{1}{2} y^2. 
\]
\item[4.] \textbf{Asymptotic normality of the jump increments of the L\'evy bridges $Y^{\x, \e} / \e	$: } We keep the notation of 3., denote by $W^{\x, \e}_i$ the $i$-th jump of the process $Y^{\x,\e}_t / \e$ and define
\begin{align*}
\bar W^{\x, \e} 
&:= \langle W^{\x, \e}_1, \x\rangle \frac{\x}{|\x|} + (\alpha-1)^{-1} \Big(W^{\x, \e}_1 -\langle W^{\x, \e}_1, \x\rangle \frac{\x}{|\x|}\Big).
\end{align*}

Then the family $(\mathbb{W}^{\x, \e})_{\e>0}$ defined by 
\[
\mathbb{W}^{\x, \e} := \sqrt{\alpha(\alpha-1) g(\e^{-(1-\rho)})^{\alpha-2}} \cdot 	\Big(\bar W^{\x, \e} - \frac{\x}{\e m_{\x, \e}}\Big)  
\]
converges in distribution as $\e\to 0+$ to a standard normal random vector in $\mathbb{R}^n$. 
\end{enumerate}
 
\end{theorem}
In particular, item~1. of Theorem \ref{th:resumen} turns out to be a corollary of Theorem~\ref{th:1Rn}, item 2. is a direct consequence of Theorem~\ref{Th:2Rn}, while items 3. and 4. are implications of Theorem~\ref{Th:3Rn} below. We proceed with the discussion of the results first in the special case formulated in Theorem \ref{th:resumen}. In item~f. below the generalizations made in our main results, Theorems \ref{th:1Rn}, \ref{Th:2Rn} and \ref{Th:3Rn} are explained. 

\begin{enumerate}
 \item[a.] Recall that $X^\e_t = \e L_{t\cdot \e^{-\rho}}$. While the case of scaling with $\rho = 1$ is covered and well-studied by classical LDP methods (see item 1. in Subsection~\ref{sss:glance} below), to our knowledge, this article is the first one 
 to study LDPs rigorously for asymptotically ``short'' time scales $\e^{-\rho}, \rho< 1$. 
 The case $\rho= 0$ represents the bridge of unscaled L\'evy process for small amplitude $\e$. That is to say, $\rho=0$ treats the situation of $(\e L_t)_{t \in [0, 1]}$ being conditioned to reach $\x\neq 0$ in time $1$. This case is one of the main motivations for our studies, we refer to item 2.a in Subsection~\ref{sss:glance}. Note that for $\rho<0$ we have $\lim_{\e\to 0} \P(X^\e_t \neq 0) = 0$ for any $t\ge0$. Thus it might come as a surprise that, despite the preceding convergence, we obtain results for all $\rho \in (-\infty, 1)$, with no difference in treatment between $\rho<0$, $\rho=0$ and $\rho\in (0,1)$.  

 \item[b.] The function $g$ defined as the solution of equation \refq{e:funktional} is essential in the presentation of our results. Even though the definition of $g$ seems to be slightly involved, its asymptotic behavior can easily  be read off from \refq{e:funktional}:
\be\label{14}
 \lim_{\Lambda\to \infty } \frac{g(\Lambda)}{\ln(\Lambda)^\frac{1}{\alpha}} = (\alpha-1)^{-\frac{1}{\alpha}}. 
\ee
 
 \item[c.] Note that in Theorem~\ref{th:resumen},  item 2., $I_{\x}(\varphi) = \infty$ for all $\varphi: [0,1]\to \bR^n$ with $\varphi([0,1]) \neq [[0, \x]]$. In this sense the sample path LDP is rather degenerate. 
Among all continuous and non-decreasing parametrizations of the straight-line segment $[[0, \x]]$ the minimizing paths of $I_{\x}$ on $\mathcal{D}_{[0, 1],\mathbb{R}^n}$ is easily seen to be the linear function $\varphi(t) = t \x, t\in [0, 1]$, connecting $0$ with $\x$ with constant velocity. 

  \item[d.] All jumps of $Y^\e/\e$ are identically distributed, while clearly not being independent. For $\x = (\x_1, 0, \dots, 0)$ and $\alpha > 1$ we have that for small values of $\e$ the distribution of a single jump $W_1^{\x, \e}$ is close to 
 \begin{align}\label{e:aspnormy}
\mathcal{N}\left(\frac{\x}{|\x| }g\Big(\frac{|\x|}{\e^{1-\rho}}\Big) ,
 \frac{1}{\alpha(\alpha-1)}g(\e^{-(1-\rho)})^{2-\alpha} 
 \left[\begin{array}{cc} 
 1 & 0\\ 
 0 & (\alpha-1) I_{n-1}      
\end{array}\right]
 \right), 
 \end{align}
where $I_{n-1}$ is the identity matrix of size $(n-1)\times (n-1)$. 
Recall that $\rho<1$, therefore the mean jump size of $Y^\e / \e $ diverges for all $\alpha>1$, as $\e\to0$. Furthermore by the choice $\alpha>1$ the exponent $2-\alpha$ in the variance of \refq{e:aspnormy} is smaller than 2. This implies that $\E[|W^{\x, \e}|] \approx \big| \E[W^{\x, \e}]\big|$ for small $\e$. Thus by \eqref{e:aspnormy} and the definition of $g$ and its asymptotics \refq{14}, the mean jump size satisfies the relation
 \begin{equation}\label{e: asymptotic mean increment}
\E[|W^{\x, \e}|] \approx  g\Big(\frac{|\x|}{\e^{1-\rho}}\Big) \approx 
 (\alpha-1)^{-\frac{1}{\alpha}} \cdot \Big(\ln \frac{|\x|}{\e^{1-\rho}}\Big)^{\frac{1}{\alpha}}, \qquad \mbox{ as } \e \to 0.  
 \end{equation}
Note that for $\alpha>2$ the variance in the direction of $\x$ is smaller than in the directions of its orthogonal compliment $\x^\perp$. In the case $\alpha\in(1,2)$ this relation is inverted, and for $\alpha=2$ the law in \eqref{e:aspnormy} has asymptotically a spatially homogeneous covariance matrix.

\item[e.] Since $I_{\x}(\varphi) = \infty$ for all paths $\varphi$ leaving the segment $[[0, \x]]$, we expect that on average $|W_1^{\x, \e}| \cdot N^{x, \e} \approx |\x| / \e$. Indeed, with the help of \eqref{e: asymptotic mean increment} and by Theorem~\ref{th:resumen}, item 3., the number of jumps satisfies asymptotically, as $\e\to 0$, 
\[
\E[ N^{\x, \e}] \approx m_{\x, \e} = \frac{|\x|}{\e} \cdot g\Big(\frac{|\x|}{\e^{1-\rho}}\Big)^{-1} 
\approx \frac{|\x|}{\e}\cdot \E[|W_1^{\x, \e}|]^{-1}.
 \]
The LDP of Theorem~\ref{th:resumen}, item 2., can be understood in the light of item 3. and 4., 
in that with very high probability the jumps of trajectories of $Y^\e$ resemble longer and longer sequences of smaller and smaller increments of order proportional to $\e (\ln\e^{-(1-\rho)})^\frac{1}{\alpha}$ with a strong prevalence to align in direction $\x/|\x|$. This illustrates the invariance of the optimal path on the segment $[[0,\x]]$ as discussed in item c.
 \item[f.] The results of Theorem~\ref{th:resumen} are generalized in Theorem~\ref{th:1Rn},~\ref{Th:2Rn} and~\ref{Th:3Rn} as follows: 
 \begin{enumerate}
  \item[i.] Instead of $\nu(dx) = e^{-|x|^\alpha} dx$, $\alpha>1$, 
  we allow for $\nu(dx) = e^{-f(|x|)}dx$ for a smoothly regularly varying function $f$ of index $\alpha>1$, see Definition~\ref{def:SR} below.  
  \item[ii.] Instead of the condition $\{\e L_{\e^{-\rho}} = \x\}$ for some $\rho <1$ we allow for $\e^{-\rho}$ being replaced by some regularly varying function $r_\e T$ with index $\rho>-1$ and some final time $T>0$. 
  \item[iii.] Our results of Theorem~\ref{th:1Rn} and \ref{Th:2Rn} remain valid in the presence of a Brownian motion  or a deterministic drift in the L\'evy process $L$. 
  \item[iv.] Our findings of Theorem~\ref{th:1Rn} and \ref{Th:2Rn} are robust under the presence of an independent perturbation of $L$ by an additional  L\'evy process $\xi$ with L\'evy measure $\nu_\xi$, as long as the densities of the sum are dominated by the density of the compound Poisson component in the sence of Hypothesis~II below. We stress that there the jump measure of $\xi$ neither needs to be rotationally invariant nor absolutely continuous with respect to the Lebesgue measure in~$\bR^n$. A set of sufficient conditions is established in Subsection~\ref{ss:SC}. For example, in the case of bounded jumps, it is sufficient that $\nu_\xi$ satisfies the Orey condition~\eqref{eq:Orey}. 
  
 \end{enumerate}
\end{enumerate}

\medskip 
\subsection{Motivation: a LDP for $\e L_{r_\e t}$ on short time scales $r_\e = o(\e^{-1})$}\label{sss:glance}

The large deviations principle is a powerful and well-understood concept 
to describe the precise asymptotic exponential decay rate of probabilities in terms of an optimization problem, understood in the physics literature as a generalized least action principle. For an introduction and a short overview over the history of the LDP we refer to the closing section of \cite{Var08}. Standard texts include \cite{AcostaAraujoGine, BorovkovMogulskii, Cramer, Dembo, Dupuis, Ellis, DeuStr89, DE97, El06, Feng, Freidlin, Gaertner, dH00, OlivieriVares05, RASep15, Varadhan,  Var67, Var84} among others. In particular, we stress the pioneering work of Freidlin and Wentzell \cite{Freidlin}, where the concept of the LDP is used first to describe the effect of perturbations of differential equations in the small noise limit in terms of LDPs, see also \cite{bovier04, CerraiSalins, CerraiFreidlinSalins17, DV75, Siegert}. 

Nowadays, it is well-known in the literature how to implement the sample path LDP for a parameterized family of L\'evy process $(Z^\e)_{\e>0}$, $Z^\e = (Z^\e_t)_{t\in [0,\infty)}$, with values in $\bR^n$ which satisfy the LDP for each marginals $\mbox{Law}(Z_t)$ with (good) rate function~$I_t$.  The contraction principle \cite[Theorem 4.2.1]{Dembo} yields directly, that for every $m\in\bN$, $m\ge2$, and every partition $0\le t_1<\dots<t_m$ the family $(Z^\e_{t_1},\dots Z^\e_{t_m})_{\e>0}$ satisfies the finite dimensional LDP with the rate function
\be
I_{(t_1,\dots,t_m)}(x_1,\dots,x_m):=I_{t_1}(x_1)+\sum_{i=2}^mI_{t_{i}-t_{i-1}}(x_i-x_{i-1}).\label{eq:1}
\ee
Finally, it remains to prove exponential tightness of  $(Z^\e)_{\e>0}$ to obtain the LDP for the family of processes on the path space of càdlàg functions equiped with the usual $\mathcal J_1$-topology. By Theorem~4.28 in Feng and Kurtz \cite{Feng} the rate function is given by 
\be
I(\varphi):=\sup_{\substack{0\le t_1<\dots < t_m\\\varphi\mbox{ is continuous in }t_1,\dots,t_m}}I_{(t_1,\dots,t_m)}(\varphi(t_1),\dots,\varphi(t_m)).\label{eq:2}
\ee
By \refq{eq:1} and \refq{eq:2} for any family of stochastic processes exhibiting the Markov property and stationary increments, the existence of a LDP boils down to a LDP of the marginals for fixed time and exponential tightness. We stress that \refq{eq:2} is valid for any finite dimensional LDP with respective rate functions $I_{(t_1,\dots,t_m)}$ not only of type \refq{eq:1}. \\
Since this paper covers Lévy bridges $Y^\e$ with fixed time horizon $T>0$ we have to adapt formula~\refq{eq:2} to the setting. Note that if we simply restricted the supremum in \refq{eq:2} to $t_m\le T$, then for any $\varphi,\phi\in\mathcal D_{[0,T],\bR^n}$ with $\phi(t)=\varphi(t)\textbf{1}_{[0,T)}(t)$ we would obtain $I(\varphi)\ge I(\phi)$, which is absurd. To circumvent this technicality we define by $\tilde Y^\e_t:=Y^\e_{t\wedge T}$ a family of random processes $(\tilde Y^\e)_{\e>0}$ on $\mathcal D_{[0,\infty),\bR^n}$   and apply \refq{eq:2} to $\tilde Y^\e$.  The contraction principle allows then to return to the original bridges $(Y^\e)_{\e>0}$.\\[2mm]
In the context of L\'evy processes the results in the literature for such marginals can be categorized into two cases: the 'sample mean' or 'inverse proportional' time scale $r_\e = \e^{-1}$, and 'large' time scales $r_\e\gg \e^{-1}$.\\[2mm]
\textbf{1. Sample mean time scales $r_\e = \e^{-1}$: }  We consider  the underlying L\'evy process and thus the jump sizes being multiplied by $\e$, while the times cale and thus the jump intensity is accelerated by a factor $r_\e=\e^{-1}$,  that is, $X^\e_t:=\e L_{t \e^{-1}}$. For the sake of argument assume $\e^{-1}$ to be integer-valued. Due to the infinite divisibility we can represent $X^\e_t$ for each $t>0$ as a sum of $\e^{-1}$ many i.i.d. summands
	\be
	X^\e_t=\e\sum_{i=1}^{\e^{-1}}(L_{it}-L_{(i-1)t}).\label{eq:3}
	\ee
Assume the existence of some finite exponential moments of $L_t$. By \refq{eq:3} the LDP can be understood as the asymptotic rate of convergence in the weak law of large numbers, that is, the sample mean regime. 
 The classical Cram\'er theorem \cite[Section 2.2]{Dembo} yields, that for every $t>0$ the family of $\mbox{Law}(X^\e_t)_{\e>0}$ satisfies a LDP. Moreover, the associated rate function is identified in terms of the Fenchel-Legendre transform of the logarithmic moment generating function of the increment distribution. If in this case exponential tightness of the family of processes $\mbox{Law}(X^\e)_{\e>0}$ can be proven, then $\mbox{Law}(X^\e)_{\e>0}$ satisfies a LDP and equations \refq{eq:1} and \refq{eq:2} allow to identify the sample path rate function. 
	
First results of this type were obtained by Lynch and Sethuraman \cite{Lynch} for real valued L\'evy processes. De Acosta \cite{DeAcosta} lifted those results to L\'evy processes on general Banach spaces. Nowadays, there is a large and fast growing literature on LDPs for L\'evy driven stochastic ordinary and partial differential equations which we cannot review here completely. Reference articles are for instance 
\cite{BrzezniakPengZhai,BudhirajaDupuisSalins, LiuSongZhaiZhan23,  MatoussiSabbaghZhang21,RocknerZhang07, Swiech, WuZhai20, WuZhai24, XiongZhai18, XuZhang09, YangZhaiZhang15, ZhaiZhang15} and the references therein. 
		
A different string of recent works initiated by A. Budhiraja and collaborators established the LDP in settings of perturbations random perturbations by accelerated Poisson random measures \cite{Budhiraja, BCD13,BD19, BDM11, BN15, AndreCatuogno, Hoegele}. 
	Here the jump sizes are multiplied by $\e$ while the jump intensity is accelerated (on average) by $\e^{-1}$. 
	Obviously, those processes do not necessarily exhibit stationary increments nor a representation as in equation~\refq{eq:3}. 
	Nonetheless, the jump sizes are multiplied by $\e$ while the time is accelerated by $\e^{-1}$. 
	Therefore even those processes can be considered as of the same spirit as \eqref{eq:3}. 
	Recent results in \cite{BazbhaBlanchetRheeZwart19} yield a LDP for scalar L\'evy processes with one-sided Weibull-type increments of index $\alpha \in (0,1)$, with an explicit rate function given as the $\alpha$-variation of the sample paths. For studies of heavy-tailed exit problems, such as $\alpha$-stable perturbations, we refer to \cite{Debussche, HP13, ImkPav06, Pav11}.
	
    \textbf{2.a Large time scales $r_\e\gg \e^{-1}$ and $\e$ spatial scaling: } The second case 
    consists of the underlying L\'evy process being multiplied by $\e$, while the time scale and thus the jump intensity is accelerated by a factor $r_\e$, that is asymptotically larger than $\e^{-1}$ and smaller than $\e^{-2}$, that is,  
	\begin{equation}\label{e:LDPscale}
	\lime\e r_\e=\lime(\e^2r_\e)^{-1}=\infty,
	\end{equation} 
	think of $r_\e = \e^{-\frac{3}{2}}$, for example. By definition we have $\E[Y_{t_\e}]= \e r_\e \E[L_1]$. Thus under \eqref{e:LDPscale} we need to assume  $\E[L_1] = 0$. To find a similar representation to \eqref{eq:3} in this case, we obtain either a number $r_\e$ of summands with $r_\e\gg \e^{-1}$, or, if we set the number of summands to $\e^{-1}$, then each of these summands can be represented as a sum of a diverging number of i.i.d. summands. We choose the latter approach. Here, we obtain a sum of $\e^{-1}$ summands and can approximate the distribution of each of these summands by the central limit theorem. If we now apply Cram\'er's theorem, we obtain a LDP which corresponds to the Brownian case. This is seen as follows.

For the sake of argument, set $r_\e := \e^{-\rho}$, $\rho \in (1,2)$. Furthermore, we consider $\e^{-1}$ and $\e^{-\frac{2\rho-2}{2-\rho}}$ to be integer valued. If $k_\e:= \e^{\frac{1}{2-\rho}}$ and $Y^\e:= X^{k_\e}$ we obtain 
	\begin{equation}\label{e:1_12}	
	Y^\e_t = \e^{\frac{1}{2-\rho}} L_{\e^{-\frac{\rho}{2-\rho}} t} = \e\sum_{i=1}^{\e^{-1}} \e^{\frac{\rho-1}{2-\rho}} \Big(L_{i \e^{-\frac{2\rho-2}{2-\rho}}t}-L_{(i-1) \e^{-\frac{2\rho-2}{2-\rho}}t}\Big).
	\end{equation}
	By the choice of $\rho$ we have that $\frac{2\rho-2}{2-\rho}>0$. Thus, as $\e\to0$, the central limit theorem for each of the summands yields 
	\be
	\e^{\frac{\rho-1}{2-\rho}} \Big(L_{i \e^{-\frac{2\rho-2}{2-\rho}}t}-L_{(i-1) \e^{-\frac{2\rho-2}{2-\rho}}t}\Big)
	\stackrel{d}{=} \e^{\frac{\rho-1}{2-\rho}} L_{\e^{-\frac{2\rho-2}{2-\rho}}t}= \e^{\frac{\rho-1}{2-\rho}} \sum_{j=1}^{\e^{-\frac{2\rho-2}{2-\rho}}} \Big(L_{jt}-L_{(j-1)t}\Big) \stackrel{d}{\longrightarrow} \mathcal{N}(0, \sigma^2\,t),\nonumber
		\ee
	where $\sigma^2$ denotes the variance of $L_1$. An application of Cram\'er's theorem to \eqref{e:1_12} implies the LDP for each $(Y^\e_t)_{\e>0}$, $t\in [0, T]$. If exponential tightness of the family $((Y^\e_t)_{t\in [0, T]})_{\e>0}$ can be established, then $((Y^\e_t)_{t\in [0, T]})_{\e>0}$ satisfies a LDP. By \eqref{eq:1} and \eqref{eq:2} we obtain the associated rate function, which obviously corresponds to the one in the Brownian case. By definition we have $X^\e = Y^{\e^{2-\rho}}$. Thus $(X^\e)_{\e>0}$ satisfies a LDP with the same rate function and the speed function $S(\e) = \e^{2-\rho}$. 

 First results of this type can be found in Mogulskii \cite{Mogulskii}, where a different parametrization to the one in item 1. is applied. There, for a parameter $T>0$ the author chooses $r(T)$ such that $\lim_{T\rightarrow\infty}r/T=0$ and $\lim_{T\rightarrow\infty} r/\sqrt {T}=\infty$. Furthermore, the author defines $\lambda= r^2/ T$ and $X^\lambda_t=r^{-1} L_{Tt}$ and finally establishes a LDP for $(X^\lambda)_{\lambda>0}$ with $\lambda\rightarrow\infty$. It is not difficult to see that Mogulskii's restrictions $\lim_{T\rightarrow\infty} r/T=0$ and $\lim_{T\rightarrow\infty} r/ \sqrt T =\infty$ are equivalent to the conditions~\eqref{e:LDPscale}. 

\textbf{2.b Large time scales $r_\e\gg \e^{-1}$ and different spatial scalings: the moderate deviation principle. } In situations where a different renormalization of the process is natural, often a so-called moderate deviation principle (MDP) can be derived. We refer to \cite{EichelsbacherLoewe03} for an introduction. To show this connection let $L$ be a (not necessarily centered) L\'evy process with finite exponential moments. Define $Y^\e_t:=\e L_{\e^{-1}t}$. Let $s_\e$ satisfy $\lime s_\e =\infty$ and $\lime\sqrt\e s_\e =0$. A typical MDP approach would be to study the behavior of $Z^\e_t:=s_\e (Y^\e_t- \E Y^\e_t)$. For convenience we restrict ourselves to the choice $s_\e =\e^{-c}$ with $c\in(0,\frac12)$. Let $(X^\e_t)_{t\in[0,T]}$ be defined by $X^\e_t=\e (L_{r_\e t}-\E L_{r_\e t})$, where $r_\e:=\e^{-\frac1{1-c}}\gg \e^{-1}$. On the one hand $(X^\e)_{\e>0}$ clearly belongs in this setting to case 2.a and on the other hand by definition we have $Z^\e = X^{\e^{1-c}}$. Therefore any MDP found for a the rescaled L\'evy process $(Z^\e)_{\e>0}$  can be used to describe the asymptotic of $(X^\e)_{\e>0}$ of item 2.a by a suitable reparametrization.

Recent MDP results cover much more sophisticated processes than the one used in the example above. For example, in Budhiraja et. al. \cite{Budhiraja} the authors find a MDP for $Z^\e_t:= s_\e (Y^\e_t-\E Y^\e_t)$, where $Y^\e$ is defined as the solution of stochastic differential equations and $s_\e$ is an appropriate scale function. As in \cite{Budhiraja, BCD13,BDM11} these processes are driven by a rescaled Poisson random measure. Thus, these processes are of the same spirit as our setting of item 2.a. MDP results naturally appear in multiscale dynamics such as singular forward backwards systems such as for instance \cite{AndreCatuogno, GasteratosSalinsSpiliopoulos23} and the substantial literature cited there, which goes beyond the scope of this introduction. 

\textbf{3. The missing case for jump L\'evy processes: short time scales $r_\e \ll \e^{-1}$. } 
We recall Schilder's celebrated theorem \cite{Dembo} for Brownian motion $(\e B_t)_{t\in [0, T]}$ which yields a LDP with speed function $S(\e) = \e^2$ 
\[
I_T(\phi) = \begin{cases} \frac{1}{2} \int_0^T |\phi'(s)|^2 ds& \mbox{ if } \phi'\in L^2([0, T]), \\ \infty & \mbox{ else.}\end{cases} 
\]
This clearly corresponds to the situation of $r_\e = 1$. 
However, there are no counterparts of this result for pure jump L\'evy processes of this type $(\e L_t)_{t\in [0, T]}$ established in the literature. To our knowledge there are even no LDP results in the literature for ``short'' time scales $r_\e$, ``short'' in comparison to the inverse spatial scaling $\e^{-1}$, $\lime\e r_\e=0$. This article proposes to fill this gap in the literature with a paradigmatic case study for L\'evy bridges with a light-tailed jump component. It contains a multidimensional generalization of LDP results of the recent Ph.D. thesis \cite[chapter 4]{Wetzel2} by the second author. This is carried out by a direct approach via asymptotic density estimates for the tails of the marginal distributions. A preliminary, and coarser kind of tail-estimates for light-tailed L\'evy processes had been established by the second author in \cite{Wetzel} and was published in \cite{Imkeller} in the context of exit times results. A detailed review is given after formula \refq{eq:1.11} in the introduction below. 

In the sequel we discuss the opportunities and limitations for LDP results for rescaled jump L\'evy processes in the short time regime. Consider $\e\mapsto r_\e$ to be of regular variation with an index in $(-1,\infty)$, such as $r_\e=\e^{-\rho}$, $\rho < 1$. In other words, the time scale is  not sufficiently large in comparison to the spatial scaling in order for classical theorems to be valid, as can be seen below. We continue to sketch why a straightforward LDP approach combined with the contraction principle is not successful in this setting.  With the help of the result \refq{e:fine} from part 1 in Theorem~\ref{th:resumen} it can be shown that
\be
\lime\e\ln\P(|X^\e|>\kappa)=-\lime\kappa g(\kappa\e^{-1}r_\e^{-1})=-\infty
\ee
for any $\kappa>0$. Hence the family $(X^\e_t)_{\e>0}$ satisfies the trivial LDP with the rate function $I_t(0)=0$ and $I_t(y)=\infty$ for any $y\neq0$. By the choice $r_\e=o(\e^{-1})$ this coincides with the following application of Cramér's theorem. A representation similar to \refq{eq:3} reads
\be
X^\e_t=\e\sum_{i=1}^{\e^{-1}}(L_{ir_\e\e t}-L_{(i-1)r_\e\e t}),
\ee
where each of these summands $(L_{ir_\e \e t} - L_{(i-1)r_\e \e t}) \stackrel{d}{=}L_{r_\e \e t} \to 0$ almost surely as $\e \to 0$. Therefore Cramér's theorem leads to the same void rate function.

Apparently, one might try to circumvent this obstacle by the use of a modified speed function~$S$. Indeed, with the choice $S(\e)=g(\e^{-1}r_\e^{-1})\e$ from \refq{e:fine} it follows, that
\be\label{e:fakeLDP}
\lim_{\kappa\to0}\lime S(\e)\ln\P(|X_t^\e-y|<\kappa)=-|y|
\ee
and for fixed $t>0$ the family $(X^\e_t)_{\e>0}$ satisfies a LDP with respect to the speed function $S$ and the rate-function  $I_t(y)=|y|$. However, in Appendix~\ref{s:LDPnonexistence} we show that for any scale-function~$S$, such that the family of marginals $\mbox{Law}(X_t^\e)_{\e>0}$ satisfies such a LDP, it follows that the family of processes $(X^\e)_{\e>0}$ is not $S-$exponentially tight on the path space of càdlàg functions $\mathcal D_{[0,T],\bR^n}$ equipped with the $J_1$-topology. Therefore any standard approach to establish a LDP for $(X^\e)_{\e>0}$ looks bound to fail in the case of $r_\e = o(\e^{-1})$. 

Due to the outlined structural challenges to establish a LDP for rescaled L\'evy process 
it seems natural to try to show a LDP for modified (rescaled) L\'evy processes. An obvious candidate for such an enterprise are rescaled L\'evy bridges. 
The density of the Brownian bridge was established as early as 1931 in the seminal article by Schr\"odinger  \cite{Schroedinger} and has been object of study ever since, for an overview see for instance \cite[Section 5.6.B]{KaratzasShreve}. For the rigorous construction of Markov bridges we refer to \cite{ChaumontUribeBravo,FitzsimmonsPitmanYor}. 
An overview on the respective reciprocal classes and the duality formulas is given in \cite{ConfortiLeonarRudigerRoelly, LeonardRoellyZambrini}. The LDP for Brownian bridges on manifolds was established in \cite{Hsu} with  further details in \cite{Wittich}. 
On Euclidean space  Brownian bridges satisfy the standard LDP due to \cite{BarczyKern13, GasbarraSottinenValkeila07} which had been generalized recently to Bernstein bridges \cite{PrivaultYangZambrini}. As in the case of Schilder's theorem, the results for pure jump processes look quite different. A LDP for a class of rescaled symmetric scalar L\'evy processes for bounded jumps 
in the sample mean regime 1. has been established in \cite{Yang14}. We refer to the end of Subsection~\ref{sss:approach} for a more detailed comparison with our results. 

\bigskip 
\subsection{A direct approach: convolution density estimates}
\label{sss:approach}

In order to circumvent the difficulties sketched in item 3.~of Subsection~\ref{sss:glance} and 
to take full advantage of the L\'evy bridge structure, 
we choose a direct approach. It consists of three conceptual steps, with tight estimates on the compound Poisson marginal density  in step I below, which we apply subsequently to the respective bridge marginal densities in step II. The finite dimensional distributions are obtained by similar techniques. Step III treats the necessary tightness results in order to apply \refq{eq:2}. 

\textbf{I.} We first consider a L\'evy process $L$ without neither a Brownian component nor a deterministic drift and the L\'evy measure $\nu$ is taken rotationally invariant, absolutely continuous and finite. That is, we treat the compound Poisson process $(L_t)_{t\in [0, 1]}$ with jump measure $\nu(dz) = e^{-f(|z|)} dz$ with $\nu(\mathbb{R}^n) = 1$ for some function $f$. We assume that $f$ is a smoothly regularly varying function with index $\alpha>1$. The concept of (smooth) regular variation is described in detail in N. H. \hbox{Bingham et. al. \cite{bingham1987}.} For convenience of the reader the key properties are gathered in Lemma~\ref{SR} below. 

Due to the compound Poisson structure of $L$ we have for all $z\neq0$ and $t>0$ a representation of the density $\mu_t$ of $L_t$ by 
\be
\mu_t(z)=\sum_{m=1}^\infty\P(N_t=m)\,\frac{\nu^{*m}(dz)}{dz},\label{eq:7}
\ee
where $\nu^{*m}$ is the $m$-fold convolution of $\nu$ with itself, with $\nu^{*1} = \nu$. Using the full strength of the smoothly regular variation property of $f$ and its convexity we start  in Proposition~\ref{l:2Rn} by showing exponentially sharp upper and lower estimates of $\nu^{*m}(dz)/dz$, for all $|z|>k m$, $m\in \bN$ and some universal constant $k>0$. In the first main result of this article, Theorem~\ref{th:1Rn}, we apply the convolution estimates of $\nu^{*m}(dz)/dz$ with \eqref{eq:7}  and obtain exponentially tight density estimates. To formulate those density estimates, we use an auxiliary function $g$, defined as the solution of a nonlinear functional equation \eqref{Def_gRn}. For the choice $f(x)=x^\alpha$, $\alpha>1$, this function $g$ coincides with the function of the same name defined as the solution of equation \eqref{e:funktional} in Theorem \ref{th:resumen}.
 
A comment about the scope and quality of our exponential density estimates \eqref{mu.tRn}: Although L\'evy processes are frequently investigated objects, 
the literature concerning the  marginal distributions densities and their tail behavior remains fragmented. 
To illustrate the increasing precision of available estimates, 
we start with a dichotomy result in Sato~\cite[Theorem 26.1]{sato}. 
It states that for every L\'evy process in~$\bR^n$ with any L\'evy measure $\nu\neq0$  the following limit is satisfied: 
\be\label{e:Sato}
\lim_{\Lambda\rightarrow\infty}\frac{\ln\P(|L_t|>\Lambda)}{\Lambda\ln \Lambda}=-\inf\{1/c\mid c>0,\nu(\bR^n\setminus\{y\in\bR^n||y|\ge c\})>0\}.
\ee
In Imkeller, Pavlyukevich, Wetzel \cite{Imkeller}  and Wetzel \cite{Wetzel} the authors studied the first passage problem for SDE driven by  L\'evy processes. They applied the concept of regular varying functions to identify those (scalar) L\'evy processes within the scope of their results. In particular, they assume the tails of the jump measure to be $\nu([x, \infty)) = e^{-f(x)}$ for $f$ some regularly varying function of index $\alpha>1$. In such a setting the authors show  a more sophisticated tail estimate.

The proof of Theorem~2.2. in \cite{Imkeller} is based on estimates of the following type. 
There is a uniform bound $\tilde D_\e$, $\tilde D_\e \to 0$, such that for every fixed $\Lambda,t,h>0$ 
the parameter $\e$ can be chosen sufficiently small such that
\be
|\ln\P(\e L_t>\Lambda)+\Lambda \tilde D_\e| \spac{8pt}\le h \tilde D_\e.\label{eq:1.11}
\ee
The upper bounds is explicitly deduced in Lemma~5.1 in \cite{Imkeller}, see also Lemma \ref{levytail}, while the lower bound appears implicitly in Subsection 4.2 of  \cite{Imkeller}. 
The results of this article (Theorem~\ref{th:1Rn}(i)) improve \eqref{eq:1.11} in three main aspects:  
(1) Our results are valid for L\'evy processes with values in $\bR^n$.  
(2) Instead of only the tails for the L\'evy process $L_t$ 
of the distribution we estimate its density $\mu_t$, $t>0$. 
(3) The exponential estimate is asymptotically sharp. Moreover, \eqref{eq:1.11} can be recovered.

In order to achieve these improvements, we use slightly stricter conditions when selecting the jump measure. Firstly, we assume $\nu$ to be rotationally invariant in $\bR^n$. Secondly, we replace the concept of regular variation by the concept of smooth regular variation, in order to allow for higher order approximations.   That is to say, the Lévy measures considered in this article have an absolutely continuous component with density   $\exp(-f(|x|))$, $x\neq0$, for $f$ some smoothly regularly varying function of index $\alpha>1$.

More precisely, in \eqref{mu.tRn} it is shown that for any $h>0$ and any $-\infty<\rho<\gamma<1$ a parameter $\Lambda_o>0$ can be chosen sufficiently large, such that for $|\Lambda| >\Lambda_0$ and $t\in[\Lambda^\rho,\Lambda^\gamma]$ 
\be
\Big|\ln\mu_t(\Lambda)+ |\Lambda|\Big(f'(g(\frack {\Lambda}{t}))+g(\frack {\Lambda}{t})^{-1}\Big)\Big|\leq \delta g(\frack {\Lambda}t)^{-1},
\label{eq:1.12}
\ee
where $g$ is the solution of the functional equation \eqref{Def_gRn}, which in the case $f(x)=x^\alpha$ coincides with the solution of \refq{e:funktional} of Theorem \ref{th:resumen}. In addition, we show that 
\be
\lim_{\Lambda \rightarrow\infty}\sup_{t\in[\Lambda^\rho,\Lambda^\gamma]}\frac{\ln(\Lambda) g(\frack {\Lambda}{t})^{-1}}{f'(g(\frack {\Lambda}{t}))}\in(0,\infty), \qquad \mbox{thus}\qquad \lim_{\Lambda \rightarrow\infty}\sup_{t\in[\Lambda^\rho,\Lambda^\gamma]}\frac{g(\frack {\Lambda}{t})^{-1}}{f'(g(\frack {\Lambda}{t}))}= 0.
\label{eq:1.12b}
\ee
\noindent In Corollary~\ref{corRn} which is given for values in $\bR^n$ the density $\mu_t(\Lambda)$ is replaced by 
the tail $\P(|L_t| > \Lambda)$. For $n=1$, it is therefore natural to compare the estimates in \refq{eq:1.11} and \refq{eq:1.12}. 
 In order to do so, we set $\Lambda=\e^{-1}$. First, we see that the results \refq{eq:1.11} and \refq{eq:1.12} are consistent. 
 For this reason, we verify the limits
\be\label{eq:1119}
\lim_{\e\rightarrow0}\tilde D_\e\Big(f'\Big(g\Big(\frac{\e^{-1}}{t}\Big)\Big)+g\Big(\frac{\e^{-1}}{t}\Big)^{-1}\Big)^{-1}=\lim_{\e\rightarrow0}\tilde D_\e f'\Big(g\Big(\frac{\e^{-1}}{t}\Big)\Big)^{-1}=1
\ee
uniformly for $t\in[\e^{-\rho},\e^{-\gamma}]$. The first identity is a direct consequence of \refq{eq:1.12b}. The second identity is justified in Remark \ref{remRn} below, once the necessary properties of $g$ are available in Lemma~\ref{g_lnRn}. 

Furthermore, the limits in \refq{eq:1.12b} yield that the results of \refq{eq:1.12} provide a considerably finer quantification of the tails than the results in \refq{eq:1.11}. We stress that on the level of precision of~\eqref{eq:1.12} we manage to identify the impact of the time variable $t$, which remains hidden on the right-hand side of \eqref{eq:1.11}, and which is an important novelty of this paper. In particular, the identification of the leading exponential order $f'(g(\frack {\Lambda}{t}))$ and the second exponential order  $g(\frac{\Lambda}{t})^{-1}$ of $(\Lambda,t)$ jointly are the key results which open the gate for LDP results for the rescaled L\'evy (bridge) process. Later on, \refq{eq:1.12} is applied with the choice $\Lambda=|\x|\e^{-1}$ and $t=r_\e T$ for sufficiently small~$\e$. Again, for those fine density estimates we use the smooth regular variation of $f$. After the compound Poisson case we show such estimates to remain valid if we allow the existence of a Brownian component, a deterministic drift and infinite L\'evy measures. 

\textbf{II.} In a second step we condition the process $(X^\e_t)_{t\in[0,1]}$ on the event $\{X^\e_1=\x\}$ for some $\x\neq 0$. We denote the resulting L\'evy bridge by $Y^\e = (Y^\e_t)_{t\in[0,1]}$. 
Note that the event $\{X^\e_1=\x\}$ is typically a null set. However, since the law of $X^\e_1$ is absolutely continuous, we can calculate the densities $\bar \mu_t$ of $Y^\e_t$ in terms of the densities $\mu_t$ of $L_t$. It is well-known by \cite{FitzsimmonsPitmanYor, Schroedinger} that for $t\in (0,1)$, $y\in \bR^n$ and $\e>0$  we have 
\be
\bar\mu^\e_t(y)=\e^{-1}\mu_{r_\e t}(y\e^{-1})\mu_{r_\e(1-t)}((\x-y)\e^{-1})\mu_{r_\e}(\x \e^{-1})^{-1}.\label{eq:8}
\ee 
We capitalize on the asymptotics \eqref{eq:1.12} and \eqref{eq:1.12b} and obtain for small $\e$ 
\ba
\ln \bar\mu^\e_t(y)
&= f'(g(\frac{|y|\e^{-1}}{t r_\e})) |y| \e^{-1}  +f'(g(\frac{|\x - y|\e^{-1}}{(1-t) r_\e}))|\x - y| \e^{-1} -f'(g(\frac{|\x|\e^{-1}}{ r_\e})) |\x| \e^{-1}  \\
&\qquad + o\Big(g\big(\frac{\e^{-1}}{r_\e}\big)|\x -y| \e^{-1}\Big).
\ea
We consider the case $y \in [[0, \x]]$ and obtain $|\x|  = |\x -y| + |y|$, 
for which we may rewrite the right-hand side as 
\begin{align}
\ln \bar\mu^\e_t(y)
&=\bigg(f'(g(\frac{|y|\e^{-1}}{t r_\e}))-f'(g(\frac{|\x|\e^{-1}}{ r_\e}))\bigg) |y| \e^{-1}  +\bigg(f'(g(\frac{|\x - y|\e^{-1}}{(1-t) r_\e}))-f'(g(\frac{|\x|\e^{-1}}{ r_\e}))\bigg)|\x - y| \e^{-1}\nonumber\\
&\qquad + o\Big(g\big(\frac{\e^{-1}}{r_\e}\big)|\x -y| \e^{-1}\Big).\label{sdf}
\end{align}
By the specific choice of the function $g$ we  obtain the following \textit{asymptotic cancellation relation}, which is formulated and proven in Lemma~\ref{g_lnRn}.(ii),  in that  for any $\gamma>0$
\begin{align}
&f'(g(|u|\Lambda))-f'(g(\Lambda)) = \frac{\ln |u|}{g(\Lambda)} + o(\frac{|\ln |u||}{g(\Lambda)}) \nonumber\\
&\qquad \mbox{uniformly for }|u|\in [\ln(\Lambda)^{-\gamma} , \ln(\Lambda)^{\gamma}],\;\mbox{as } \Lambda \to \infty. \label{cancel}
\end{align}
Using \refq{cancel} for the parameters $\Lambda = (\e r_\e)^{-1}$, $u_1 = \frac{y}{t}$,  $u_2=\x$, $u_3 = \frac{\x-y}{1-t}$ and $u_4=\x$ in \refq{sdf} we obtain for $|y| \in [|\ln \e|^{-\gamma} , |\x| - |\ln \e |^{-\gamma}]$ and small values of $\e$ 
\begin{align}
-\ln \bar\mu^\e_t(y)
&\approx \bigg(\ln \frac{|y|}{t |\x|}\bigg) g((\e r_\e)^{-1})^{-1} |y| \e^{-1} 
+ \bigg(\ln \frac{|\x- y|}{(1-t) |\x|}\bigg) g((\e r_\e)^{-1})^{-1} |\x - y| \e^{-1} \nonumber\\
&= \bigg(|y| \ln \frac{|y|}{t} + |\x -y| \ln \frac{|\x -y|}{1-t} - |\x| \ln |\x| \bigg) g((\e r_\e)^{-1})^{-1} \e^{-1}. \label{eq:asymptoticcancellation}
\end{align}
Consequently, 
with speed function $S(\e) = \e g((\e r_\e)^{-1})$ we obtain the limit 
\begin{align}\label{eq:123}
\lim_{\e\to 0} -S(\e) \ln \bar\mu^\e_t(y) = |y| \ln \frac{|y|}{t} + |\x -y| \ln \frac{|\x -y|}{1-t} - |\x| \ln |\x|.
\end{align}
In the case of $y\notin [[0, \x]]$ we have $|\x|<|y|+|\x-y|$. Therefore the cancellation \refq{cancel}  in  the asymptotics of $ \ln \bar\mu^\e_t(y)$ in \refq{sdf} is incomplete and instead of \refq{eq:123} we obtain in this case 
\be\label{eq:124}
\lim_{\e\to 0} -S(\e) \ln \bar\mu^\e_t(y) = \infty.
\ee
Obviously this implies the LDP for $(Y^\e_t)_{\e>0}$ with speed function $S$ and a rate function given by the right-hand side of \refq{eq:123} and \refq{eq:124}. On the other hand, let $I_{\x}$ be the rate function from Theorem~\ref{th:resumen}.2, then the right hand side of \refq{eq:123} and \refq{eq:124} equals the infimum $\inf_{\varphi(t)=y}I_{\x}(\varphi)=I_{\x}(\varphi^*_{t,y})$, where $\varphi^*_{t,y}$ is the linear interpolation of the data points $\varphi^*_{t,y}(0)=0$, $\varphi^*_{t,y}(t)=y$ and $\varphi^*_{t,y}(1)=\x$.

Being a bridge process though, $(Y^\e)_{t\in[0,1]}$ does not have independent increments such that equation~\refq{eq:1} does not apply directly. Therefore in Proposition~\ref{l:43Rn} we directly state a multidimensional LDP for $(Y^\e_{t_1},\dots,Y^\e_{t_m})_{\e>0}$ and any set of times $0<t_1<\dots<t_m<1$. The proof of Proposition~\ref{l:43Rn} is obtained by the use of a $m$-fold $n$-dimensional version of \refq{eq:8}. The resulting rate function can be described by an adapted version of the right-hand sides of \refq{eq:123} and \refq{eq:124}.

\noindent \textbf{III.} Finally, we show exponential tightness of the family $\mbox{Law}(Y^\e)_{\e>0}$ with respect to the speed function $S$ established in Lemma \ref{l:44Rn}. As a consequence of our procedure we directly obtain the desired LDP for $(Y^\e)_{t\in[0,1]}$ on the path space. An application of equation \refq{eq:2} determines the integral shape of the associated rate function. 

The arguably the closest results in the existing literature so far are given in \cite{BazbhaBlanchetRheeZwart19} and \cite{Yang14}. We stress, that our results differ considerably from their findings. The results in \cite{BazbhaBlanchetRheeZwart19} establish a LDP for one dimensional and one sided rescaled pure jump L\'evy processes with Weibull jumps of exponent $\alpha\in (0,1)$ in the sample mean regime 1. of Subsection~\ref{sss:glance}, whereas we treat $\alpha >1$ and the ``short time'' regime 3. in Subsection~\ref{sss:glance}. 
 Their rate function is neatly given in terms of the $\alpha$-variation of the paths, which seems incomparable to the entropy we obtain. We refer to \cite{Imkeller} for a complete discussion of the different regimes which emerge for $\alpha\in(0,1)$ and $\alpha>1$.
 In~\cite{Yang14} the author establishes a LDP for speed function $S(\e) = \e$ for one dimensional rescaled pure jump L\'evy bridges in the sample mean regime of item 1 in Subsection~\ref{sss:glance} for bounded jumps (whereas our results treat Weibull distributed jumps for index $\alpha>1$). The rate function obtained is abstract and given as an integral of the path derivative in the Fenchel-Legendre transform of the characteristic exponent, which has no known closed form in general. Our results instead establish a LDP for multidimensional rescaled L\'evy bridges (not necessarily pure jump, and may include a Brownian component) with a nontrivial rotationally invariant jump component with Weibull jumps of exponent $\alpha>1$ in the short times regime $r_\e \ll \e^{-1}$ as explained in item 3 of Subsection \ref{sss:glance}. Our rate function $I_{\x}$ is given explicitly and shown to be a simple variant of the entropy on the path $[[0, {\x}]]$. Moreover, our results are robust under independent Brownian or L\'evy perturbations (satisfying for instance some Orey condition) as long as the asymptotics of the density tails remains untouched in the sense of Hypothesis II or the sufficient conditions given in Lemma \ref{hinr_bed2} below.

\bigskip 

\subsection{Organization of the article }
The article is organized as follows. After the basic definitions in Subsection~\ref{ss:formal} 
we formulate the main hypothesis in Subsection~\ref{sss:hypo} and state the first main result Theorem~\ref{th:1Rn}, which generalizes item 1. of Theorem~\ref{th:resumen} in Subsection~\ref{sss:density}. 
In Subsection~\ref{sss:bridgeLDP} we state the second main result Theorem~\ref{Th:2Rn}, 
which generalizes item 2. of Theorem~\ref{th:resumen}. 
In Subsection~\ref{sss:empirical} we generalize item 3. and 4. of Theorem~\ref{th:resumen} in Theorem~\ref{Th:3Rn}. Finally, we formulate in Subsection~\ref{ss:SC} 
sufficient conditions for the main Hypotheses~II, which is stated in Subsection~\ref{sss:hypo}. 

Most of the proofs are all gathered in Section~\ref{s:proofs}. We start with the exponential jump convolution power density estimates of $\nu^{*m}$ in Subsection~\ref{ss:convolution}. In Subsection~\ref{ss:density} these results are used in order to show first main result given by the density estimates of \eqref{eq:8} in Theorem~\ref{th:1Rn} in different settings (a), (b) and (c) of increasing generality. In Subsection~\ref{ss:LDP} we establish the LDP for the L\'evy-bridges stated in Theorem~\ref{Th:2Rn}. In Subsection~\ref{sss:N} we establish the LDP for the number of jumps given in Theorem~\ref{Th:3Rn}(i), while  Theorem~\ref{Th:3Rn}(ii) is shown in Subsection~\ref{sss:W}. 
In Appendix~\ref{a:ProofsSC} we give the proof of the sufficient conditions stated in Subsection~\ref{ss:SC}. 
In Appendix~\ref{s:LDPnonexistence} we show a negative result for the exponential tightness of rescaled L\'evy processes on path space which excludes a classical LDP for $r_\e\ll\e^{-1}$, as announced above.

\section{Object of study and main results}\label{s:OoSaMR}
\setcounter{equation}{0}

\subsection{The fundamental concepts and the basic notation}\label{ss:formal}\setcounter{equation}{0}

\subsubsection{Smoothly regularly varying functions} 

An important tool in our analysis turns out to be the concept of regular varying and smoothly regularly varying functions. The definition and many properties of both can be found in Chapters~1.4 and 1.8,  respectively, of Bingham \cite{bingham1987}. In contrast to the original definition we allow those functions to have negative values at some starting intervals.	

\begin{defin}\label{def:SR}
\begin{enumerate}
 \item[(i)] Let $z\in\bR$. A function $f:(z,\infty)\rightarrow\bR$ is called \textbf{regular varying with index $\alpha\in\bR$}, if $\sup\{\Lambda\ge z\mid f(\Lambda)\le0\}<\infty$ and $\lim_{\Lambda\rightarrow\infty}\frac{f(\lambda \Lambda)}{f(\Lambda)}=\lambda^\alpha$ holds for every $\lambda>0$. We denote by $R_\alpha$ the class of regular varying functions with index $\alpha$.
 \item[(ii)] Let $z\in\bR$. A function $f:(z,\infty)\rightarrow\bR$ is called \textbf{smoothly regularly varying with index $\alpha\in\bR$}, if $f$ is infinitely often differentiable, $z_o:=\max\{1,\sup\{\Lambda\ge z\mid f(\Lambda)\le0\}\}<\infty$, and $h:(\ln z_o,\infty)\rightarrow\bR$, $h(\cdot):=\ln f(\exp(\cdot))$ satisfies $\lim_{\Lambda\rightarrow\infty}h'(\Lambda)=\alpha$ and $\lim_{\Lambda\rightarrow\infty}h^{(m)}(\Lambda)=0$ for any $m=2,3,\dots$. We denote by $SR_\alpha$ the class of smoothly regularly varying functions with index $\alpha$.
 \item[(iii)] Let $z>0$. A function $f:(0,z)\rightarrow\bR$ is called \textbf{(smoothly) regularly varying in $0$ with index $\alpha\in\bR$}, if $(\Lambda\mapsto f(\Lambda^{-1})^{-1})$ is (smoothly) regular varying in $0$ with index $\alpha$.
\end{enumerate}
\end{defin}
Usually there is no risk of confusion between (smooth) regular variation in 0 and (smooth) regular variation in $\infty$ and this distinction is omitted. For convenience of the reader we gather the most important properties.

\begin{lemma}[Properties of smoothly regularly varying functions]\label{SR}
Let $\alpha,\beta\in\bR$,  $f\in SR_\alpha$ and $g\in SR_\beta$. Then the following statements are valid:
$$\begin{array}{l}
\begin{array}{ll}
\hspace{-12pt}(i)\;f\in R_\alpha&(ii)\,\lim_{\Lambda\rightarrow\infty}\frac{\Lambda f'(\Lambda)}{f(\Lambda)}=\alpha\\
\hspace{-12pt}(iii)\mbox{ If }\alpha\ge\beta\mbox{ then }f+g\in SR_{\alpha}.\hspace{18mm}\left.\right. &(iv)\mbox{ If }\alpha>\beta\mbox{ then }f-g\in SR_{\alpha}.\\
\hspace{-12pt}(v)\;f\cdot g\in SR_{\alpha+\beta}&(vi)\mbox{ If }\lim_{\Lambda\rightarrow\infty} g(\Lambda)=\infty\mbox{  then }f\circ g\;\in SR_{\alpha\beta}.\\
\hspace{-12pt}(vii)\;1/f\in SR_{-\alpha}&(viii)\mbox{ If }\alpha\neq0\mbox{ then }|f'|\in SR_{\alpha-1}.
\end{array}\\
\hspace{-7pt}(ix)\mbox{ If }\alpha>-1,\mbox{ then $z$ can be chosen sufficiently large, such that $\Lambda\mapsto\int_z^\Lambda f(y)dy$ exists.}\\
\mbox{This function belongs to $SR_{\alpha+1}$}.\\
\hspace{-7pt}\mbox{$(x)$ If $\alpha>0$, then $z$ can be chosen sufficiently large, such that $f$ is invertible on $[z,\infty)$ and its}\\
\mbox{inverse function $f^{-1}$ belongs to $SR_{\frac1\alpha}$}.\\
\end{array}$$
\end{lemma}
The proof of Lemma~\ref{SR} can be found in Chapter~1.8 of \cite{bingham1987}.

\begin{rem}\label{alpha0}
In the case $\alpha=0$, part $(viii)$ of the preceding lemma does not allow any statement about the derivative of $f$ in general. However from part $(ii)$ it follows, that $|f'|$ in this case can be estimated from above by an $SR_{-1}$ function: For any $\kappa>0$ we obtain $|f'(x)|<\kappa x^{-1}f(x)$ for $x$ sufficiently large, where $(x\mapsto x^{-1}f(x))\in SR_{-1}$, which is clear by part $(iii)$. More specifically, we consider $\beta\neq0$ and $f(x)=\ln g(x)$, $g\in SR_\beta$. Then $f\in SR_0$ and parts $(v)$, $(vii)$ and $(viii)$ imply that $|f'(x)|=\Big|\frac{g'(x)}{g(x)}\Big|\in SR_{-1}$.
\end{rem}

\bigskip 

\subsubsection{L\'evy processes with values in $\bR^n$}\label{sss:Levy}

Given a complete filtered probability space $(\Omega, \cA, \P, (\cF_t)_{t\geq 0})$ 
we consider a L\'evy process $L = (L_t)_{t\geq 0}$, see \cite[Definition 1.6]{sato}. 
By the L\'evy-Khintchine formula \cite[Theorem 8.1]{sato} the characteristic function of the marginals $L_t$, $t>0$, of $L$ and thus its distribution can be described by a generating triplet $(\sigma^2,\nu,\Gamma)$. For further details we refer to \cite{applebaum}. Here $\sigma^2\in \bR^{n\times n}$ is a symmetric and non-negative definite square matrix. 
Moreover, $\nu: \cB(\bR^n) \to [0, \infty]$ is a sigma finite measure (the so-called L\'evy measure) satisfying 
\[\nu(\{0\})=0\mb{10mm}{and}\int_{\bR^n} (1\wedge |y|^2)\nu(dy)<\infty\]
and $\Gamma\in \bR^n$ represents a deterministic drift. 
For any $t\geq 0$ the characteristic function has the following representation
\ba
\lambda \mapsto &\E\exp(i\langle \lambda, L_t\rangle)\\
&=\exp\Big(t\Big(-\frac{1}{2} \langle \lambda, \sigma^2 \lambda\rangle + i\langle \lambda, \Gamma\rangle
+\int_{\bR^n}(\exp(i\langle \lambda, y\rangle )-1-(i\langle \lambda, y\rangle)\textbf{1}_{\{|y|\leq 1\}}(y))\nu(dy)\Big)\Big).
\ea
If $\nu$ is finite, the characteristic function can be written as 
\be
\E\exp(i\langle \lambda, L_t\rangle)=\exp\Big(t\Big(-\frac{1}{2} \langle \lambda, \sigma^2 \lambda\rangle + i\langle \lambda,\tilde \Gamma\rangle+\int_\bR(\exp(i\langle \lambda, y\rangle)-1)\nu(dy)\Big)\Big),\label{eq:14}
\ee
where 
\be\label{tildegamma}
\tilde \Gamma=\Gamma+\int_{\{|y|\le 1\}} y \nu(dy).
\ee 
Whenever in this paper a L\'evy measure is known to be finite, we tacitly agree to refer to representation \refq{eq:14}. In particular, if a L\'evy process exhibits a generating triplet $(0,\nu,0)$ with $\nu$ being finite, then we agree that this L\'evy process shall be equal to the sum of its jumps.

Let $\eta$ and $\xi$ be two independent L\'evy processes with generating triplets $(\sigma^2_\eta,\nu_\eta,\Gamma_\eta)$ and $(\sigma^2_\xi,\nu_\xi,\Gamma_\xi)$. Then $\eta+\xi$ is known to be a L\'evy process with generating triplet $(\sigma^2_\eta+\sigma^2_\xi,\nu_\eta+\nu_\xi,\Gamma_\eta+\Gamma_\xi)$. 

When we study a L\'evy process $L$ with generating triplet $(\sigma^2,\nu,\Gamma)$ we may define by the L\'evy-It\^o decomposition \cite[Theorem 19.2+19.3]{sato} two independent L\'evy processes $\eta$ and $\xi$, such that $\eta+\xi$ is identically distributed as $L$. Within this paper the L\'evy measure $\nu_\eta$ are defined, such that $\nu_\eta(\bR^n)<\infty$. In order to use the representation \refq{eq:14} for the characteristic function of $\eta_t$, we choose~$\Gamma_\eta$ and $\Gamma_\xi$ such that
$\Gamma_\xi+\Gamma_\eta-\int_{\{|y|\le 1\}} y \nu_\eta(dy)=\Gamma$.

\subsubsection{Large deviations principles}\label{sss:LDPs}
As mentioned above, we need to extend the classical concept of a large deviation principle by a more general speed function. For this purpose, the concept of exponential tightness needs to be adapted, too.
\begin{defin}\label{def:LDP}
\begin{enumerate}
 \item[(i)] A function $S:(0,\infty)\rightarrow(0,\infty)$ is called a \textbf{speed function}, if $\lim\limits_{\e\rightarrow0+} S(\e)=0$ \mbox{and there exists a continuous invertible function $S_o:(0,\infty)\rightarrow(0,\infty)$, such that $\lim\limits_{\e\rightarrow0+}\frac{S_o(\e)}{S(\e)}=1$.}
 \item[(ii)] Let $(\mathfrak{X},\mathcal{T})$ be a topological space equipped with its Borel-$\sigma$-algebra $\mathcal{B}$, and $(X^\e)_{\e>0}$ be a family of random elements with values in $\mathfrak{X}$. $\mbox{Law}(X^\e)_{\e>0}$ is said to satisfy a \textbf{large deviations principle (LDP) on $(\mathfrak{X},\mathcal{T})$ with respect to a speed function $S$ and a rate function~$I$}, if for every open subset $A\subset\mathfrak{X}$
\be
\liminf_{\e\rightarrow0}S(\e)\ln\P(X^\e\in A)\ge-\inf_{x\in A}I(x)
\ee 
is valid and for every closed subset $A\subset\mathfrak{X}$ 
\be
\limsup_{\e\rightarrow0}S(\e)\ln\P(X^\e\in A)\le-\inf_{x\in A}I(x).
\ee 
\item[(iii)] Let $(\mathfrak{X},\mathcal{T})$ be a topological space equipped with its Borel-$\sigma$-algebra $\mathcal{B}$, let $(X^\e)_{\e>0}$ be a family of random elements with values in $\mathfrak{X}$ and $S$ a speed function. $(X^\e)_{\e>0}$ is said to be \mbox{\textbf{$S$-exponentially tight},} if for every $k>0$ there exists a compact subset $A_k\subset \mathfrak{X}$, such that
\be
\limsup_{\e\rightarrow0}S(\e)\ln\P(X^\e\notin A_k)\le-k.
\ee 
\end{enumerate}
\end{defin}

\begin{rem}
Let $S$ be a speed function. By definition there is a continuous invertible function $S_o:(0,\infty)\rightarrow(0,\infty)$, such that $\lim\limits_{\e\rightarrow0}\frac{S_o(\e)}{S(\e)}=1$. It is important to notice, that the following two statements are equivalent:
\begin{itemize}
	\item A family $(X^\e)_{\e>0}$ of random elements satisfies a LDP on a topological space $(\mathfrak{X}, \mathcal{T})$ equipped with its Borel-$\sigma$-algebra $\mathcal{B}$, 
	with respect to the speed function $S$ and a rate function $I$. 
	\item The family $(Y^\e)_{\e>0}$ of random elements with $Y^\e:=X^{S_o^{-1}(\e)}$ satisfies a LDP with the usual speed function $\e\mapsto\e$ on a topological space $(\mathfrak{X}, \mathcal{T})$ with respect a rate function $I$.
\end{itemize}
It is this equivalence which allows us to use all those results, the contraction principle and many more, which have been developed inside the common concept of large deviations. Definition \ref{def:LDP} can be understood as a continuous version of definition 1.12 in \cite{OlivieriVares05}.
\end{rem}

\bigskip 
\subsection{The hypotheses and the main results} 

\subsubsection{The hypotheses on the L\'evy process $L$}\label{sss:hypo}

Throughout the present paper we consider $\e\mapsto r_\e$ to be a regular varying function with its index in $(-1,\infty)$. The process $L$ denotes a L\'evy processes with values on $\bR^n$ in the sense of Subsection~\ref{sss:Levy}. The generating triplet of $L$ is denoted by  $(\sigma^2,\nu,\Gamma)$. \vspace{5mm}

\noindent\textbf{Hypotheses on the L\'evy process $L$:} 
The elements of the generating triplet $(\sigma^2,\nu,\Gamma)$ of $L$ satisfy the following properties. 
The L\'evy measure $\nu$ can be written as $\nu=\nu_\eta+\nu_\xi$ and satisfies: 
\begin{enumerate}
  \renewcommand{\labelenumi}{\Roman{enumi}}\label{BedingungenRn}
	\item  The L\'evy measure $\nu_\eta$ is finite, $\nu_\eta(\bR^n)<\infty$, and has a density on $\bR^n\setminus\{0\}$ of the form 
	\be\label{27}
	\nu_\eta(dz)/dz=\exp(-f(|z|)),
	\ee
	where~$f$ is a smoothly regularly varying function $f\in SR_\alpha$ for some $\alpha>1$. 
	\item Let $\xi$ denote a L\'evy process with generating triplet $(\sigma^2,\nu_\xi,\Gamma_\xi)$ with $\Gamma_\xi=\Gamma+\int_{\{|y|\le 1\}} y \nu_\eta(dy)$. There is a $\tilde s\in\bR$, such that~$\xi_t$ has a density $\mu_{\xi,t}$ on $\bR^n\setminus[-\tilde s t,\tilde s t]^n$ for every $t>0$. Furthermore a parameter $\aleph>1-\frac1\alpha$ exists, such that for every $\gamma<1$ the following limit is valid 
	\be
	\lim_{\Lambda\rightarrow\infty}\sup_{\substack{t<\Lambda^\gamma\\|y|= \Lambda}}\frac{\ln\mu_{\xi,t}(y)}{\Lambda(\ln \Lambda)^\aleph}=-\infty.\label{BEDIIRn}
	\ee
\end{enumerate}

\begin{rem}\label{rem:cpp}
\begin{enumerate}
 \item Hypothesis~II is clearly satisfied if $\nu_\xi=0$ and $\mbox{rank}(\sigma^2) \in \{0, n\}$. 
 In this case the generating triplet of $L$ equals $(\sigma^2,\nu_\eta,\Gamma)$. Therefore those conditions can be interpreted in the following way: The asymptotic is defined by $\nu_\eta$ which obeys Hypothesis~I. Hypothesis~II defines the limitations within which an additional jump activity does not disturb the asymptotics of~I: On the one hand, the L\'evy measure of this additional jump activity has to have lighter tails than $\nu_\eta$. This is encoded in \eqref{BEDIIRn}. 
 On the other hand, $\nu_\xi$ needs not satisfy smoothness criteria and symmetries as strict as $\nu_\eta$ given by \refq{27} with $f\in SR_\alpha$. We refer to Subsection~\ref{ss:SC} for sufficient conditions, which are easier to verify.

\item In many parts of our analysis we will derive the distribution of $L_t$ as a convolution of the distributions of $\eta_t$ and $\xi_t$. Under the assumption that $\eta_t$ and $\xi_t$ exhibit a density $\mu_{\eta,t}$ and $\mu_{\xi,t}$, respectively, with respect to the Lebesgue measure on $\bR^n$,  we can describe the density of $L_t$ as

\be
\mu_{L,t}(x)=\int_{\bR^n}\mu_{\eta,t}(x-y)\mu_{\xi,t}(y)dy.\label{muL}
\ee

By Hypothesis~I it follows that $\eta_t$ has a density on $\bR^n\setminus\{0\}$ and by Hypothesis~II it follows that $\xi_t$ possesses a density on $\bR^n\setminus[-\tilde st,\tilde st]^n$. For any $\gamma<1$ and $|x|$ sufficiently large it follows that $L_t$ possesses a density on the position $x$ for any $t\le|x|^\gamma$. Throughout this paper we agree to understand \refq{muL} as a simplified way of writing the mathematically precise representation

$$\mu_{L,t}(x)=\frac{d}{dx}\int_{\bR^n}\upsilon_{\eta,t}(dx-y)\upsilon_{\xi,t}(dy),$$

where $\upsilon_{\eta,t}$ and $\upsilon_{\xi,t}$ denote the law of the distributions of $\eta_t$, $\xi_t$ respectively.

 \item Sufficient conditions for Hypothesis~II are stated in Subsection~\ref{ss:SC} and proved in Appendix~\ref{a:ProofsSC}. 
\end{enumerate}
\end{rem}

\subsubsection{First main result: exponential density estimates of $\mbox{Law}(L_t)$ (Theorem~\ref{th:1Rn})}\label{sss:density}

Before we present the desired density estimate, we state an estimate for the $m$-th convolution power~$\nu^{*m}$, $m\in \bN$, motivated by Subsection~\ref{sss:approach}.I. 

\begin{prop}[$m$-th convolution density tail estimate]\label{l:2Rn} 
Let Hypothesis~I be satisfied for $\nu = \nu_\eta$ with $\nu(\bR^n)=1$. For $m\in\bN$ we denote the $m$-fold convolution of $\nu$ with itself by~$\nu^{*m}$. Then for any $\delta>0$ there is a $k>0$ such that for all $m\in\bN$ and all $|x|>km$ it follows  
\ba
\qquad \frac{\nu^{*m}(dx)}{dx} &\le (\alpha-1)^{\frac{(n-1)(m-1)}{2}}(2\pi f''(\frack{|x|}{m})^{-1})^{\frac{n(m-1)}2}\cdot \exp(-m(f(\frack{|x|}{m})-\delta)),\\
\frac{\nu^{*m}(dx)}{dx} &\ge \frac{(\alpha-1)^{\frac{(n-1)(m-1)}{2}}}{m^{\frac n2}}(2\pi f''(\frack{|x|}{m})^{-1})^{\frac{n(m-1)}2}\cdot \exp(-m(f(\frack{|x|}{m})+\delta)).\\
\label{Faltung_FormelRn}
\ea
\end{prop}

\noindent 
The proof Proposition~\ref{l:2Rn} is given in Subsection~\ref{ss:convolution}. Estimate \refq{Faltung_FormelRn} applied to equation \refq{eq:7} leads then to the desired estimate of the density of $L_t$. 
\begin{theorem}[Exponential density estimate of $\mbox{Law}(L_t)$ by an auxiliary functional equation] \label{th:1Rn}
Let $L$ be a L\'evy process which satisfies Hypotheses~I and II and denote by $\mu_t$, $t>0$, the density of the marginal $L_t$. 
\begin{enumerate}
 \item[(i)] 
Then for every $\delta>0$ and every $\rho<\gamma<1$, there is some $k>0$, such that for every $|x|>k$ and every $t\in[|x|^\rho,|x|^\gamma]$ it follows 
\be
\hspace{-2mm}-|x|\Big(f'(g(\frack {|x|}t))-(1-\delta)g(\frack {|x|}t)^{-1}\Big)\le\ln\mu_t(x)\le\sup\limits_{s\le t}\ln\mu_s(x)\le-|x|\Big(f'(g(\frack {|x|}t))-(1+\delta)g(\frack {|x|}t)^{-1}\Big),\label{mu.tRn}
\ee
where $g:[\Lambda_0,\infty)\rightarrow\bR$ is the unique solution of the nonlinear functional equation 
\be
g(\Lambda)f'(g(\Lambda))-f(g(\Lambda))+\ln g(\Lambda)-\frack n2\ln f''(g(\Lambda))\spac{7pt}=\ln((2\pi)^{-\frac n2}(\alpha-1)^{-\frac{n-1}2} \Lambda)\label{Def_gRn}
\ee 
for some suitably chosen $\Lambda_0>0$. 

\item[(ii)] Let $\rho<\gamma<1$ and $c\in(\alpha^{-1},1)$. For every $\delta>0$ the value $|x|$ can be chosen sufficiently large, such that for every $t\in[|x|^\rho,|x|^\gamma]$ we have
\begin{equation}\label{e:cjspe}
\P(\sup_{s\le t}|L_s-L_{s-}|>|x|^c~|~L_t=x) \leq \exp(-|x|^{\alpha c-\delta}).
\end{equation}
\end{enumerate}
\end{theorem}
The proof is given in Subsection~\ref{ss:density} and Lemma~\ref{g_lnRn}(i). 
\noindent Obviously the density estimates of Theorem~\ref{th:1Rn}$(i)$ directly imply the following respective tail estimate. 

\begin{cor}[Exponential tail estimate of $\mbox{Law}(L_t)$]\label{corRn} 
Let the assumptions of Theorem~\ref{th:1Rn}$(i)$ be satisfied. 
Then for every $\delta>0$ and every $\rho<\gamma<1$, $k$ can be chosen sufficiently large, such that for every $\Lambda>k$ and every $t\in[\Lambda^\rho,\Lambda^\gamma]$ it follows 
\begin{align*}
&-\Lambda\Big(f'(g(\frack {\Lambda}t))+(1+\delta)g(\frack {\Lambda}t)^{-1}\Big)\\
&\qquad \le\ln \P(|L_t|>\Lambda)\le\ln \P\Big(\sup\limits_{s\le t}|L_s|>\Lambda\Big)
\le\ln \P\Big(\sup\limits_{0\le s_o\le s_1\le t}|L_{s_o}-L_{s_1}|>\Lambda\Big)\\
&\qquad \qquad \le-\Lambda\Big(f'(g(\frack {\Lambda}t))+(1-\delta)g(\frack {\Lambda}t)^{-1}\Big).
\end{align*}
\end{cor}

\noindent 	
The existence, uniqueness and the properties of the solution function $g$ of \eqref{Def_gRn} are crucial for our results. For a better understanding we give an asymptotic approximation of $g(\Lambda)$ in terms of the Fenchel-Legendre transformation. In order to assess the result of Theorem~\ref{th:1Rn} we then summarize the key properties of $g$ in Lemma \ref{g_lnRn} below.
\begin{rem}
Note that in the classical Cramér Theorem \cite[Section 2.2]{Dembo} the rate function is stated in terms of the Fenchel-Legendre transform. To give a better understanding of the definition of the function $g$ and the main term $f'(g(\cdot))$ in~\refq{mu.tRn} we sketch a connection to the Fenchel-Legendre transform of the logarithmic jump density $f$.

Consider $f$ to be given in Hypothesis I and for the sake of argument we assume $f$ to be convex. Denote by $f^*$ its Fenchel-Legendre transform $f^*(u):=\sup_x(ux-f(x))=x_uu-f(x_u)$  for an optimizer $x_u$. Hence we have $f^*(u)=x_uf'(x_u)-f(x_u)$. The optimizer $x_u$ satisfies $x_u=(f')^{-1}(u)$, thus $f^*(f'(u))=uf'(u)-f(u)$. A comparison with the defining equation of $g$ in \refq{Def_gRn} and the identification of its main terms yields the asymptotics
$$f^*(f'(g(\Lambda)))\spac{7pt}=g(\Lambda)f'(g(\Lambda))-f(g(\Lambda))\spac{7pt}=\ln\Lambda+O(\ln g(\Lambda)).$$
Thus for large values of $\Lambda$ we have 
$$f'(g(\Lambda))\spac{7pt}\approx(f^*)^{-1}(\ln\Lambda)\mb{4mm}{and}g(\Lambda)\spac{7pt}\approx (f')^{-1}\circ(f^*)^{-1}(\ln\Lambda).$$
\end{rem}\bigskip

\begin{lemma}[Key properties of $g$]\label{g_lnRn}
Let $\alpha>1$ and $f\in SR_\alpha$. Let $b<\alpha$ and $k\in SR_{b}$.
\begin{enumerate}
\item[$(i)$] \textbf{Existence, uniqueness and regularity: } There is some $r_o>0$ such that for every $\Lambda>r_o$ the equation
\begin{equation}
g(\Lambda)f'(g(\Lambda))-f(g(\Lambda))+k(g(\Lambda))\spac{8pt}=\ln \Lambda\label{e:NLFE}
\end{equation}
has a unique solution $g\in SR_o$. Let $h:(\ln r_o,\infty)\rightarrow\bR$ be defined by $h(\Lambda)=g(\exp(\Lambda))$. Then~$h\in SR_{\frac1\alpha}$. 
\item[$(ii)$] \textbf{The asymptotic cancellation relation: } For each $\gamma>0,$ $\delta>0$ there is $z>r_o$ such that for every $\Lambda>z$ and $y\in[(\ln \Lambda)^{-\gamma},(\ln \Lambda)^\gamma]$ the following estimate is valid
\be
\Big|g(\Lambda)[f'(g(y\Lambda))-f'(g(\Lambda))]-\ln y\Big|\le\delta|\ln y|.\label{eq.ln.gRn}
\ee
\item[$(iii)$] \textbf{Asymptotic behavior: } The function $g$ satisfies the following limits: 
\begin{align}
&\lim_{\Lambda\rightarrow\infty}\frack{\displaystyle g'(\Lambda)\Lambda\ln \Lambda}{\displaystyle g(\Lambda)}=\frack1\alpha,\label{e:lim4}\\
&\lim_{\Lambda\rightarrow\infty}f(g(\Lambda))(\ln \Lambda)^{-1}=\frack{1}{\alpha-1},\label{e:lim3}\\
&\lim_{\Lambda\rightarrow\infty}g(\Lambda)f'(g(\Lambda))(\ln \Lambda)^{-1}=\frack{\alpha}{\alpha-1},\label{e:lim2}\\
&\lim_{\Lambda\rightarrow\infty}(g(\Lambda)f'(g(\Lambda))-f(g(\Lambda)))(\ln \Lambda)^{-1}=1.\label{e:lim1}
\end{align}
\item[$(iv)$] \textbf{Sensitivity: } Given $k_0,k_1\in SR_b$ such that $c:=\lim_{\Lambda\rightarrow\infty}(k_1(\Lambda)-k_o(\Lambda))\in\bR$. 
In addition let $g_i$ be uniquely defined by $$g_i(\Lambda)f'(g_i(\Lambda))-f(g_i(\Lambda))+k_i(g(\Lambda))=\ln \lambda \qquad \mbox{ for }i\in\{0,1\}.$$ 
Then we have for $i\in \{0,1\}$ 
$$\lim_{\Lambda\rightarrow\infty}g_i(\Lambda)(f'(g_o(\Lambda))-f'(g_1(\Lambda)) =c.$$
\end{enumerate}
\end{lemma}
The proof is given at the beginning of Subsection~\ref{ss:density}. 

\begin{rem}\label{remRn}
\begin{enumerate}
 \item Lemma~\ref{g_lnRn} helps to understand the terms $g(\cdot)$ and $f'(g(\cdot))$ in estimate \refq{mu.tRn}. 
Part $(i)$ of the lemma and $f\in SR_\alpha$ lead to the following estimate: For each $\delta<0$ there is a $z>r_o$ such that for all $\Lambda>z$
\be
(\ln \Lambda)^{\frac1\alpha-\delta}<g(\Lambda)<(\ln \Lambda)^{\frac1\alpha+\delta}\mb{5mm}{and}(\ln \Lambda)^{\frac{\alpha-1}\alpha-\delta}<f'(g(\Lambda))<(\ln \Lambda)^{\frac{\alpha-1}\alpha+\delta}.
\ee 

 \item We stress the following important application of \eqref{e:lim2} in Lemma~\ref{g_lnRn}(iii), which is used in several parts of the proof. 
 Let $z_\Lambda>0$ such that $\lim_{\Lambda\rightarrow\infty}\frac{z_\Lambda\ln \Lambda}{\Lambda}=0$. 
 Then $\Lambda$ can be chosen sufficiently large, such that upper and lower bounds of $\mu_t(y)$ in \refq{mu.tRn} hold for any $\delta>0$ uniformly for all $y\in[\Lambda-z_\Lambda,\Lambda+z_\Lambda]$. The same is valid for the estimates of Corollary~\ref{corRn}.

\item\label{D_eps} Furthermore, Lemma~\ref{g_lnRn} puts us in the situation to prove the limit \refq{eq:1119} and thus the consistence of the results in Theorem \ref{th:1Rn} with the tail estimates used in \cite{Imkeller}. We start with the definition of $\tilde D_\e$. Let $q_{\e}:=\sup\{y\in(0,\infty)|f(y)<\ln |\e|\}$, $d_\alpha:=\alpha(\alpha-1 )^{-(1-\frac{1}{\alpha})}$ and $\tilde D_{\e}:=d_\alpha\frac{\ln |\e|}{q_{\e}}\e^{-1}$.

By item $(i)$ we have $\lim\limits_{\e\rightarrow0}f'(g(\frac{\e^{-1}}{t}))/f'(g(\e^{-1}))=1$ uniformly for $t\in[\e^{-\rho},\e^{-\gamma}]$. The definition of $q_\e$ yields $\lim\limits_{\e\rightarrow0}f(q_\e)/|\ln \e| =1$. Item $(iii)$ of the lemma then yields $\lim\limits_{\e\rightarrow0}f(g(\e^{-1}))/f(q_\e)=(\alpha-1)^{-1}$. By $f\in SR_\alpha$ we get  $\lim\limits_{\e\rightarrow0}\frac{g(\e^{-1})}{q_\e}=(\alpha-1)^{-\frac1\alpha}$. Combining these findings we get
\begin{align*}
\lim_{\e\rightarrow0}\frac{f'(g(\frac{\e^{-1}}{t}))\e^{-1}}{\tilde D_\e}&=\lim_{\e\rightarrow0}\frac{g(\e^{-1})f'(g(\e^{-1}))\e^{-1}}{g(\e^{-1})\tilde D_\e}\\
&=(\alpha-1)^{\frac1\alpha}\lim_{\e\rightarrow0}\frac{g(\e^{-1})f'(g(\e^{-1}))\e^{-1}}{q_\e\tilde D_\e}\\
&=\frac{(\alpha-1)^{\frac1\alpha}}{d_\alpha}\lim_{\e\rightarrow0}\frac{g(\e^{-1})f'(g(\e^{-1}))\e^{-1}}{\e^{-1}|\ln \e|}~=1.
\end{align*}
Here we use the definition of $d_\alpha$ and the limit \refq{e:lim2} in part $(iii)$ of the lemma for the last step.
\end{enumerate}
\end{rem}

\bigskip
\subsubsection{Second main result: a sample path LDP of the bridge $Y^\e$ (Theorem~\ref{Th:2Rn})}\label{sss:bridgeLDP}

 We state the following standing assumptions: Define $(X^\e)_{t\in[0,T]}$ by $X^\e_t:=\e L_{r_\e t}$ for a  L\'evy process~$L$ satisfying Hypotheses~I and II and  $\e\mapsto r_\e$ being of regular variation with a index in $(-1,\infty)$. Thus parameters $\rho<\gamma<1$ can be chosen, such that for every $\x\in\bR^n\setminus\{0\}$, and $\e$ sufficiently small we have $r_\e T\in[(|\x|\e^{-1})^\rho,(|\x|\e^{-1})^\gamma]$. In this sense we are set to apply Theorem~\ref{th:1Rn} and estimate the density of $X^\e_t$ for any $t\in(0,T]$. 
 
 \noindent Let $\x\in\bR^n\setminus\{0\}$ and $T>0$. For any $\e>0$ sufficiently small we condition the process $(X^\e)_{t\in[0,T]}$ on the event $\{X^\e_T=\x\}$ and denote the resulting bridge process by $(Y^{\x, \e}_t)_{t\in[0,T]}$. We denote by $\mathcal D_{[0,T],\bR^n}$ the space of c\`adl\`ag functions $[0,T]\rightarrow\bR^n$. In Theorem~\ref{Th:2Rn} below we establish the LDP for $(Y^{\x, \e})_{\e>0}$ on $\mathcal D_{[0,T],\bR^n}$ equipped with the uniform norm $||\,\cdot\,||_\infty$. Denote $[[0,\x]]=\{s\x\mid s\in[0,1]\}$ and set 
\ba
D_{\x,T}:=\{\varphi\in\mathcal D_{[0,T],\bR^n}\mid |\varphi(\cdot)|\mbox{ is continuous and nondecreasing with }\,\varphi([0,T])=[[0,\x]]\}.
\ea
Note that by the monotonicity assumption we have that $|\varphi(\cdot)|$ is differentiable almost everywhere on $[0,T]$ for $\varphi\in D_{\x,T}$.

\begin{theorem}[The LDP with speed function $S$ and rate function $I_{\x}$ for $Y^{\x, \e}$]\label{Th:2Rn} \hfill\\
Fix the notations and assumptions of this subsection. Fix $\x\in \bR^n$, $\x\neq 0$. Then the family $(\P^{\x, \e})_{\e>0}$, $\P^{\x, \e} = \Law(Y^{\x, \e})$ satisfies a LDP on $(\mathcal D_{[0,T],\bR^n},||\,\cdot\,||_\infty)$ with speed function $S(\e):=\e\cdot g(\e^{-1}r_\e^{-1})$, where $g$ is defined by \refq{Def_gRn} and the rate function
\begin{equation}\label{entropyrate}
I_{\x}(\varphi)=\left\{
\begin{array}{ll}
\displaystyle\int_0^T|\varphi|'~\ln|\varphi|' ~dt-|\x|\ln\frack{|\x|}T,&\mbox{if 
$\varphi \in D_{\x,T}$},\\[5mm]
\infty,&\mbox{otherwise}.
\end{array}
\right.
\end{equation}
Here we denote by $|\varphi|(t) = |\varphi(t)|$, $|\varphi|'(t) = \frac{d}{dt} |\varphi(t)|$ the total derivative, whenever it exists and set it equal to $0$ otherwise. We set $r \ln r = 0$, whenever $r = 0$.
\end{theorem}
\begin{rem}
The rate function in Theorem~\ref{th:1Rn} shows an interesting connection to the well-known Sanov theorem (for example Theorem~6.2.10 in \cite{Dembo} by Dembo and Zeitouni). To see this connection we consider the one dimensional case with the choice $\x=1$. In that case the rate function $I$ of Theorem~\ref{Th:2Rn} can be written
\be
I(\varphi)=\int_o^T\varphi'(t)\ln\varphi'(t)dt+\ln T
\ee
for every continuous non decreasing function $\varphi:[0,T]\rightarrow\bR$ with $\varphi(0)=0$ and $\varphi(T)=1$. Every such function can be interpreted as the cumulative distribution function of a probability measure~$\mu_\varphi$. Let $\mu_o$ denote the probability measure that is given by a uniform distribution on $[0,T]$ and by $H(\,\cdot\,|\,\cdot\,)$ denote the relative entropy from Sanov's theorem. Then the rate function \eqref{entropyrate} of Theorem~\ref{Th:2Rn} 
has the following representation 
\be
I(\varphi)=\left\{
\begin{array}{ll}
H(\mu_\varphi|\mu_o)&\mbox{if $\varphi$  non decreasing with $\varphi(0)=0$ and $\varphi(T)=1$},\\[3mm]
\infty&\mbox{otherwise}.
\end{array}
\right.
\ee
We stress that we do \textnormal{not} use Sanov's theorem during the proof of Theorem~\ref{Th:2Rn}. 
\end{rem}

\subsubsection{Third main result: asymptotic empirical path properties (Theorem~\ref{Th:3Rn})}\label{sss:empirical}

Finally, we analyze the path properties and respective jumps characteristics that actually lead to the event $\{\e L_{r_\e T}=\x\}$. Recall that by Hypothesis II the jump measure $\nu_\xi$ needs not be finite. In the case $\nu_\xi(\bR^n)=\infty$ it is known that the process $\xi$ almost surely has an infinite number of jumps on any positive time interval. Obviously, under this setting it is not particularly insightful to estimate the number of jumps and the distribution of the jump sizes of the event $\{\e L_{r_\e T}=\x\}$. In addition the jump tails of $\xi$ are lighter than those of $\eta$. Hence it is natural to assume $\xi=0$ in this subsection. Thus we set $L=\eta$ to be a L\'evy process with  generating triplet $(0,\nu_\eta,0)$, where~$\nu_\eta$ satisfies the Hypothesis~I. In order to circumvent unnecessary technicalities, we assume  the function~$f$ to be convex and monotonically non-decreasing.

\begin{defin}\label{defN}
We condition the process $L$ on the event $\{\e L_{r_\e T}=\x\}$ and keep the notation of $X^\e$ and $Y^\e$ from  Subsection \ref{sss:bridgeLDP}.
\begin{enumerate}
 \item Let $N^{\x,\e}$ denote the \textbf{number of jumps of $Y^\e$ on the interval of time $[0,T]$}. 
 \item For $i=1,\dots,N^{\x,\e}$ let $W^{\x,\e}_i\in\bR^n$ denote the \textbf{$i$-th jump increment of $Y^\e/\e$}. 
\end{enumerate}
\end{defin}
It is well-known that the family $(W^{\x,\e}_i)_{i=1,\dots,N^{\x,\e}}$ conditioned on $\{N^{\x,\e}=m\}$ is identically distributed for any fixed $m\in\bN$. In our case it can be read for instance from the common density of all jumps in formula \refq{4111Rn}. 
In order to determine the limiting laws of $N^{\x,\e}$ and $W^{\x,\e}_i$, both random variables must be scaled suitably. 
\begin{theorem}[The asymptotic empirical path properties of the bridge $Y^\e$]\label{Th:3Rn}\hfill
\begin{enumerate}
 \item[(i)] Let the speed function $S$ be defined as in Theorem~\ref{Th:2Rn}. Set $m_{\x,\e}:=g(\frac{|\x|\e^{-1}}{r_\e T})^{-1}|\x|\e^{-1}$, $k_\e:=\alpha |\x|^{-1} \e  g(\frac{|\x|\e^{-1}}{r_\e T})|\ln\e|$ and $\bar N^{\x,\e}:=\sqrt{k_\e}(N^{\x,\e}-m_{x,\e})$. Then the family $(\mathbf{Q}^{\x, \e})_{\e>0}$, $\mathbf{Q}^{\x,\e} = \Law(\sqrt{S(\e)}\bar N^{\x,\e})$ satisfies a LDP on $(\bR,|\cdot|)$ with speed function $S(\e)$ and good rate function \hbox{$J(y)=\frac12 y^2$.}
 \item[(ii)] For 
 \be\label{eq:223}
 \bar W^{\x,\e}:=\sqrt{f''(g(\frack{|\x|\e^{-1}}{r_\e T}))}~\Big(W^{\x,\e}-g(\frack{|\x|\e^{-1}}{r_\e T})\frac{\x}{|x|}\Big) \ee
 let the family of random variables 
 $(\mathbb{W}^{\x,\e})_{\e>0}$ be defined by 
 \be\label{224}
\mathbb{W}^{\x,\e}:=\langle\bar W_1^{\x,\e},\x\rangle\frac{\x}{|\x|}+(\alpha-1)^{-1}\Big(\bar W_1^{\x,\e}-\langle\bar W_1^{\x,\e},\x\rangle\frac{\x}{|\x|}\Big)
.
\ee
Then $\mathbb{W}^{\x,\e}$ converges as $\e\to 0$ in distribution to a standard normally distributed random vector on $\bR^n$.
\end{enumerate}
\end{theorem}

\begin{rem}
If for any $m\in \bN$, we condition $Y^{\x, \e}$ on exactly $N^{\x, \e} = m$ jumps,  
we have the representation \eqref{4111Rn} of the common density $(W^{\x,\e}_1, \dots,  W^{\x, \e}_m)$, 
which is invariant under index permutations. This implies the stationarity of the increments. 

In the case $\xi \neq 0$ one might define $N^{\x, \e}$ and $W^{\x, \e}_i$ as the jump frequency and jump increments of the compound Poisson component $\eta$ instead of $L = \eta +\xi$. Under this consideration the results of Theorem~\ref{Th:3Rn} remain valid. Since the proof is lengthy and mainly of technical nature without additional insights, we have omitted this result.  

\end{rem}

\bigskip

\subsection{Sufficient conditions for L\'evy perturbations $\xi$ of the CPP $\eta$ 
(Hypothesis~II)}\label{ss:SC}

For L\'evy processes $\xi$ on $\bR$ it is known that, if there exists a $\beta\in(0,2)$, such that the corresponding measure $\nu_\xi$ satisfies the Orey condition
\be
\liminf_{x\rightarrow0}x^{-\beta}\int_{-x}^xy^2\nu_\xi(dy)>0,\label{eq:Orey}
\ee
then for any $t>0$ the distribution of $\xi_t$ has a density. This result has been proven by S.~Orey~\cite{Orey}. For generalizations we refer to \cite{HK07} and references therein. In order to formulate and prove sufficient conditions for Hypothesis~II we need \hbox{a $\bR^n$ version} of~\refq{eq:Orey}. Thus we start with the following lemma.

\begin{lemma}\label{l:Orey_Rn}
Let $n\in\bN$ and $(v_1,\dots,v_n)$ be an orthonormal base of $\bR^n$. Let $\beta\in(0,2)$ and a L\'evy process $\xi$ on $\bR^n$ has a generating triplet $(\sigma^2,\nu,\Gamma)$, such that for every $i=1,\dots,n$ the following estimate is satisfied
\be\label{eq:Orey_Rn}
\lim_{r\rightarrow0+}r^{-\beta}\int_{|y|<r}\langle v_i,y\rangle^2\nu(dy)>0.
\ee
Then for any $t>0$ the distribution of $\xi_t$ has a density $\mu_t$, for which we have 	
\be
\lim_{t\rightarrow0+}t^{\frac{n}{2-\beta}}\sup_{y\in\bR^n}\mu_t(y)\spac{7pt}<\infty.
\ee
\end{lemma}
Obviously, \refq{eq:Orey_Rn} is a $\bR^n$-valued version of the Orey condition~\refq{eq:Orey}. Not surprisingly, the main argument of the proof given in Appendix~\ref{a:ProofsSC} is a direct translation of the original calculation by Orey~\cite{Orey} to the $\bR^n$ setting.

\begin{lemma}\label{hinr_bed2} (i) Let two independent L\'evy processes $\xi^1$ and $\xi^2$ each satisfy Hypothesis~II with parameters $\aleph_1$ and~$\aleph_2$ respectively. 
Then the process $\xi=\xi^1+\xi^2$ satisfies Hypothesis~II with any parameter $\aleph\in(\frac1\alpha,\aleph_1\wedge \aleph_2)$.\\

\noindent $(ii)$ Let $\xi$ be a L\'evy process with generating triplet $(\sigma^2,\nu_\xi,\Gamma)$. If one of the following conditions is satisfied, then $\xi$ also satisfies Hypothesis~II. 
\begin{enumerate}
	\item[III] There exists $\Lambda>0$, such that $\nu_\xi(\{y\in\bR^n\mid|y|>\Lambda\})=0$. Furthermore, one of the following condition is satisfied: 
	\begin{enumerate}
			\item $\det \sigma^2>0$.
			\item There is a parameter \hbox{$\beta\in(0,2)$,} such that $\nu_\xi$ satisfies \refq{eq:Orey_Rn}.
	\end{enumerate}
	\item[IV]	There exists a parameter $\alpha_\xi>\alpha$ and a non decreasing  function $f_\xi\in R_{\alpha_\xi}$, such that for each~$\Lambda>0$ we have 
	\be\label{4_15}
	\nu_\xi(\{y\in\bR^n\mid|y|>\Lambda\})\le\exp(-f_\xi(\Lambda)).
	\ee
	Moreover, $\det \sigma^2>0$ and there are  $\Delta,K>0$ such that for any subset $A\subset\{y\in\bR^n||y|>\Delta	\}$  we have the implication
\be\label{4_16}
\nu_\xi(A)\le K\lambda_{n-1}(A).
\ee
\end{enumerate}
\end{lemma}

\begin{rem} $(i)$ To understand the importance of condition \refq{4_16} we sketch an example of a Lévy process $\xi$ that satisfies \refq{4_15} and $\det\sigma^2>0$, but violates \refq{BEDIIRn} in Hypothesis II:  consider a Lévy process $\xi$ on $\bR^n$, $n>2$, with generating triplet $(\sigma^2,\nu_\xi,0)$, $\det \sigma^2>0$.  Assume that $\nu_\xi(\bR^n)<\infty$  and there is a sequence $(y_i)_{i\in\bN}\in\bR^n$ with $\lim_{i\rightarrow\infty}|y_i|=\infty$ and $\nu_\xi(\{y_i\})>0$. Note that this assumption does not contradict condition~\refq{4_15}. Let $N_t$ denote the number of jumps of $\xi$ in the time interval~$(0,t)$ and $\mu_{W,t}$ the marginal density of the Brownian component. We have $\lim_{t\to0+}t^{-1}P(N_t=1)>0$ and $\lim_{t\to0+}t^{\frac n2}\mu_{W,t}(0)>0$.  For every $i\in\bN$ we obtain
\be\label{231}
\lim_{t\to0+}\mu_{\xi,t}(y_i)\spac{7pt}\ge\lim_{t\to0+}P(N_t=1)\nu_{\xi}(\{y_i\})\mu_{W,t}(0)\spac{7pt}=\infty.
\ee
By $\lim_{i\rightarrow\infty}|y_i|=\infty$ the limit \refq{231} contradicts \refq{BEDIIRn}. Therefore \refq{4_15} together with $\det\sigma^2>0$ is not sufficient to obtain Hypothesis II.

$(ii)$ In the case $n=1$ the condition \refq{4_16} is void and for $n\ge2$ rotational invariance of $\nu_\xi$ is sufficient to grant condition \refq{4_16}. At the same time, the scope of the condition~\refq{4_16} is much larger than those examples.

$(iii)$ Note that, in the counterexample given in $(i)$ the condition \refq{4_16} is violated even if the Lebesgue measure $\lambda_{n-1}$ is replaced by $\lambda_{n-2}$. It remains an open question, whether the condition~\refq{4_16} can be relaxed and to which extent precisely.

\end{rem}

\bigskip
\section{Proofs of the main results}\label{s:proofs}\setcounter{equation}{0}

\subsection{Convolution estimates (Proposition~\ref{l:2Rn})}\label{ss:convolution}
\noindent\textbf{Proof of Proposition~\ref{l:2Rn}:} 
We establish \eqref{Faltung_FormelRn} for all $m\in \bN$. 
In the case $m=1$ we have $\frac{\nu^{*m}(dx)}{dx} = e^{-f(|x|)}$ and statement \eqref{Faltung_FormelRn} is clearly satisfied.

We continue with the case $m\ge2$. Due to the rotational invariance of $\nu$, which implies the rotational invariance of $\nu^{*m}$ for any $m\ge 2$, 
it is enough to consider values $x$ on the positive semi-axis $x = |x| e_1$. Before we actually target the convolution integral, we define for every $m\ge 2$ and $x = |x| e_1$ the important auxiliary function $f_{x,m}:\bR^n\setminus\{-\frac xm\}\rightarrow\bR$:
\begin{equation}
f_{x,m}(z)=f(|\frack xm+z|)-(f(\frack{|x|}m)+z_1f'(\frack{|x|}m)).\label{deffxm}
\end{equation}
Using $f_{x,m}$ we may formulate the convolution integral as 
\begin{\eq}
&&\frac{\nu^{*m}(dx)}{dx}=\int_{\bR^n\setminus\{0\}}\cdots\int_{\bR^n\setminus\{0\}}\nonumber
\exp\Big(-\sum_{i=1}^{m-1}f(|y_i|)-f\Big(\big|x-\sum_{i=1}^{m-1}y_i\Big|\Big)\Big)dy_1\dots dy_{m-1}\\
&=& \int_{\bR^n\setminus\{\frack{x}{m}\}}\cdots\int_{\bR^n\setminus\{\frac{x}{m}\}} 
\exp\bigg(-\sum_{i=1}^{m-1}f\Big(|\frac{x}{m}+y_i|\Big)-f\Big(\big|\frac{x}{m}-\sum_{i=1}^{m-1}y_i\big|\Big)\bigg)
dy_1\dots dy_{m-1}\label{4.5rn}\\
&=&\exp\Big(-mf(|\frack{x}{m}|)\Big)\int_{\bR^n\setminus\{\frack{x}{m}\}}\cdots\int_{\bR^n\setminus\{\frac{x}{m}\}}
\exp\Big(-\sum_{i=1}^{m-1}f_{x,m}(y_i)-f_{x,m}\Big(-\sum_{i=1}^{m-1}y_i\Big)\Big)dy_1\dots dy_{m-1},\nonumber
\end{\eq}
where the seemingly missing summands $\sum_{i=1}^{m-1} (-y_i) f'(|\frack{x}{m}|)  -(-\sum_{i=1}^{m-1} y_i) f'(|\frack{x}{m}|)$ add up to zero.

Thus, in order to establish the desired upper and lower bounds for the convolutions of $\nu$, we first establish an upper and lower bound for $f_{x,m}$ in the vicinity of the origin. Set $c\in(1-\frac\alpha2,1)$ and $\delta>0$. For a sufficiently large $k>0$ we assume $\frac{|x|}m>k$ and estimate the value of $f_{x,m}(z)$ as follows. By the choice of $x = |x| e_1$ and $c<1$, and for $|z|\le(\frac{|x|}m)^c$ we obtain $|\frac xm+z_1e_1|=\frac{|x|}m+z_1$, and by definition of $f_{x,m}$ we have
\ba
f_{x,m}(z) &= f(|\frack xm+z|) -f(\frack{|x|}{m})-z_1f'(\frack{|x|}m))\\
&=(f(|\frack xm+z|)-f(|\frack{x}m+z_1e_1|))+(f(\frack{|x|}m+z_{1})-(f(\frack{|x|}m)+z_1f'(\frack{|x|}m))).\label{331}
\ea
Hypothesis~I on $f$ yields that for any $\delta>0$, $c\in(1-\frac\alpha2,1)$ there is $\frac{|x|}{m}$ sufficiently large such that for~$|z|\le(\frac{|x|}m)^c$ it follows
\be\label{332}
f(\frack{|x|}{m}+z_1)-(f(\frack{|x|}{m})+z_1 f'(\frack{|x|}{m}))\spac{7pt}\le\sup_{|q|\le(\frac{|x|}{m})^c}\frack{1}{2} f''(\frack{|x|}{m}+q)z_1^2\spac{7pt}\le(1+\delta)\frack{1}{2}f''(\frack{|x|}{m})z_1^2.
\ee
To estimate the first summand on the right side of \refq{331} we apply  Lemma~\ref{SR}(ii) together with the estimate $\sqrt{r+s}-\sqrt{r}\leq\frac{1}{2}\frac{s}{\sqrt{r}}$ for all $r, s\geq 0$. For $\frac{|x|}m$ sufficiently large and $|z|<(\frac{|x|}m)^c$ we have
\ba\label{333}
f(|\frack xm+z|)&-f(|\frack{x}m+z_1e_1|)\spac{7pt}\le\sup_{|q|\le(\frac{|x|}{m})^c}\frack{1}{2} f'(\frack{|x|}{m}+q)(|\frack xm+z|-|\frack{x}m+z_1e_1|)\\
&\spac{7pt}\le(1+\delta)f'(\frack{|x|}m)\frack12\frack m{|x|}|(0,z_2,\dots,z_n)|^2\spac{7pt}\le(1+2\delta)\frack1{\alpha-1}\frack12f''(\frack{|x|}m)\sum_{i=2}^nz_i^2.
\ea
Inserting \refq{332} und \refq{333} into \refq{331} we obtain an upper bound for $f_{x,m}$. \noindent The corresponding lower bound of $f_{x,m}$ can be estimated similarly. Hence for any $\delta>0$, $c\in(1-\frac\alpha2,1)$ there is $k>0$ such that for any $\frac{|x|}{m}>k$ and $|z|<(\frac{|x|}m)^c$ we have 
\begin{align}
\Big|f_{x,m}(z)-\frack{1}{2}f''(\frack{|x|}{m})\Big(z_{1}^2+\frack{1}{\alpha-1}\sum_{i=2}^{n}z_{i}^2\Big)\Big|
&\le\delta|z|^2f''(\frack{|x|}{m}).
\label{421}
\end{align}\begin{rem}
The different asymptotics in direction $\x$ and $\x^\perp$, which are present  in \refq{eq:223} of Theorem~\ref{Th:3Rn}(ii) originate in the preceding estimate~\eqref{421} with the different regimes of the coordinates~$z_1$ and $z_{2}, \dots, z_{n}$.
\end{rem}

\noindent Estimate \refq{421} was established for any $c\in(1-\frac\alpha2,1)$, $|z|<(\frac{|x|}m)^{c}$ and $x=(x_1,0,\dots,0)$ with $\frac{|x|}m$ sufficiently large. Obviously \refq{421} remains valid for any such $z\in[-(\frac{|x|}m)^c,(\frac{|x|}m)^c]^n$ and $c\in(1-\frac\alpha2,1)$.\vspace{2mm}

\noindent We start with the direct proof the upper bound, followed by the more involved estimate from below.\\ 

\noindent\textbf{Upper bound of \eqref{Faltung_FormelRn}: } 
By \refq{4.5rn} we have 
\ba\label{eq:4.29}
\frac{\nu^{*m}(dx)}{dx}&\spac{7pt}\le\exp(-mf(|\frack xm|))\int_{\bR^n\setminus\{\frac{|x|}m\}}\cdots\int_{\bR^n\setminus\{\frac{|x|}m\}}\exp\Big(-\sum_{i=1}^{m-1}f_{x,m}(y_i)\Big)dy_1\dots dy_{m-1}\\
&\spac{7pt}=\exp(-mf(|\frack xm|))\Big(\int_{\bR^n\setminus\{\frac{|x|}m\}}\exp(-f_{x,m}(y))dy\Big)^{m-1}.
\ea
We estimate the value of the integral for $y\in [-(\frac{|x|}m)^c,(\frac{|x|}m)^c]^n$ and $y\notin [-(\frac{|x|}m)^c,(\frac{|x|}m)^c]^n$ separately. In the first case we may apply \refq{421} to estimate $f_{x,m}(y)$. For $\frac{|x|}{m}$ large enough, and \hbox{$y\in [-(\frac{|x|}m)^c,(\frac{|x|}m)^c]$,} we estimate $f_{x,m}(y)>(1-\delta)\frac12f''(\frac{|x|}{m})(y_1^2+\frac1{\alpha-1}\sum_{i=2}^ny_i^2)$. Then the Gaussian renormalization $\sqrt{2\pi a} = \int_{\bR} \exp(\frac{-s^2}{2 a}) ds$, $a>0$, implies 
\ba\label{eq:4.30}
\int_{\big[-(\frac{|x|}m)^c,(\frac{|x|}m)^c\big]^n}\exp(-f_{x,m}(y))dy&\spac{7pt}\le\int_{\bR^{n}}\exp\Big(-(1-\delta)\frack12f''(\frack{|x|}{m})\Big(y_1^2+\frack1{\alpha-1}\sum_{i=2}^ny_i^2\Big)\Big)dy\\
&\spac{7pt}=(\alpha-1)^{\frac{n-1}2}(2\pi((1-\delta)f''(\frack{|x|}m))^{-1})^{\frac n2}.
\ea

\noindent It remains to estimate the integral for $y\in \bR^n\setminus[-(\frac{|x|}m)^c,(\frac{|x|}m)^c]^n$ for $|y|\ge(\frac{|x|}m)^c$. By the definition of $f_{x,m}$ we have
\ba\label{3_39}
\frac{d^2}{ds^2}f_{x,m}(sy)&=\frac{d^2}{ds^2}\Big(f(|\frack xm+sy|)-f(|\frack xm|)-sy_1f'(|\frack xm|)\Big)\\
&=\Big(\frac{d^2}{ds^2}|\frack xm+sy|\Big)f'(|\frack xm+sy|)+\Big(\frac{d}{ds}|\frack xm+sy|\Big)^2f''(|\frack xm+sy|),
\ea
\noindent for any $y\in\bR^n\setminus\{0\}$, $s>0$, such that $\frac xm+sy\neq0$. We start with the case of $f$ being convex and non-decreasing. In that case the right side of \refq{3_39} is non-negative. And hence together with $f_{x,m}(0)=0$ we have
\be\label{convexfxm}
f_{x,m}(sz)\spac{8pt}\ge sf_{x,m}(z)
\ee
for any $z\in\bR^n$ and $s>1$. Next we lift the condition of $f$ being non-decreasing, while still being convex. By Hypothesis~I, in particular, $f\in SR_\alpha$ and $\alpha>1$, we have that $r\mapsto f(r)$ is monotonically non-decreasing for $r>r_o$ and some $r_o>0$. Let $\tilde f:[0,\infty)\mapsto\bR$ be a monotonically non-decreasing function with $f\ge\tilde f$ on $[0,r_o]$ and $f=\tilde f$ on $[r_o,\infty)$, and let $\tilde f_{x,m}$ be defined accordingly. Then $\tilde f_{x,m}$ satisfies \refq{convexfxm} by construction. For $\frac{|x|}m$ sufficiently large we obtain

\be\label{311}
f_{x,m}(sz)\spac{8pt}\ge\tilde f_{x,m}(sz)\spac{8pt}\ge s\tilde f_{x,m}(z)\spac{8pt}= sf_{x,m}(z)
\ee
for $|z|\le(\frac{|x|}m)^c$ and $s>1$. We apply \refq{311} to estimate $f_{x,m}(y)$ for $|y|\ge(\frac{|x|}m)^c$. By the choice of~$c$ we have $\alpha-2+2c>0$. Choose $\delta\in(0,\frac12(\alpha-2-2c))$. Set $s=|y|(\frack{|x|}m)^{-c}>1$ and $z=\frac ys$. We obtain $y=sz$ and $|z|=(\frack{|x|}m)^c$, which allows to estimate $f_{x,m}(z)$ by~\refq{421}. Together with $f''\in SR_{\alpha-2}$ we obtain
\ba
f_{x,m}(y)&\spac{8pt}\ge sf_{x,m}(z)\spac{8pt}=|y|(\frack{|x|}m)^{-c}f_{x,m}(z)\spac{8pt}\ge(1-\delta)|y|(\frack{|x|}m)^{-c}\min\{1,\frack1{\alpha-1}\}\frack12f''(\frack{|x|}m)|z|^2\\
&\spac{8pt}\ge(\frack{|x|}m)^{\alpha-2+c-\delta}|y|\spac{8pt}=(\frack{|x|}m)^{\alpha-2+2c-\delta}+(\frack{|x|}m)^{\alpha-2+c-\delta}(|y|-(\frack{|x|}m)^c).
\ea
By the choice of $c$ and $\delta$ we may estimate the exponents $\alpha-2+2c-\delta>\delta$ and $\alpha-2+c-\delta>-1$. For $\frac{|x|}m$ sufficiently large we obtain
\ba\label{eq:4.31}
&\int_{\bR^n}\exp(- f_{x,m}(y))\textbf{1}_{[(\frac{|x|}m)^c,\infty)}(|y|)dy\\
&\le\int_{\bR^n}\exp\Big(-\Big((\frack{|x|}m)^\delta+(\frack{|x|}m)^{-1}(|y|-(\frack{|x|}m)^c)\Big)\Big)\textbf{1}_{[(\frac{|x|}m)^c,\infty)}(|y|)dy\spac{8pt}\le\exp(-(\frack{|x|}m)^{\frac\delta2}).
\ea
\noindent Finally, we lift the condition of $f$ being convex. By $f\in SR_\alpha$, $\alpha>1$, there is some $r_o>0$, such that~$f$ is convex on $[r_o,\infty)$. Let $\tilde f:[0,\infty)\to\bR$ be convex and $\tilde f=f$ on $[r_o,\infty)$. Let $\tilde f_{x,m}$ be defined analogously to~$f_{x,m}$. Then we have
\ba
&\int_{\bR^n}\exp(- f_{x,m}(y))\textbf{1}_{[(\frac{|x|}m)^c,\infty)}(|y|)dy\\
\le&\int_{\bR^n}\exp(-\tilde f_{x,m}(y))\textbf{1}_{[(\frac{|x|}m)^c,\infty)}(|y|)dy+\int_{\bR^n}\exp(- f_{x,m}(y))\textbf{1}_{[0,r_o]}(|y+\frack xm|)dy,
\ea
where the first integral on the right-hand side can be estimated by from above by the right-hand side of \refq{eq:4.31}. To estimate the second integral we use the definition of $f_{x,m}$ together with $\int_{\bR^n}\exp(-f(|x|))dx=\nu(\bR^n)=1$ and Lemma \ref{SR}$(ii)$. For $\frac{|x|}m$ sufficiently large we obtain 
\ba\label{315}
&\int_{\bR^n}\exp(- f_{x,m}(y))\textbf{1}_{[0,r_o]}(|y+\frack xm|)dy\\
=&\int_{\bR^n}\exp(-f(|y+\frack xm|)+f(\frack {|x|}m)+y_1f'(\frack {|x|}m))\textbf{1}_{[0,r_o]}(|y+\frack xm|)dy\\
\le&\exp(-\frack12(\alpha-1)f(\frack {|x|}m))\int_{\bR^n}\exp(-f(|y+\frack xm|))dy\spac{7pt}\le\exp(-\frack {|x|}m).
\ea
The desired upper bound of \refq{Faltung_FormelRn} follows combining the estimates \refq{eq:4.29}, \refq{eq:4.30} and \refq{eq:4.31} - \refq{315}.\vspace{2mm}

\noindent\textbf{Lower bound of \eqref{Faltung_FormelRn}: } To calculate the convolution $\frac{\nu^{*m}(dx)}{dx}$, originally an $(m-1)$-fold integral must  be calculated, where each of the integration parameters $y_1,\,y_2,\dots,y_{m-1}$  is $\bR^n$-valued. We will transform this into an $n$-fold integral whose integration parameters $\tilde y_1,\,\tilde y_2,\dots,\tilde y_n$ are each $\bR^{m-1}$-valued. We then carry out a substitution of the integration parameters: let $\bar y_1=\tilde y_1$ and  for $i=2,\,3,\dots,n$ let $\bar y_i=\frac1{\sqrt{\alpha-1}}y_i$. Finally, we succeed in estimating this $n$-fold integral as the $n$-th power of the simple integral over $\bR^{m-1}$. To simplify the notation we use an auxiliary function
$$\theta_{m-1}:\bR^{m-1}\rightarrow[0, \infty), \qquad 
\theta_{m-1}(z):=\sum_{i=1}^{m-1}z_i^2+\Big(\sum_{i=1}^{m-1}z_i\Big)^2.$$ 
Fix $c\in(1-\frac\alpha2,1)$, $\tilde c\in(1-\frac\alpha2,c)$. Let $U_{x,m}:=\{y\in\bR^{m-1}||\sum_{i=1}^{m-1}y_i|<(\frac{|x|}m)^{\tilde c}\}$. By \refq{4.5rn}, \refq{421} and the nonnegativity of the integrand we have
\ba\label{eq:342}
&\exp(mf(|\frack xm|))\cdot\frac{\nu^{*m}(dx)}{dx}\spac{7pt}=\int\limits_{\bR^n}\cdots\int\limits_{\bR^n}\exp\Big(-\sum_{i=1}^{m-1}f_{x,m}(y_i)-f_{x,m}\Big(-\sum_{i=1}^{m-1}y_i\Big)\Big)dy_1\dots dy_{m-1}\\
\ge&\int\limits_{\big[-(\frac{|x|}m)^c,(\frac{|x|}m)^c\big]^n}\cdots\int\limits_{\big[-(\frac{|x|}m)^c,(\frac{|x|}m)^c\big]^n}\textbf{1}_{\big[-(\frac{|x|}m)^c,(\frac{|x|}m)^c\big]^{n}}\Big(\sum_{i=1}^{m-1}\big(y_{i,1}, \dots, y_{i,n}\big)\Big)\cdot\\
&\cdot\exp\bigg[-\frac{1+\delta}2f''(\frac{|x|}m)\Big(\theta_{m-1}(y_{1,1},\dots,y_{m-1,1})+\frac1{\alpha-1}\sum_{i=2}^n\theta_{m-1}(y_{i,1},\dots,y_{m-1,i})\Big)\bigg]dy_1\dots dy_{m-1}\\
=&\int\limits_{\big[-(\frac{|x|}m)^c,(\frac{|x|}m)^c\big]^{m-1}}\cdots\int\limits_{\big[-(\frac{|x|}m)^c,(\frac{|x|}m)^c\big]^{m-1}}\textbf{1}_{\big[-(\frac{|x|}m)^c,(\frac{|x|}m)^c\big]^{n}}\Big(\sum_{i=1}^{m-1}\big(\tilde y_{1,i}, \dots,\tilde y_{n,i}\big)\Big)\cdot\\
&\cdot\exp\bigg[-\frac{1+\delta}2f''(\frac{|x|}m)\Big(\theta_{m-1}(\tilde y_{1,1},\dots,\tilde y_{1,{m-1}})+\frac1{\alpha-1}\sum_{i=2}^n\theta_{m-1}(\tilde y_{i,1},\dots,\tilde y_{i,{m-1}})\Big)\bigg]d\tilde y_1\dots d\tilde y_{n}\\
\ge&(\alpha-1)^{\frac12(n-1)(m-1)}\cdot\int\limits_{\big[-(\frac{|x|}m)^{\tilde c},(\frac{|x|}m)^{\tilde c}\big]^{m-1}}\cdots\int\limits_{\big[-(\frac{|x|}m)^{\tilde c},(\frac{|x|}m)^{\tilde c}\big]^{m-1}}\prod_{j=1}^n\textbf{1}_{\big[-(\frac{|x|}m)^{\tilde c},(\frac{|x|}m)^{\tilde c}\big]}\Big(\sum_{i=1}^{m-1}\bar y_{i,j}\Big)\cdot\\
&\cdot\exp\bigg[-\frac{1+\delta}2f''(\frac{|x|}m)\sum_{i=1}^n\theta_{m-1}(\bar y_{i,1},\dots,\bar y_{i,{m-1}})\bigg]d\bar y_1\dots d\bar y_{n}\\
\ge&(\alpha-1)^{\frac12(n-1)(m-1)}\Bigg(\int\limits_{\big[-(\frac{|x|}m)^{\tilde c},(\frac{|x|}m)^{\tilde c}\big]^{m-1}}\exp\bigg[-\frac{1+\delta}2f''(\frac{|x|}m)\theta_{m-1}(z_{1},\dots,z_{m-1})\bigg]\textbf{1}_{U_{x,m}}(z)dz\Bigg)^n.
\ea
Note that the integral on the right-hand side can be compared to the well-known convolution of normally distributed random variables. Indeed, it is not hard to see that 
\be
\int_{\bR^{m-1}}\big(\frack a{2\pi}\big)^\frac m2\exp(-\frack a2\theta(y))dy=\sqrt{\frack a{2m\pi}},\mb{10pt}{thus}\int_{\bR^{m-1}}\exp(-\frack a2\theta(y))dy=\frack 1{\sqrt m}\big(\frack{2\pi}a\big)^\frac{m-1}2
\ee
for any $a>0$. By the choice $a=(1+\delta)f''(\frac{|x|}m)$ the proof of \refq{Faltung_FormelRn} boils down to examine the effect of the reduced area of integration $[-(\frac{|x|}m)^{\tilde c},(\frac{|x|}m)^{\tilde c}]^{m-1}\cap U_{x,m}$.

Consider the $m-2$ dimensional linear subspace $V_{m-1}:=\{y\in\bR^{m-1}\mid\sum_{i=1}^{m-1}y_i=0\}$ of $\bR^{m-1}$. We split the integration in the direction of $V_{m-1}$ and its orthogonal complement.  For this purpose let $(v_{1},\dots,v_{m-2})$ be an orthonormal basis for $V_{m-1}$. Let $v_{m-1}:=\frac1{\sqrt{m-1}}(1,1,\dots,1)$ and $\tilde v_{m-1}:=\frac{1}{\sqrt{m}}v_{m-1}$. Thus, by construction we obtain  an orthonormal basis $(v_{1},\dots,v_{m-1})$  for $\bR^{m-1}$ and $(v_{1},\dots,v_{m-2},\tilde v_{m-1})$  is an orthogonal basis for $\bR^{m-1}$ such that for $y \in\bR^{m-1}$ and $z=\sum_{i=1}^{m-2}y_iv_i+y_{m-1}\tilde v_{m-i}$ we have
\ba
\theta_{m-1}\Big(\sum_{i=1}^{m-2}\bar y_i v_{i}+\bar y_{m-1} \tilde v_{m-1}\Big)&=|z|^2+\Big(\sum_{i=1}^{m-1}z_i\Big)^2\\
&=\Big(\sum_{i=1}^{m-2}y_i^2+\frack 1{m}y_{m-1}^2\Big)+\Big(y_{m-1}\sqrt{\frack{m-1}m}\Big)^2=|y|^2.\label{eq:343}
\ea
We continue to estimate the remaining integral. By the definition of $v$ and $\tilde v$ together with \refq{eq:343} we obtain
\ba\label{eq:3.44}
&\int\limits_{\big[-(\frac{|x|}m)^{\tilde c},(\frac{|x|}m)^{\tilde c}\big]^{m-1}}
\exp\Big(-(1+\delta)f''(\frack{|x|}m)\theta_{m-1}( y)\Big)\textbf{1}_{U_{x,m}}( y)dy\\
&=\int_{\bR^{m-1}}\exp\Big(-(1+\delta)f''(\frack{|x|}m)\theta_{m-1}\Big(\sum_{i=1}^{m-1} y_iv_i\Big)\Big)\textbf{1}_{\big[-(\frac{|x|}m)^{\tilde c},(\frac{|x|}m)^{\tilde c}\big]^{m-1}}\Big(\sum_{i=1}^{m-1} y_iv_{i}\Big)\textbf{1}_{U_{x,m}}\Big(\sum_{i=1}^{m-1} y_iv_{i}\Big)dy\\
&\ge\int_{\bR^{m-1}}\frac{\exp\Big(-(1+\delta)f''(\frack{|x|}m)|y|^2\Big)}{\sqrt m}\textbf{1}_{\big[-(\frac{|x|}m)^{\tilde c},(\frac{|x|}m)^{\tilde c}\big]^{m-1}}\Big(\sum_{i=1}^{m-2}y_iv_{i}+y_{m-1}\tilde v_{m-1}\Big)\textbf{1}_{[0,(\frac{|x|}m)^{\tilde c}]}(|y_{m-1}|)dy.
\ea
In the last step we carry out base change of $\bR^{m-1}$ from $(v_1,\dots,v_{m-1})$ to $(v_1,\dots,v_{m-2},\tilde v_{m-1})$ and apply \refq{eq:343}. Furthermore, we use that by the definition of $U_{x,m}$ and $\tilde v_{m-1}$ we have
$$\textbf{1}_{U_{x,m}}\Big(\sum_{i=1}^{m-2} y_iv_{i}+y_{m-1}\tilde v_{m-1}\Big)\spac{7pt}=\textbf{1}_{[0,(\frac{|x|}m)^{\tilde c}]}\Big(\sqrt{\frac{m-1}{m}}|y_{m-1}|\Big)\spac{7pt}\ge\textbf{1}_{[0,(\frac{|x|}m)^{\tilde c}]}(|y_{m-1}|).$$
We start the integration in the $\tilde v_{m-1}$-coordinate. Fix $\bar c\in(1-\frac\alpha2,\tilde c)$. Then, $\frac{|x|}m$ can be chosen sufficiently large, such that for any $|y_{m-1}|<(\frac{|x|}m)^{\bar c}$ and $z\in V\cap[-(\frac{|x|}m)^{\bar c},(\frac{|x|}m)^{\bar c}]^{m-1}$  it follows that $(z+y_{m-1}\tilde v_{m-1})\in[-(\frac{|x|}m)^{\tilde c},(\frac{|x|}m)^{\tilde c}]^{m-1}$. Then we obtain 
\ba\label{eq:3.45}
&\int_{\bR^{m-1}}\exp\Big(-(1+\delta)f''(\frack{|x|}m)|y|^2\Big)\textbf{1}_{\big[-(\frac{|x|}m)^{\tilde c},(\frac{|x|}m)^{\tilde c}\big]^{m-1}}\Big(\sum_{i=1}^{m-2}y_iv_{i}+y_{m-1}\tilde v_{m-1}\Big)\textbf{1}_{[0,(\frac{|x|}m)^{\tilde c}]}(|y_{m-1}|)dy\\
&\ge\int_{-(\frac{|x|}m)^{\bar c}}^{(\frac{|x|}m)^{\bar c}}\Bigg(\int_{\bR^{m-2}}\exp\Big(-(1+\delta)f''(\frack{|x|}m)|\tilde y|^2\Big)\textbf{1}_{\big[-(\frac{|x|}m)^{\bar c},(\frac{|x|}m)^{\bar c}\big]^{m-1}}\Big(\sum_{i=1}^{m-2}\tilde y_iv_{i}\Big)d\tilde y\Bigg)\cdot\\
&\hspace{100mm}\cdot\exp(-(1+\delta)f''(\frack{|x|}m)y_{m-1}^2)dy_{m-1}\\
&\ge(1-\delta)\sqrt{\frac{2\pi}{f''(\frac{|x|}{m})}}\int_{\bR^{m-2}}\exp\Big(-(1+\delta)f''(\frack{|x|}m)|\tilde y|^2\Big)\textbf{1}_{\big[-(\frac{|x|}m)^{\bar c},(\frac{|x|}m)^{\bar c}\big]^{m-1}}\Big(\sum_{i=1}^{m-2}\tilde y_iv_{i}\Big)d\tilde y,\\
\ea
where in the last step we use $f''\in SR_{\alpha-2}$ together with the choice $\bar c>1-\frac\alpha2$. Indeed, similar to the result in \eqref{eq:4.31}, the integral on $\bR\setminus[-(\frac{|x|}m)^{\bar c},(\frac{|x|}m)^{\bar c}]$ is negligibly small such that for any $\delta>0$ and $\frac{|x|}{m}$ sufficiently large, we have (for $z = y_{m-1}$) that
\be\label{321}
\int_{-(\frac{|x|}m)^{\bar c}}^{(\frac{|x|}m)^{\bar c}}\exp(-(1+\delta)f''(\frack{|x|}m)z^2)dz\ge (1-\frack\delta2)\int_\bR\exp(-(1+\delta)f''(\frack{|x|}m)z^2)dz  \ge (1-\delta)\sqrt{\frac{2\pi}{f''(\frack{|x|}m)}}.
\ee
Finally we remember that $(v_1,\dots,v_{m-2})$ is an orthonormal base of $V$. Then we obtain
\ba\label{eq:3.46}
&\int_{\bR^{m-2}}\exp\Big(-(1+\delta)f''(\frack{|x|}m)|y|^2\Big)\textbf{1}_{\big[-(\frac{|x|}m)^{\bar c},(\frac{|x|}m)^{\bar c}\big]^{m-1}}\Big(\sum_{i=1}^{m-2}y_iv_{i}\Big)dy\\
=&\int_{\bR^{m-2}}\exp\Big(-(1+\delta)f''(\frack{|x|}m)\Big|\sum_{i=1}^{m-2}y_iv_{i}\Big|^2\Big)\textbf{1}_{\big[-(\frac{|x|}m)^{\bar c},(\frac{|x|}m)^{\bar c}\big]^{m-1}}\Big(\sum_{i=1}^{m-2}y_iv_{i}\Big)dy\\
\ge&\int_{\bR^{m-2}}\exp\Big(-(1+\delta)f''(\frack{|x|}m)|y|^2\Big)\textbf{1}_{\big[-(\frac{|x|}m)^{\bar c},(\frac{|x|}m)^{\bar c}\big]^{m-2}}(y)dy\\
=&\Big(\int_{-(\frac{|x|}m)^{\bar c}}^{(\frac{|x|}m)^{\bar c}}\exp\Big(-(1+\delta)f''(\frack{|x|}m)\tilde y^2\Big)d\tilde y\Big)^{m-2}\spac{7pt}\ge(1-\delta)^{m-2}(2\pi f''(\frack{|x|}m)^{-1})^{\frac{m-2}{2}}.
\ea
The second step of \refq{eq:3.46} relies on the following elementary fact: Let $n\in\bN$ and denote by $W_n$ the unit cube $[-1,1]^n$ in $\bR^n$. 
For $i=1,\dots,n$ let $V_i=\{y\in\bR^n|y_i=0\}$. Then for any linear $n-1$ dimensional linear subspace $V$ of $\bR^n$, any $k>0$ and any $i=1,\dots,n$ it follows 
$$\lambda_{n-1}(\{y\in V\cap W_n\mid |y|<k\})\ge\lambda_{n-1}(\{y\in V_i\cap W_n\mid |y|<k\}).$$
The last step of \refq{eq:3.46}, similarly to the last estimate of \refq{eq:3.45}, relies on \refq{321} together with the choice $\bar c>1-\frac\alpha2$.\\

\noindent The concatenation of estimates \refq{eq:342}, \refq{eq:3.44}, \refq{eq:3.45} and \refq{eq:3.46} establishes the lower bound of~\refq{Faltung_FormelRn}. Combining the upper and the lower bounds of \refq{Faltung_FormelRn} completes the proof of  Proposition~\ref{l:2Rn}.

\subsection{Proof of the estimate of the density $\mu_t$ (Theorem~\ref{th:1Rn})}\label{ss:density}
As crucial auxiliary results for the proof of Theorem~\ref{th:1Rn} we first establish Lemma~\ref{g_lnRn}.\vspace{2mm}

\noindent \textbf{Proof of Lemma~\ref{g_lnRn}. Item (i):} Combining Lemma~\ref{SR}$(v)$ and $(viii)$ yields $(\Lambda\mapsto \Lambda f''(\Lambda))\in SR_{\alpha-1}$. Due to $\frac{d}{d\Lambda}(\Lambda f'(\Lambda)-f(\Lambda))=\Lambda f''(\Lambda)$ combined with Lemma~\ref{SR}$(ix)$ it follows that $(\Lambda\mapsto \Lambda f'(\Lambda)-f(\Lambda))\in SR_\alpha$. By Lemma~\ref{SR}$(iii)$ it follows that $(\Lambda\mapsto \Lambda f'(\Lambda)-f(\Lambda)+k(\Lambda))\in SR_\alpha$. That is, by Lemma~\ref{SR}$(x)$ the left-hand side of \eqref{e:NLFE} is eventually invertible and the functions $h$ and $g$ are well defined with $h\in SR_{\frac1\alpha}$. By Lemma~\ref{SR}$(vi)$ and $g(\Lambda)=h(\ln \Lambda)$ we finally obtain $g\in SR_o$.\\

\noindent\textbf{Item (ii): } By Item $(i)$ we have that $g$ is ultimately non-decreasing with $\lim_{\Lambda\to\infty}g(\Lambda)=\infty$, $\lim_{\Lambda\to\infty}\frac{g(\Lambda^\gamma)}{g(\Lambda)}=\gamma^{\frac1\alpha}$, thus $\lim_{\Lambda\to\infty}\frac{g((\ln\Lambda)^\rho\Lambda)}{g(\Lambda)}=1$  for every $\gamma>0$, $\rho\in\bR$. Moreover, Lemma~\ref{SR}$(v)$ and $(viii)$ imply that $(x\mapsto xf''(x))\in SR_{\alpha-1}$ and $(x\mapsto|k'(x)|)\in SR_{\beta-1}$ with $\alpha>\beta$, therefore $\lim_{\Lambda\to\infty}\frac{k'(g(\Lambda))}{g(\Lambda)f''(g(\Lambda))}=0$. For any $\delta,\rho>0$ there is $\Lambda$ sufficiently large, such that for all $y\in[(\ln\Lambda)^{-\rho},(\ln\Lambda)^\rho]$ it follows 
\ba
\frack{d}{dy}\ln y=\frack{d}{dy}\ln(\Lambda y)&=\frack{d}{dy}\Big(g(\Lambda y)f'(g(\Lambda y))-f(g(\Lambda y))+k(g(\Lambda y))\Big)\\
&=\Lambda g'(\Lambda y)\Big(g(\Lambda y)f''(g(\Lambda y))+k'(g(\Lambda y))\Big)\\
&\ge(1-\delta)\Lambda g(\Lambda)g'(\Lambda y)f''(g(\Lambda y))\\
&=(1-\delta)g(\Lambda)\frack{d}{dy}f'(g(\Lambda y)).
\ea
By same arguments a similar upper bound is established. In case of $y=1$ both sides of \refq{eq.ln.gRn} equal $0$. Thus \refq{eq.ln.gRn} is satisfied for $y=1$. Therefore by $|g(\Lambda)\frack{d}{dy}f'(g(\Lambda y))-\frack{d}{dy}\ln y|\le\delta\frack{d}{dy}\ln y$ it follows that \refq{eq.ln.gRn} is valid all $y\in[(\ln \Lambda)^{-\gamma},(\ln \Lambda)^\gamma]$. This completes the proof of Lemma~\ref{g_lnRn}$(ii)$.\\

\noindent \textbf{Item (iii): } Again we use $(\Lambda\mapsto \Lambda f'(\Lambda)-f(\Lambda))\in SR_\alpha$. Together with $b<\alpha$, the limit~\eqref{e:lim1} is a direct consequence of the definition of $g$ in \refq{e:NLFE}. By Lemma~\ref{SR}$(ii)$ we have $$\lim\limits_{\Lambda\rightarrow\infty}\frac{g(\Lambda)f'(g(\Lambda))}{f(g(\Lambda))}=\alpha.$$ Thus \eqref{e:lim1} directly implies \eqref{e:lim2} and \eqref{e:lim3}. 
Let $h$ be defined as in the statement of Item (i). Then we have $g(\Lambda)=h(\ln \Lambda)$ and thus $g'(\Lambda)=\frac{1}{\Lambda}h'(\ln \Lambda)$. Lemma~\ref{g_lnRn}$(i)$ and Lemma~\ref{SR}$(ii)$ yield $$\lim\limits_{\Lambda\rightarrow\infty}\frac{g'(\Lambda)\Lambda\ln \Lambda}{g(\Lambda)}=\lim\limits_{\Lambda\rightarrow\infty}\frac{h'(\ln \Lambda)\ln \Lambda}{h(\ln \Lambda)}=\frac1\alpha,$$ which proves \eqref{e:lim4}. \\

\noindent\textbf{Item (iv):} For every $\delta>0$ there is $\Lambda_\delta$ sufficiently large, such that $k_1(\Lambda)\in[k_0(\Lambda)+c-\delta,k_0(\Lambda)+c+\delta]$ for every $\Lambda>\Lambda_\delta$. By definition of $g_i$ we obtain $g_o(x\exp(-c-\delta))<g_1(x)<g_o(x\exp(-c+\delta))$. By the statement of item $(i)$ we already have that $\lim\limits_{\Lambda\rightarrow\infty}\frac{g_o(\Lambda)}{g_1(\Lambda)}=1$. Consequently, it is sufficient to determine the limit value for one of the values $i\in\{0, 1\}$. Using part
$(ii)$, we obtain
$$g_o(\Lambda)(f'(g_o(\Lambda))-f'(g_1(\Lambda)))\spac{7pt}\ge g_o(\Lambda)(f'(g_o(\Lambda))-f'(g_o(\Lambda\exp(-c+\delta)))\spac{7pt}\ge c-2\delta.$$
for $\Lambda$ sufficiently large. By same arguments a similar upper bound is established and $\delta\to0$ yields the assertion. This completes the proof of Lemma~\ref{g_lnRn}.\\

\bigskip 
\noindent \textbf{Proof of Theorem~\ref{th:1Rn}(i): } Choose some arbitrary constants $\rho<\gamma<1$ and let $L$ be a L\'evy process with its generating triplet $(\sigma^2,\nu,\Gamma)$. Let $\nu_\eta$ and $\nu_\xi$ be L\'evy measures, such that $\nu=\nu_\eta+\nu_\xi$, for which Hypotheses~I and II are satisfied. 
We prove the theorem in three subsequent settings of increasing generality: 
\begin{enumerate}
	\item[\textbf{(a)}] $\sigma^2=0$, $\Gamma = 0$, $\nu_\xi=0$ and $\nu_\eta(\bR^n)=1$ (the pure compound Poisson case).
	\item[\textbf{(b)}] $\sigma^2=0$, $\Gamma = 0$, $\nu_\xi=0$ and $\nu_\eta(\bR^n)<\infty$. 
	\item[\textbf{(c)}] No further restrictions besides Hypotheses~I and II.
\end{enumerate}

\bigskip
\noindent\textbf{Proof of \eqref{mu.tRn} in setting (a)}
\noindent We consider the solution $g_o$ of the auxiliary modified functional equation 
\begin{align}
&g_o(\Lambda)f'(g_o(\Lambda))-f(g_o(\Lambda))+\ln g_o(\Lambda)-\frack n2\ln f''(g_o(\Lambda))+\frack{ng_o(\Lambda)f'''(g_o(\Lambda))}{2f''(g_o(\Lambda))}\nonumber\\
&\hspace{8cm}+\frack n2\ln(2\pi)+\frack{n-1}2\ln(\alpha-1)=\ln \Lambda.\label{Def_g0Rn}
\end{align}
Existence, uniqueness and the properties of $g_o$ follow by Lemma~\ref{g_lnRn}. \\

\noindent\textbf{Claim 1: } For any $\delta>0$ and $\rho<\gamma<1$, there exists a constant $k>0$ sufficiently large, such that for every $|x|>k$ and every $t\in[|x|^\rho,|x|^\gamma]$ the following estimate is valid:
\ba
-|x|\big(f'(g_o(\frack {|x|}t))+&(n(\frack\alpha2-1)-1+\delta)g_o(\frack {|x|}t)^{-1}\big)\spac{7pt}\le\ln\mu_t(|x|)\\
&\spac{7pt}\le\sup\limits_{s\le t}\ln\mu_s(|x|)\spac{7pt}\le-|x|\big(f'(g_o(\frack {|x|}t))+(n(\frack\alpha2-1)-1-\delta)g_o(\frack {|x|}t)^{-1}\big).\label{431Rn}
\ea
By Lemma~\ref{SR}$(ii)$ we have 
\begin{equation}\label{e:tripleprimef}
\lim_{\Lambda\rightarrow\infty}\frack{ng_o(\Lambda)f'''(g_o(\Lambda))}{2f''(g_o(\Lambda))}=n(\frack\alpha2-1).
\end{equation}
A comparison of the definition \eqref{Def_gRn} of $g$ with the definition \eqref{Def_g0Rn} of $g_o$ 
combined with the limit~\refq{e:tripleprimef} and the uniqueness property in Lemma ~\ref{g_lnRn}$(iv)$ yields that \refq{431Rn} implies the desired estimate~\eqref{mu.tRn} in setting (a). This establishes the main statement \eqref{mu.tRn} in Theorem~\ref{th:1Rn}(i) in setting~(a). It is therefore enough to show Claim 1 in the sequel.\\ 

\noindent\textbf{Proof of Claim 1, lower bound:} For $n\in\bN$, $x\in\bR^n$ and $t\in[|x|^\rho,|x|^\gamma]$ let $m_{x,t}:=|x|g_o(\frac{|x|}t)^{-1}$ and $\tilde m_{x,t}:=\lfloor m_{x,t}\rfloor+1$. Since $\rho <\gamma<1$ we have by the definition of $m_{x,t}$ and Lemma \ref{g_lnRn}$(i)$ that
$$\lim\limits_{|x|\rightarrow\infty}\inf\limits_{t\in[|x|^\rho,|x|^\gamma]}\frac{|x|}{\tilde m_{x,t}}=\infty\qquad\mbox{and}\qquad\lim\limits_{|x|\rightarrow\infty}\sup\limits_{t\in[|x|^\rho,|x|^\gamma]}\frac{t}{\tilde m_{x,t}}=0.$$
The first limit in the preceding display implies that the lower bound of \refq{Faltung_FormelRn} of Proposition~\ref{l:2Rn} can be applied to estimate $\nu^{*\tilde m_{x,t}}(dx)/dx$. The second limit is used in the inequality of \eqref{327}. 
We use equation~\refq{eq:7} and recall that $N_t$ has a Poisson distribution with expectation~$t$ (since we are in setting (a)). For $\delta\in(0,1)$, $|x|$~can be chosen sufficiently large such that for all $t\in[|x|^\rho,|x|^\gamma]$ the following estimate holds:
\begin{align}
&\ln\mu_t(x)\ge\ln\Big(\P(N_t=\tilde m_{x,t})\frac{\nu^{*\tilde m_{x,t}}(dx)}{dx}\Big) \nonumber\\
&=-t+\ln\Big(\frac{t^{\tilde m_{x,t}}}{\tilde m_{x,t}!}\nu^{*\tilde m_{x,t}}(dx)\Big)\nonumber\\
&\ge- m_{x,t}\Big(\ln\frac { m_{x,t}}{t}-1+f\Big(\frac{x}{ m_{x,t}}\Big)+\frack n2\ln f''\Big(\frac{x}{ m_{x,t}}\Big)-\frack{n-1}{2}\ln(\alpha-1)-\frack n2\ln2\pi+\frack\delta2\Big)\label{327}\\
&=-|x|g_o(\frack{|x|}t)^{-1}\Big(\ln\frack {|x|}t-\ln g_o(\frack{|x|}t)-1+f(g_o(\frack{|x|}t))+\frack n2\ln f''(g_o(\frack{|x|}t))-\frack{n-1}{2}\ln(\alpha-1)-\frack n2\ln2\pi+\frack\delta2\Big)\nonumber\\
&=-|x|g_o(\frack{|x|}t)^{-1}\Big(g_o(\frack{|x|}t)f'(g_o(\frack{|x|}t))+\frac n2\frac{g_o(\frack{|x|}t)f'''(g_o(\frack{|x|}t))}{f''(g_o(\frack{|x|}t))}-1+\frack\delta2\Big).\label{432Rn}
\end{align}
Here we apply Stirling's formula to approximate $\tilde m_{x,t}!$ and in the last step we use the definition~\refq{Def_g0Rn} of~$g_o$. Finally by \eqref{e:tripleprimef} inserted in \eqref{432Rn} we obtain the lower bound of \refq{431Rn}.\\

\noindent \textbf{Proof of Claim 1, upper bound:} \\
Fix some $\delta\in (0,\frac13)$ and set 
\begin{equation}\label{d:sx}
s_x:=\lfloor 2(1-\gamma)^{-1}(1-\rho)^{\frac{\alpha-1}\alpha}\frac\alpha{\alpha-1}g(|x|)^{-1}|x|+1\rfloor.\end{equation} By \refq{eq:7} we have
\be
\mu_t(x)\spac{7pt}\le s_x\max_{m\le s_x}\P(N_t=m)\frac{\nu^{*m}(dx)}{dx}+\sum_{m=s_x}^\infty\P(N_t=m)\frac{\nu^{*m}(dx)}{dx}.\label{433Rn}
\ee
We start with the estimate of the first summand on the right-hand side of \eqref{433Rn}. Since $g(\Lambda) \to \infty$ as $\Lambda\to \infty$ we have $\lim\limits_{x\rightarrow\infty}\frac{|x|}{s_x}=\infty$. Therefore, Proposition~\ref{l:2Rn}, in particular the upper bound of~\refq{Faltung_FormelRn}, can be applied to estimate $\nu^{*m}(dx)/dx$ for all $m\le s_x$. Again we apply Stirling's formula to estimate the Poisson distribution $\P(N_t=m)$. Set $r:=\frac{|x|}{m}$. For $|x|$ sufficiently large we obtain
\begin{\eq}
-\ln\P(&\hspace{-9pt}N_t&\hspace{-9pt}=m)\frac{\nu^{*m}(dx)}{dx}
\spac{7pt}\ge m(\ln\frack mt-1+f(\frack{|x|}m)+\frack n2\ln f''(\frack {|x|}m)-\frack n2\ln2\pi-\frack{n-1}2\ln(\alpha-1)-\delta)\nonumber\\
&\ge&m(\ln\frack mt-1+f(\frack{|x|}m)+\frack n2\ln f''(\frack {|x|}m)-\frack n2\ln2\pi-\frack{n-1}2\ln(\alpha-1))-\delta s_x\label{434}\\
&=&\frack{|x|}{r}(\ln\frack{|x|}{t}-\ln r-1+f(r)+\frack n2\ln f''(r)-\frack n2\ln2\pi-\frack{n-1}2\ln(\alpha-1))-\delta s_x.
\label{434.Rn}
\end{\eq}
We differentiate the right-hand side of \eqref{434.Rn} with respect to $r$ and obtain 
\ba
\frack d{dr}&\frack{|x|}{r}\Big(\ln\frack {|x|}{t}-\ln r-1+f(r)+\frack n2\ln f''(r)-\frack{n-1}{2}\ln(\alpha-1)-\frack n2\ln2\pi\Big)\\
=&\frack{|x|}{r}\Big(-r^{-1}+f'(r)+\frack n2\frack{f'''(r)}{f''(r)}\Big)\\
&-\frack{|x|}{r^2}\Big(\ln\frack {|x|}{t}-\ln r-1+f(r)+\frack n2\ln f''(r)-\frack{n-1}{2}\ln(\alpha-1)-\frack n2\ln2\pi\Big)\\
=&\frack{|x|}{r^2}\Big(rf'(r)-f(r)+\ln r+\frack n2\frack{rf'''(r)}{f''(r)}-\frack n2\ln f''(r)+\frack{n-1}{2}\ln(\alpha-1)+\frack n2\ln2\pi-\ln\frack {|x|}{t}\Big).\label{354}
\ea
From the definition of $g_o$ in \eqref{Def_g0Rn} it follows that this derivative has a unique zero at $r_{x,t}=g_o(\frac{|x|}t)$. In both cases, $r\rightarrow0$ and $r\rightarrow\infty$, the right-hand side of \refq{434.Rn} tends to $\infty$. Therefore, $r_{x,t}=g_o(\frac{|x|}t)$ is the minimizer on the right-hand side of \refq{434.Rn}. Recall that $m=\frac{|x|}r$. Thus, by construction it follows that $m_{x,t}$ is the minimizer on the right-hand side of \refq{434}. The corresponding value of the minimum is obtained by a term-by-term comparison of \refq{327} with \refq{434.Rn} combined with the identity \refq{432Rn}. From the definitions of $s_x$ and $m_{x,t}$ together with Lemma~\ref{g_lnRn}$(i)$, it follows that $k$ can be chosen sufficiently large, such that $s_x\delta<km_{x,t}\delta$ for all $t\in[|x|^\rho,|x|^\gamma]$. We obtain
\ba
-\ln\max_{m\le s_x}\frack{t^m}{m!}\frac{\nu^{*m}(dx)}{dx}\spac{7pt}\ge& |x|(f'(g_o(\frack{|x|}t))+g_o(\frack{|x|}t)^{-1}(n(\frack\alpha2-1)-1-2k\delta))
\ea
again for $|x|$ sufficiently large and all $t\in[|x|^\rho,|x|^\gamma]$. Therefore the first summand in \eqref{433Rn} satisfies the claimed upper bound.\\

\noindent We continue with the second summand on the right-hand side of \eqref{433Rn}. Note that in case of $m>s_x$ (in contrast to $m\le s_x$) we cannot assume $\frac{|x|}m$ to be large enough to use the upper bound of \refq{Faltung_FormelRn}. Instead, we use that $\nu^{*m}(dx)/dx$ is uniformly bounded. For every $\kappa>0$ we have $\sup\limits_{|y|>\kappa,m\in\bN}\nu^{*m}(dy)/dy<\infty$, moreover, for $|x|$ sufficiently large we have \mbox{$\sup\limits_{|y|>|x|,m\in\bN}\nu^{*m}(dy)/dy<1$.} Thus
\be
\sup_{t\in[|x|^\rho,|x|^\gamma]}\sum_{m=s_x}^\infty\P(N_t=m)\frac{\nu^{*m}(dx)}{dx} 
\le\sum_{m=s_x}^\infty\frac{(|x|^\gamma)^m}{m!}\spac{8pt}\le\frac{(|x|^\gamma)^{s_x}}{s_x!}\sum_{m=0}^\infty\Big(\frac{|x|^\gamma}{s_x}\Big)^m 
\le 2\frac{(|x|^\gamma)^{s_x}}{s_x!}.\label{438Rn}
\ee
By Lemma~\ref{g_lnRn}$(i)$, the limit \refq{e:lim2} of Lemma~\ref{g_lnRn}$(iii)$, and Lemma~\ref{SR}$(vi)$ and $(viii)$ we have
\be
\lim_{\Lambda\rightarrow\infty}\frac{\sup_{t\in[\Lambda^\rho,\Lambda^\gamma]}g_o(\Lambda)f'(g_o(\frac{\Lambda}{t}))}{\ln \Lambda}\spac{8pt}=\lim_{\Lambda\rightarrow\infty}\frac{g_o(\Lambda)f'(g_o(\Lambda))}{\ln \Lambda}\cdot\frac{f'(g_o(\Lambda^{1-\rho}))}{f'(g_o(\Lambda))}=\frack\alpha{\alpha-1}(1-\rho)^{\frac{\alpha-1}{\alpha}}.\label{437Rn}
\ee
By \refq{438Rn}, \refq{437Rn} and the definition of $s_x$ in \refq{d:sx} we obtain that $|x|$ can be chosen sufficiently large, such that 
\begin{align}
-\ln\sup_{t\in[|x|^\rho,|x|^\gamma]}\sum_{m=s_x}^\infty\P(N_t=m)\frac{\nu^{*m}(dx)}{dx}
&\ge s_x\big(\ln\frack{s_x}{|x|^\gamma}-2\big) 
\ge(1-\delta)(1-\gamma)s_x\ln|x|\nonumber\\
&\ge(2-3\delta)\sup_{t\in[|x|^\rho,|x|^\gamma]}f'(g_o(\frack{|x|}t))|x|\label{e:771}.
\end{align}
Recall that by the initial choice of $\delta<\frac13$ we have $2-3\delta>1$. Thus, by the limits of \eqref{eq:1.12b} it follows that the term on the right-hand side of \refq{e:771} is sufficiently small to satisfy the upper bound of~\refq{431Rn}. This completes the proof of \hbox{Theorem~\ref{th:1Rn}} in the setting~(a).\\

\noindent \textbf{Proof of \refq{mu.tRn} in setting (b):} 
Let $a:=\nu(\bR^n)$. Consider the L\'evy process $\tilde L_t:=L_{a^{-1}t}$. This process has the jump measure 
\be\label{eq:4.44}
\tilde\nu(dx)=\exp(-\tilde f(|x|)) dx\mb{2mm}{with}\tilde f:=f+\ln a \in SR_{\alpha},\mb{2mm}{thus}\tilde\nu(\bR)=\int_o^\infty\exp(-\tilde f(y))dy=1.
\ee
Therefore $\tilde L$ falls under the setting (a) and its marginal density $\tilde \mu_t(x)$ can be estimated by \refq{mu.tRn}, where $g$ is replaced by $\tilde g$, and $\tilde g$ is defined by \eqref{Def_gRn} with $f$ being replaced by $\tilde f$. Note that the respective derivatives of $f$ and $\tilde f$ coincide.\\
By definition of $\tilde L$ we have $\mu_t(x)=\tilde\mu_{at}(x)$. Therefore $\mu_t(x)$ can be estimated by \refq{mu.tRn} with $g(\frac{|x|}{t})$ being replaced by $\tilde g(\frac{|x|}{at})$. A comparison of the respective defining equations, \refq{Def_gRn} for $g$ and the correspondingly adjusted equation for $\tilde g$, yields that $g(\frac{|x|}t)=\tilde g(\frac{|x|}{at})$. This completes the proof of \hbox{Theorem~\ref{th:1Rn}}(i) in the setting (b).\\

\noindent \textbf{Proof of \refq{mu.tRn} in setting (c):} Let $\mu_{\eta,t}$ (and $\mu_{\xi,t}$, resp.) denote the density of the distribution of $\eta_t$ (and $\xi_t$, resp.). In the proof of setting (b) it is shown that estimate \refq{mu.tRn} holds for $\mu_t$ being replaced by $\mu_{\eta,t}$. By Hypothesis~II we know that \refq{BEDIIRn} holds and by construction we know, that $\mu_t$ equals the convolution of the distributions of $\eta_t$ and $\xi_t$.\\

\noindent \textbf{Proof of the lower bound of \refq{mu.tRn} in setting (c):} Let $\vartheta\in(\gamma,1)$. 
By the convolution density formula $\mu_t$ has the obvious lower bounds 
\ba
\mu_t(x) 
\spac{7pt}= (\mu_{\xi, t} * \mu_{\eta,t})(x)\spac{7pt}= \int_{\bR^n} \mu_{\xi, t}(y) \mu_{\eta,t}(x-y) dy\spac{7pt}\geq \P(|\xi_t|<|x|^\vartheta) \inf_{|y-x|\le|x|^\vartheta}\mu_{\eta,t}(y)
\ea
and hence 
\be
-\ln \mu_t(x) \le-\ln\Big(\P(|\xi_t|<|x|^\vartheta)\inf_{|y-x|\le|x|^\vartheta}\mu_{\eta,t}(y)\Big).\label{438bRn}
\ee
The choice of $\vartheta$ yields that 
$\lim_{x\rightarrow\infty}\inf_{t<|x|^\gamma}\P(|\xi_t|<|x|^\vartheta)=1$. By the proof of setting~(b) and Remark~\ref{remRn}, item~2. for $\Lambda = |x|$ and $z_\Lambda = |x|^\vartheta$, the lower bound of \refq{mu.tRn} can be applied to estimate $\inf_{|y-x|\le|x|^\vartheta}\mu_{\eta,t}(y)$ for $|x|$ large enough. 
Inserting those estimates in \refq{438bRn} we obtain that $|x|$ can be chosen sufficiently large, such that $\mu_t(x)$ satisfies the lower bound of \refq{mu.tRn} for all $t\in[|x|^\rho,|x|^\gamma]$.\\

\noindent \textbf{Proof of the upper bound of \refq{mu.tRn} in setting (c):} Choose $c\in(1-\frac1\alpha,1)$, such that \refq{BEDIIRn} in Hypothesis~II is satisfied. Let $\kappa\in(\gamma,1)$, $p\in(0,c+\frac1\alpha-1)$ and $q(x):=|x|(\ln|x|)^{-p}$. We have 
\be
\mu_t(x) \le\sup_{|y|\le|x|^\kappa}\mu_{\eta_t}(x-y)+(2q(x))^n\sup_{|y|\in[|x|^\kappa,q(x)]}\mu_{\eta,t}(x-y)\mu_{\xi,t}(y)+\sup_{|y|\ge q(x)}\mu_{\xi,t}(y).\label{439Rn}
\ee
We estimate the terms on the right-hand side of \eqref{439Rn} one by one. 
We start with the first term. By case (b) the choice of $\kappa<1$ combined with Remark~\ref{remRn}, item~2., the upper bound of~\refq{mu.tRn} is valid for $\sup_{|y|\le|x|^k}\mu_{\eta_t}(x-y)$. 

\noindent We continue with the third summand on the right-hand side of \eqref{439Rn}. 
By~\refq{BEDIIRn}, the definition of $q(\cdot)$ and the choice of $p$, in particular $c-p>1-\frack1\alpha$, it follows that the third summand $\sup_{|y|\ge q(x)}\mu_{\xi,t}(y)$ in \eqref{439Rn} is sufficiently small to satisfy the upper bound of~\refq{mu.tRn}.

\noindent To estimate the second term on the right side of \refq{439Rn} we apply the already proven case~(b) to estimate \hbox{$\mu_{\eta_t}(x-y)$} and then apply \refq{BEDIIRn} to estimate $\mu_{\xi,t}(y)$. For any $\delta>0$ we can choose $|x|$ sufficiently large, such that the following estimate holds uniformly for all $t\in[|x|^\rho,|x|^\gamma]$.
\begin{align}
\hspace{-5mm}-\sup_{|y|\in[|x|^\kappa,q(x)]}\ln\mu_{\eta,t}(x-y)\mu_{\xi,t}(y)
&=\inf_{|y|\in[|x|^\kappa,q(x)]}\Big(|x-y| \tilde h_\delta(\ln\frac{|x-y|}{t}))+|y|(\ln|y|)^c\Big),\label{453Rn}
\end{align}
where $\tilde h_\delta(r):=f'(g(\exp(r)))-(1+\delta)g(\exp(r))^{-1}$. By Lemma~\ref{g_lnRn}(i) and Lemma~\ref{SR}(iii), (vi), (vii) and (viii) we obtain that $\tilde h_\delta\in SR_{1-\frac1\alpha}$. Thus for $|x|$ sufficiently large and any $t\in[|x|^\rho,|x|^\gamma]$ we have for $|y|\in[q(x),|x|^\kappa]$ that the first summand on the right-hand side of \eqref{453Rn} is monotonically growing as a function of $|x-y|$ while the second summand is monotonically growing as a function of $|y|$. Consequently the desired minimum is obtained for $y \in [[0,x]]$, where we have $|x-y|=|x|-|y|$. Set $\Lambda=|y|$. With $\tilde h_\delta\in SR_{1-\frac1\alpha}$, thus $\tilde h_\delta'\in SR_{-\frac1\alpha}$, and for $|x|$ sufficiently large we obtain

\begin{\eq}
\frack{d}{d\Lambda}\Big((|x|-\Lambda)\tilde h_\delta(\ln\frack{|x|-\Lambda}{t})+\Lambda(\ln \Lambda)^c\Big)&=&-\tilde h_\delta(\ln\frack{|x|-\Lambda}{t})-\tilde h_\delta'(\ln\frack{|x|-\Lambda}{t})+(\ln \Lambda)^c+c(\ln \Lambda)^{c-1}\nonumber\\
&\ge&-2\tilde h_\delta(\ln\frack{|x|-\Lambda}{t})+\kappa^c(\ln|x|)^c\spac{10pt}>0\label{343}
\end{\eq}
for all $t\in[|x|^\rho,|x|^\gamma]$ and $\Lambda\in[|x|^k,q(x)]$. The positivity in the last step follows by $\tilde h_\delta\in SR_{1-\frac1\alpha}$ together with the choice of $c>1-\frac1\alpha$. Now, having identified the minimizer of the right-hand side of \refq{453Rn} as $\Lambda=|y|=|x|^\kappa$ for $y\in[[0,x]]$, we can finally estimate the second summand of the right-hand side of \refq{439Rn}:
\begin{\eq}
-\sup_{|y|\in[|x|^\kappa,q(x)]}\ln(2q(x))^n\mu_{\eta,t}(x-y)\mu_{\xi,t}(y)&\ge&\inf_{\Lambda\in[|x|^\kappa,q(x)]}\Big((|x|-\Lambda)\tilde h_\delta(\ln\frack{|x|-\Lambda}{t})+\Lambda(\ln \Lambda)^c\Big)-n\ln|x|\nonumber\\
&\ge&(|x|-|x|^\kappa)\tilde h_\delta(\ln\frack{|x|-|x|^\kappa}{t}).
\end{\eq}
Finally , by the definition of $\tilde h_\delta$ and Remark~\ref{remRn}, item~2., we have that the second summand of~\refq{439Rn} also satisfies the upper bound of \refq{mu.tRn}. This completes the proof of Theorem~\ref{th:1Rn}(i) in setting~(c) and hence under the hypotheses in full generality.\vspace{2mm}

\noindent \textbf{Proof of Theorem~\ref{th:1Rn}(ii): } 
We show the conditional jump size probability estimate \eqref{e:cjspe}. Fix some $c\in(\alpha^{-1},1)$ and $\delta\in(0,\alpha c-1)$. For $|x|$ sufficiently large let $\xi^x$ and $\eta^x$ be L\'evy processes with generating triplets $(\sigma^2,\nu_{\xi,x},\Gamma)$, and $(0,\nu_{\eta,x},0)$, respectively, 
where 
$$\nu_{\xi,x}:=\nu_{|_{\{y\in\bR^n~:~|y|\le|x|^c\}}}\qquad \mbox{ and }\qquad \nu_{\eta,x}:=\nu-\nu_{\xi,x}.$$ 
For $t>0$, recall that $\mu_t$ (and $\mu_{\xi^x,t}$) denotes the density of the distribution of $L_t$ (and $\xi^x_t$, respectively). Let $N^x_t$ be the number of jumps of $\eta^x$ in the time interval $(0,t]$ and $W^x_i$ denote the 
$i$-th jump of $\eta^x$. Furthermore, for $m=1,2,\dots $ and $t>0$ denote by $\mu^{x,m}_t$ the density of the distribution of $\xi^x_t+\sum_{i=1}^mW_i^x$, and $\mu^{x,o}_t=\mu_{\xi^x,t}$. We have
\begin{\eq}
\P\Big(\sup_{s\le t}|L_s-L_{s-}|>|x|^c\Big|L_t=x\Big)&=&\Big(\sum_{m=1}^\infty\P(N^x_t=m)\mu^{x,m}(x)\Big)\mu_t(x)^{-1}\nonumber\\
&\le&\P(N^x_t>0)\Big(\sup_{m\ge1}\mu^{x,m}(x)\Big)\mu_t(x)^{-1}.\label{441Rn}
\end{\eq}
By the asymptotic of the Poisson distribution of $N^x_t$ and the lower bound of Theorem~\ref{th:1Rn}(i) we know that $|x|$ can be chosen sufficiently large, such that for the first and the third factor on the right-hand side of \refq{441Rn} we have
\be
\P(N^x_t>0)\spac{8pt}\le t\nu_{\eta,x}(\bR^n)\spac{8pt}\le\exp(-|x|^{\alpha c-\frac\delta2})\mb{8mm}{and}\mu_t(x)\spac{8pt}\ge\exp(-|x|\ln|x|)\label{442Rn}
\ee
is valid for all $t\in[|x|^\rho,|x|^\gamma]$. To estimate the remain	ing second factor on the right-hand side of~\refq{441Rn}, we denote by $\tilde\mu^x$ the density of the distribution of each single jump of~$\eta^x$. By definition, for $m=1,2,\dots $ we have $\mu^{x,m}=\mu^{x,m-1}*\tilde\mu^x$. Thus by $f\in SR_\alpha$ we obtain
\ba
\sup_{m\ge1}\mu^{x,m}(x)&\spac{6pt}\le\sup_{y\in\bR^n}\tilde\mu^x(y)\spac{6pt}= \exp(-f(|x|^c))\Big( \int_{|x|^c}^\infty e^{-f(r)} dr\Big)^{-1}\le|x|^k\label{443Rn}
\ea
for a suitably chosen exponent $k>0$ and $|x|$ sufficiently large. Finally we obtain the assertion $(ii)$ by a combination of \refq{441Rn}, \refq{442Rn} and \refq{443Rn}. This completes the proof of \hbox{Theorem~\ref{th:1Rn}.}\vspace{3mm}

\noindent\textbf{Proof of Corollary~\ref{corRn}:} By definition we have $\P(|L_t|>\Lambda)=\int_{\bR^n}\mu_t(x)\textbf{1}_{[\Lambda,\infty)}(|x|)dx$. 
Thus the upper and the lower bounds of~\refq{mu.tRn} yield a bound on $\P(|L_t|>\Lambda)$, where $|x|$ is replaced by $\Lambda$. 
Obviously the inequalities 
$$
\P(|L_t|>\Lambda)\spac{7pt}\le\P\Big(\sup_{s\le t}|L_s|>\Lambda\Big)\spac{7pt}\le\P\Big(\sup_{0\le s_o\le s_1\le t}|L_{s_1}-L_{s_o}|>\Lambda\Big)
$$ 
are satisfied. First we establish an upper bound for $\P(\sup_{s\le t}|L_s|>\Lambda)$. This is used later on to obtain the desired upper bound of $\P(\sup_{0\le s_o\le s_1\le t}|L_{s_1}-L_{s_o}|>\Lambda)$. 
Note that 
\be
\P\Big(|L_t|>(1-(\ln\Lambda)^{-2})\Lambda\Big|\sup_{s\le t}|L_s|>\Lambda\Big)\P\Big(\sup_{s\le t}|L_s|>\Lambda\Big)\le
\P(|L_t|>(1-(\ln\Lambda)^{-2})\Lambda).\nonumber
\ee
Thus by the Markov property and an adapted reflection principle for $L$ we obtain 
\begin{align*}
\P\Big(\sup_{s\le t}|L_s|>\Lambda\Big)&\le
\P\Big(|L_t|>(1-(\ln\Lambda)^{-2})\Lambda\Big) \P\Big(|L_t|>(1-(\ln\Lambda)^{-2})\Lambda\Big|\sup_{s\le t}|L_s|>\Lambda\Big)^{-1}\nonumber\\
&\le\P\Big(|L_t|>(1-(\ln\Lambda)^{-2})\Lambda\Big)\Big(\inf_{s\le t}\P(|L_s|\le(\ln\Lambda)^{-2}\Lambda)\Big)^{-1}.
\end{align*}
By assumption we have $t\le\Lambda^\gamma$ for $\gamma<1$. Thus, for $\Lambda$ sufficiently large, the second factor is smaller than 2. As noted above we may use the upper bound of~\refq{mu.tRn} to estimate $\P(|L_t|>\Lambda)$ and 
by Remark~\ref{remRn},~item~2., the same upper bound holds for $\P(|L_t|>(1-(\ln(\Lambda)^{-2})\Lambda)$. 
Thus the upper bound of \refq{mu.tRn} remains valid if $\sup_{s\le t}\mu_t(x)$ is replaced by $\P(\sup_{s\le t}|L_s|>\Lambda)$ and $|x|$ by $\Lambda$. That is to say, for $\Lambda$ sufficiently large we have 
\begin{equation}\label{keineAhnung}
\ln \P\Big(\sup_{s\le t} |L_s| > \Lambda\Big) \le - \Lambda \Big(f'(g(\frack{\Lambda}{t})) - (1+ \delta) g(\frack{\Lambda}{t})^{-1}\Big). 
\end{equation}
We now establish the upper bound for $\P(\sup_{0\le s_o\le s_1\le t}|L_{s_1}-L_{s_o}|>|x|)$. Set $k_x=n\lfloor 1+\ln|x|\rfloor^2$, $\tilde D_x=\{y=(y_1,\dots,y_n)\mid|y_i|\in\{0,\frac 1{k_x},\frac 2{k_x},\dots,1\}\}$ and \hbox{$D_x=\{\frac y{|y|}\mid y\in\tilde D_x\setminus\{0\}\}$.} By definition we have $|y|=1$ for every $y\in D_x$ and 
\be
\inf_{z\neq 0}\max_{y\in D_x}\langle y,\frac{z}{|z|}\rangle\spac{8pt}\ge1-\sup_{|z|=1}\min_{y\in D_x}|y-z|\spac{8pt}\ge1-(\ln|x|)^{-2}.
\ee
Thus by a simple union bound we estimate
\ba
\P\Big(\sup_{0\le s_o\le s_1\le t}|L_{s_1}-L_{s_o}|>&|x|\Big)\spac{5pt}=	\P\Big(\sup_{0\le s_o\le s_1\le t}\,\sup_{|y|=1}\langle y,L_{s_1}-L_{s_o}\rangle>|x|\Big)\\
\spac{5pt}\le&|D_x|\max_{y\in D_x}\P\Big(\,\sup_{0\le s_o\le s_1\le t}\langle y,L_{s_1}-L_{s_o}\rangle>(1-(\ln|x|)^{-2})|x|\Big).\label{362}
\ea
For $y\in D_x$ and $k=0,1, \dots, k_x-2$ we define
\begin{align*}
\varphi_{x,y}&:=\inf\{s>0\mid\sup_{0<s_o<s}\langle y,L_s-L_{s_o}\rangle\ge(1-(\ln|x|)^{-2})|x|\}\\ 
\tau_k^{x,y}&:=\inf\{s>0\mid\langle y,L_s\rangle\le-\frack{k}{k_x}|x|\}.
 \end{align*}
 Using these stopping times and the strong Markov property of $L$ we continue to estimate the probability on the right-hand side of \refq{362} 
\ba
&\P\Big(\,\sup_{0\le s_o\le s_1\le t}\langle y,L_{s_1}-L_{s_o}\rangle>(1-(\ln|x|)^{-2})|x|\Big)=\P(\varphi_{x,y}\le t)\\
&\le\sum_{k=0}^{k_x-2}\P\Big(\varphi_{x,y}\le t,\inf_{s<\varphi_{x,y}}\langle y,L_s\rangle\in[-\frack{k+1}{k_x}x,-\frack k{k_x}x]\Big)+\P\Big(\inf_{s<t}\langle y,L_s\rangle\le(1-\frack 1{k_x})x\Big)\\
&\le\sum_{k=0}^{k_x-2}\P\Big(\tau_k^{x,y}\le t,\sup_{s<t}\langle y,L_{\tau_k^{x,y}+s}-L_{\tau_k^{x,y}}\rangle\ge(1-\frack2{k_x})x\Big)+\P\Big(\sup_{s<t}|L_s|\le(1-\frack 1{k_x})x\Big)\\
&\le k_x\P\Big(\sup_{s<t}|L_s|\ge(1-\frack2{k_x})x\Big).\label{362b}
\ea
For $|x|$ sufficiently large by construction we have $|D_x|\le(3\ln|x|)^n\le|x|$ and $k_x\le(\ln|x|)^3$. Again, by the already proven upper bound of $\P(\sup_{s\leq t} |L_s| > |x|)$ \eqref{keineAhnung} in combination with Remark~\ref{remRn}, item~2., we obtain that the upper bound of \refq{mu.tRn} holds for $\P(\sup_{s\leq t} |L_s| > (1- 2 |(\ln |x|)^{-2}) |x|)$. Clearly, the coefficients $k_x$ and $|D_x|$ are sufficiently small such that by \eqref{362b} inserted in \eqref{362} we obtain the desired upper bound for $\P(\sup_{0\le s_o\le s_1\le t}|L_{s_1}-L_{s_o}|>|x|)$. This completes the proof of Corollary~\ref{corRn}.

\bigskip 

\subsection{Proof of the LDP for the L\'evy bridge $Y^\e$ (Theorem~\ref{Th:2Rn})}\label{ss:LDP}

\noindent In order to prove Theorem~\ref{Th:2Rn} we combine the following Proposition~\ref{l:43Rn}, Lemma~\ref{l:44Rn} and the representation \eqref{eq:2} of Theorem~4.28 in \cite{Feng}. 

\begin{prop}[Finite dimensional LDP]\label{l:43Rn}
Fix $\x\in \mathbb{R}^n$, $\x\neq 0$, let $T>0$ and $(Y^\e)_{\e>0}$ and $S$ be defined as in Theorem~\ref{Th:2Rn}. For any $m\in\bN$ and $\tau=(t_1,\dots,t_m)$ with $0<t_1<\dots<t_m<T$ the family $(\Law(Y^\e_{t_1},\dots,Y^\e_{t_m}))_{\e>0}$ satisfies a LDP \hbox{on $((\bR^n)^m,|\,\cdot\,|)$} with the speed function $S$ and the good rate function
\begin{equation}\label{finiteentropy}
I_{\x, \tau}(y_1,\dots,y_m)=\left\{
\begin{array}{ll}
\sum\limits_{i=1}^{m+1} |y_i-y_{i-1}| \ln\frac{|y_i-y_{i-1}|}{t_i-t_{i-1}}-|\x|\ln\frac{|\x|}{T}&\mbox{if $y_i\in[[0,\x]]$, $|y_1|\le\dots\le|y_m|$,}\\[3mm]
\infty&\mbox{otherwise}.
\end{array}\right.
\end{equation}
Here we denote $y_o=t_o=0$, $y_{m+1}=\x$, $t_{m+1}=T$ and whenever $y_i=y_{i-1}$ let \hbox{$|y_i-y_{i-1}|\ln\frac{|y_i-y_{i-1}|}{t_i-t_{i-1}}$} be defined as $0$.
\end{prop}

\begin{rem}\label{rem:convI}
Note that for any $x>0$ we have $\inf_{\varphi}\int_o^T\varphi'(t)\ln\varphi'(t)=x\ln\frac xT$, where the infimum is taken over all continuous non-decreasing functions $\varphi:[0,T]\to\bR$ with $\varphi(0)=0$ and $\varphi(T)=x$. With this identity a comparison of the rate functions $I_{\x}$ and $I_{\x,\tau}$ yields 
\be
I_{\x,(t_1,\dots,t_m)}(y_1,\dots,y_m)=\inf\{I_{\x}(\varphi)\mid\varphi\in\mathcal D_{[0,T],\bR^n},\varphi(t_i)=y_i, i=1,\dots,m\}.
\ee
\end{rem}

\begin{lemma}[$S$-exponential tightness]\label{l:44Rn}
Let $T>0$ and $(Y^\e)_{\e>0}$ and $S$ be defined as in Theorem~\ref{Th:2Rn}. Let $\mathcal J_1$  denote the Skorohod topology on $\mathcal D_{[0,T],\bR^n}$. Then $(Y^\e)_{\e>0}$ is S-exponentially tight in $(\mathcal D_{[0,T],\bR^n},\mathcal J_1)$.
\end{lemma}

\noindent Before the proof of Proposition~\ref{l:43Rn} and Lemma~\ref{l:44Rn} show that those results imply Theorem~\ref{Th:2Rn}.\vspace{2mm}

\noindent \textbf{Proof of Theorem~\ref{Th:2Rn}: } In order to apply Theorem~4.28 in J. Feng, T. G. Kurtz \cite{Feng} we define $\tilde Y^\e_t=Y^\e_{t\wedge T}$. By construction, for $0<t_1<\dots<t_m<T$ this process satisfies the same LDPs formulated in Proposition~\ref{l:43Rn}. Furthermore we know that $\tilde Y^\e_o=0$ and $\tilde Y^\e_t=\x$ for any $t\ge T$ almost surely. Thus, for $0=t_1<\dots<t_m<T$ the family $(\tilde Y^\e_{t_1},\dots,\tilde Y^\e_{t_m})$ satisfies the LDP with the speed function $S$ and the rate function 
\be\label{390a}
I_{\x, \tau}(y_1,\dots,y_m):=\left\{
\begin{array}{ll}
I_{\x, (t_2,\dots,t_m)}(y_2,\dots,y_m)&\mbox{if }y_1=0\\
\infty&\mbox{if }y_1\neq0
\end{array}\right.
\ee
and for $0\le t_1<\dots<t_m$ let $m_o=\max(\{0\}\cup\{i=1,\dots,m|t_i<T\})$. The family $(\tilde Y^\e_{t_1},\dots,\tilde Y^\e_{t_m})$ satisfies the LDP with the speed function $S$ and the rate function
\be\label{391a}
I_{\x, \tau}(y_1,\dots,y_m):=\left\{
\begin{array}{ll}
I_{\x, (t_1,\dots,t_{m_o})}(y_1,\dots,y_{m_o})&\mbox{if }y_{m_o+1}=y_{m_o+2}=\dots=y_{m}=\x\\
\infty&\mbox{otherwise},
\end{array}\right.
\ee
for $m_o\ge1$. For $m_o=0$ let $I_{\x, \tau}(y)=0$ for $y_i=\x$ and $I_{\x, \tau}(y)=\infty$ otherwise.

\noindent By construction we have that $S$-exponential tightness of $(Y^\e)_{\e>0}$ on $(\mathcal D_{[0,T],\bR^n},\mathcal J_1)$, which is given by Lemma \ref{l:44Rn}, directly implies $S$-exponential tightness of $(\tilde Y^\e)_{\e>0}$ on $(\mathcal D_{[0,\infty),\bR^n},\mathcal J_1)$. Now, by Theorem~4.28 in J. Feng, T. G. Kurtz \cite{Feng} we obtain that $(\tilde Y^\e)_{\e>0}$ satisfies the LDP on $(\mathcal D_{[0,\infty),\bR^n},\mathcal J_1)$ with the speed function $S$.  With Remark \ref{rem:convI} the comparison between the rate functions $I_{\x}$ and $I_{\x,\tau}$  yields
\be\label{392}
\hspace{-7mm}\tilde I_{\x}(\varphi)=\hspace{-5pt}\sup_{\substack{0\le t_1<\dots<t_m\\\mbox{$\varphi$ continuous in $t_1,\dots,t_m$}}}\hspace{-5pt}I_{\x,(t_1,\dots,t_m)}(\varphi(t_1),\dots,\varphi(t_m))=
\left\{
\begin{array}{ll}
I_{\x}\big(\varphi_{\big|[0,T]}\big)\hspace{5mm}&\mbox{for $\varphi_{\big|[T,\infty)}\equiv\x$}\\
\infty&\mbox{otherwise,}
\end{array}\right.
\ee
where $I_{\x}$ is the rate function on $\mathcal D_{[0,T],\bR^n}$ stated in Theorem~\ref{Th:2Rn}.

Note that the rate function $\tilde I_{\x}$ assigns the value $\tilde I_{\x}(\varphi)=\infty$ to every discontinuous function $\varphi$. By Theorem~4.13 and Theorem~4.14 in Feng and Kurtz \cite{Feng}, we get that the family $(Y^\e)_{\e>0}$ is $S$-exponentially tight also with respect to the supremum norm. It follows that the large deviation principle therefore holds on $(\mathcal D_{[0,\infty),\bR^n},\|\cdot\|_\infty)$ with the same speed function $S$ and rate function $\tilde I_{\x}$. 

\noindent By the contraction principle the continuous embedding $(\mathcal D_{[0,T],\bR^n},\|\cdot\|_\infty)\hookrightarrow(\mathcal D_{[0,\infty),\bR^n},\|\cdot\|_\infty)$ finaly yields the LDP for $(Y^\e)_{\e>0}$ on $(\mathcal D_{[0,T],\bR^n},\|\cdot\|_\infty)$ with the speed function $S$ and rate function $I_{\x}$ given in \refq{entropyrate} in Theorem~\ref{Th:2Rn}. 
Therefore we have shown that Lemma 3.1 and Proposition \ref{l:43Rn} imply Theorem \ref{Th:2Rn}. 

\noindent For later use note that the procedure of applying \refq{eq:2} upon a process $Y^\e_{\cdot\wedge T}$ on $\mathcal D_{[0,\infty),\bR^n}$ and afterwards transfering the LDP upon $\mathcal D_{[0,T],\bR^n}$ via contraction principle yields a rate function with $I_{\x}(\varphi)=\sup_{0<t1<\dots<t_m<T}I_{\x,(t_1,\dots,t_m)}(\varphi(t_1),\dots,\varphi(t_m))$ for every continuous function $\varphi$.

\noindent The next lemma, which shall be proven first, 
is the essential tool for the proofs of Proposition~\ref{l:43Rn} and Lemma~\ref{l:44Rn}.
We show that asymptotically the paths of $Y^\e$ which fall outside the continuous parametrizations of the segment $[[0, \x]]$ or which have not increasing norms are $S$-exponentially negligible.

\begin{lemma}[$S$-exponential negligibility of interval breakouts and decreasing paths]\label{l:45Rn}\hfill\\
\noindent Let $\x\in \bR^n\setminus\{0\}$, $T>0$ and $Y^\e$ and $S$ be defined as in Theorem~\ref{Th:2Rn}. Then for any $\kappa>0$ we have
\be\label{e:negl}
\lime S(\e)\ln\Big(\P\Big(\sup_{0\le t\le T}\inf_{y\in[[0,\x]]}|Y^\e_t-y|\ge\kappa\Big)+\P\Big(\sup_{0\le s<t\le T}|Y^\e_s|-|Y^\e_t|\ge\kappa\Big)\Big)\spac{8pt}=-\infty.
\ee
\end{lemma}

\noindent \textbf{Proof of Lemma~\ref{l:45Rn}: } 
To estimate the first summand let $A\subset\bR^n$ be a closed subset, such that 
\begin{equation}\label{e:geo}
\inf_{y\in A}|y|+\inf_{y\in A}|\x-y|>|\x|
\end{equation} 
and let $\sigma_A:=\inf\{t>0~|~\e L_t\in A\}$. Let $k_1:=\inf_{y\in A}|y|$ and $k_2:=\inf_{y\in A}|\x-y|$. 
By the strong Markov property of $L$, the definition of the bridge density of $Y^\e$ we have 
\ba\label{457Rn}
\P(\exists t\in[0,T]:Y^\e_t\in A)&
\spac{5pt}=\P(\sigma_A \le r_\e T~|~\e L_{r_\e T} = \x)\\[2mm]
&\spac{5pt}\le \P(\sigma_A \leq r_\e T) \sup_{t < T, y \in A} \mu_{r_\e (T-t)}(\e^{-1}(\x-y)) \mu_{r_\e T} (\e^{-1} \x)^{-1}\\
&\spac{5pt}\le \P(\sup_{t\le r_\e T} |L_t| \ge k_1 \e^{-1}) \sup_{t < T, |z| \ge k_2} \mu_{r_\e t}(\e^{-1} z) 
\mu_{r_\e T} (\e^{-1} \x)^{-1}.
\ea
By Theorem~\ref{th:1Rn} and Corollary~\ref{corRn} each of the terms $\P(\sup_{t\le r_\e T}|\e L_t|\ge k_1)$, $\sup_{t<T, z\in A}\mu_{r_\e t}(\e^{-1}z)$ and $\mu_{r_\e T}(\e^{-1}\x)$ can be estimated as in inequality~\eqref{mu.tRn} from above and from below. Combining the mentioned estimates with the definition of $k_1$ and $k_2$, and \eqref{e:geo} we obtain 
\be\label{458Rn}
\lime f'(g(\frack{\e^{-1}}{r_\e}))^{-1}\e\ln\P(\exists t\in[0,T]:Y^\e_t\in A)\spac{8pt}\le|\x|-k_1-k_2\spac{8pt}<0.
\ee
Note that by definition of $S(\e)$ it is easy to see that $\lim_{\e\to 0 } S(\e) (\e f'(g(\frack{\e^{-1}}{r_\e})))^{-1} = \infty$. Thus 
\be\label{459Rn}
\lime S(\e)\ln\P(\exists t\in[0,T]:Y^\e_t\in A)\spac{5pt}=-\infty.
\ee
Let $M$ be a finite collection of sets, such that each $A\in M$ satisfies \eqref{e:geo} and $\{y\in\bR^n\mid\inf_{z\in[[0,\x]]}|y-z|\ge\kappa\}\subseteq\bigcup_{A\in M}A$. 
For instance, think of $M = \{A_o, A_1, \dots, A_\ell\}$ for some $\ell < \infty$ with $A_o := \{y\in \mathbb{R}^n~|~|y|\ge 2|\x|\}$ and $\{A_1, \dots, A_\ell\}$ is a finite cover of the precompact set $\{y\in \mathbb{R}^n~|~\inf_{z\in [[0, \x]]} |y-z| \geq \kappa\}\setminus A_o$, where each of set $A_i$, $i=1, \dots, \ell$ satisfies 
\eqref{e:geo}.
Then we have
\be
\hspace{-5mm}\lime S(\e)\ln\P\Big(\sup_{0\le t\le T}\inf_{y\in[[0,\x]]}|Y^\e_t-y|\ge\kappa\Big)\le\max_{A\in M}\lime S(\e)\ln\P(\exists\, t\in[0,T]:Y^\e_t\in A)\spac{2pt}=-\infty.\label{460Rn}
\ee
We continue to estimate the second probability in the statement \eqref{e:negl}. Let $K:=\lfloor 3|\x|\kappa^{-1}\rfloor+1$ and $\vartheta:=\inf\{t>0\mid|\e L_t|\le\sup_{s\in[0,t]}|\e L_s|-\kappa\}$. We have
\ba\label{461Rn}
&\P\Big(\sup_{0\le s<t\le T}|Y^\e_t|-|Y^\e_s|\ge\kappa\Big)\\
\spac{8pt}\le&\P\Big(\sup_{0\le t\le T}\inf_{y\in[[0,\x]]}|Y^\e_t-y|\ge\frack\kappa3\Big)+\sum_{i=0}^K\P\Big(\vartheta<r_\e T,\e|L_\vartheta|\in[\frack i3\kappa,\frack{i+1}{3}\kappa]\Big|\e L_{r_\e T}=\x\Big).
\ea
By \refq{460Rn} it follows that the first summand on the right side is sufficiently small. To estimate the remaining sum we argue similarly to \refq{457Rn} with the help of the strong Markov property
\ba
\P\Big(&\vartheta<r_\e T,\e|L_\vartheta|\in[\frack i3\kappa,\frack{i+1}3\kappa]~\Big|~\e L_{r_\e T}=\x\Big)\\
\spac{8pt}\le& \P\Big(\vartheta < r_\e T, \e |L_{\vartheta}|\in [\frack{i}{3}\kappa, \frack{i+1}{3} \kappa]\Big)\Big(
\sup_{\substack{|y|\le (\frac{i+1}{3})\kappa\\t< T}} \mu_{(T-t)r_\e} (\e^{-1} (\x -y))\Big) \mu_{T r_\e}(\e^{-1} \x)^{-1}\\
\spac{8pt}\le&\P\Big(\sup_{t<T}\e|L_{r_\e t}|\ge(\frack{i}3+1)\kappa\Big)\Big(\sup_{\substack{|z| \ge  |\x| -\frac{i+1}{3} \kappa\\s< T}} \mu_{sr_\e} (\e^{-1} z)\Big) \mu_{T r_\e}(\e^{-1} \x)^{-1}.
\ea
Similarly to\eqref{458Rn} and \eqref{459Rn} with $\tilde k_1 = (\frac{i}{3}+1) \kappa$ and $\tilde k_2 = |\x| - \frac{i+1}{3}\kappa$ it follows that 
\be\label{kA1}
\lime S(\e)\ln\P\Big(\vartheta<r_\e T,\e|L_\vartheta|\in[\frack i3\kappa,\frack{i+1}3\kappa]\Big|\e L_{r_\e T}=\x\Big)\spac{8pt}=-\infty.
\ee
Inserting \eqref{kA1} in \eqref{461Rn} we obtain the $S$-exponential negligibility of the second term in \eqref{e:negl}. This completes the proof of Lemma~\ref{l:45Rn}.\\

\noindent \textbf{Proof of Proposition~\ref{l:43Rn}: } Let $m\in\bN$ and $\tau=(t_1,\dots,t_m)\in\bR^m$ such that $0<t_1<\dots<t_m<T$ and $\x\in \bR^n$. By Lemma~\ref{l:45Rn} it follows that $(Y^\e_{t_1},\dots,Y^\e_{t_m})_{\e>0}$ is $S$-exponentially tight on $((\bR^n)^m,|\;\cdot\;|)$. Therefore it is sufficient to show the following limit. \vspace{2mm}

\noindent \textbf{Claim 2: } For every $y\in(\bR^n)^m$ we have
\be
\lim_{\kappa\rightarrow0}\lime S(\e)\ln\P(\max_{i=1,\dots,m}|Y^\e_{t_i}-y_i|\le\kappa)=-I_{\x, \tau}(y_1,\dots,y_m).\label{eq:74Rn} \ee

\noindent The $S$-exponential tightness of $(Y^\e_{t_1},\dots,Y^\e_{t_m})_{\e>0}$ together with Claim 2 implies that the LDP stated in Proposition~\ref{l:43Rn} holds (see Lemma~1.2.18 and Theorem~4.1.11 in Dembo Zeitouni~\cite{Dembo}). \\

\noindent The proof of \refq{eq:74Rn} in Claim 2 is carried out in three different settings:
\begin{itemize}
	\item[(A)] \textbf{Unordered norms or exceptional cases: } 
	$y\in(\bR^n)^m$, such that \textbf{either} there is some $i\in\{1,\dots,m\}$ with $y_i\notin[[0,\x]]$ \textbf{or} $y$ does not satisfy $|y_1|\le\dots\le|y_m|\le|\x|$.
	\item[(B)] \textbf{Strictly ordered, positive norms: } $y\in(\bR^n)^m$  satisfies $y_i\in[[0,\x]]$ and $0<|y_1|<\dots<|y_m|<|\x|$.
	\item[(C)] \textbf{Ordered, but not strictly ordered norms: } $y\in(\bR^n)^m$  satisfies $y_i\in[[0,\x]]$ and $|y_1|\le\dots\le|y_m|\le|\x|$, \textbf{but not} $0<|y_1|<\dots<|y_m|<|\x|$.
\end{itemize}
On the one hand we use Lemma \ref{l:44Rn} in the proof of Claim 2 setting (C). On the other hand, we use Claim 2 settings (A) and (B) in the proof of Lemma \ref{l:44Rn}. Thus the proof will be carried out in the following order: Proof of Claim 2 settings (A) and (B), proof of Lemma \ref{l:44Rn}, and finally proof of Claim~2 setting (C).\vspace{2mm}

\noindent \textbf{Proof of Claim 2 in setting (A):} In setting (A), by the choice of $y$ in combination with Lemma~\ref{l:45Rn} it follows that $\lim_{\kappa \to 0}\lim_{\e \to 0} S(\e) \ln \P(\max_{i=1, \dots, m} |Y^\e_{t_i} - y_i| \leq \kappa) = -\infty$. By the definition of $I_{\x, \tau}$ we have $I_{\x, \tau}(y_1, \dots, y_m) = \infty$ for any $y$ from setting (A). Thus, the limit \refq{eq:74Rn} from Claim 2 is satisfied.\\

\noindent \textbf{Proof of Claim 2 in setting (B):} 
Let $\tau=(t_1,\dots,t_m)$ with $0<t_1<\dots<t_m<T$. We consider $y\in(\bR^n)^m$ in setting (B).  
 Remember the notation $y_o=t_o=0$, $y_{m+1}=\x$ and $t_{m+1}=T$. By the choice of $y$ there is $\kappa>0$ sufficiently small, such that $\kappa<\frac13\min\{|y_{i+1}-y_i|\mid i=0,1,\dots,m\}$. \\

\noindent \textbf{Proof of the upper bound: } Let $C_{y,\kappa}:=\{z\in(\bR^n)^m\mid\max_{i=1,\dots,m}|y_i-z_i|<\kappa\}$. By construction we have 
\be\label{465Rn}
\P\Big(\max_{i=1,\dots,m}|Y^\e_{t_i}-y_i|\le\kappa\Big)\spac{8pt}=\int_{C_y,\kappa}\e^{-mn}\mu_{r_\e T}(\e^{-1}\x)^{-1}\prod_{i=0}^m\mu_{r_\e(t_{i+1}-t_i)}(\e^{-1}(z_{i+1}-z_i))dz,
\ee
where $z_0:=0$ and $z_{m+1}:=\x$. By the choice of $\kappa$ we have $\min_{i=0,1,\dots,m+1}|z_{i+1}-z_i|>\kappa$ for all \hbox{$z\in C_{y,\kappa}$.} Thus Theorem~\ref{th:1Rn} is applicable to estimate all those densities simultaneously. Let $\lambda_{n,m}$ denote the Lebesgue measure on $(\bR^n)^m$. Note that $\lim\limits_{\e\rightarrow0}S(\e)\ln(\e^{-mn}\lambda_{n,m}(C_{x,\kappa}))=0$. Thus we obtain
\ba\label{466Rn}
&-\limsup_{\e\rightarrow0}S(\e)\ln\P\Big(\max_{i=1,\dots,m}|Y^\e_{t_i}-y_i|\le\kappa\Big)\\
\hspace{8pt}&\ge-\limsup_{\e\rightarrow0}g(\frack{\e^{-1}}{r_\e})\e\sup_{z\in C_y,\kappa}\ln\Big(\Big(\prod_{i=0}^m\mu_{r_\e(t_{i+1}-t_i)}(\e^{-1}(z_{i+1}-z_i))\Big)\mu_{r_\e T}(\e^{-1}\x)^{-1}\Big)\\
\hspace{8pt}&\ge
\limsup_{\e\rightarrow0}g(\frack{\e^{-1}}{r_\e})\inf_{z\in C_y,\kappa}\Big(\sum_{i=0}^m|z_{i+1}-z_i|(f'(g(\frack{|z_{i+1}-z_i|\e^{-1}}{r_\e(t_{i+1}-t_i)}))-(1-\delta)g(\frack{\e^{-1}}{r_\e})^{-1})\\
&\hspace{40mm}-|x|(f'(g(\frack{|x|\e^{-1}}{T r_\e}))-(1+\delta)g(\frack{\e^{-1}}{r_\e})^{-1})\Big),
\ea
where in the last step we use Theorem~\ref{th:1Rn}$(i)$ combined with the fact, that by Lemma~\ref{g_lnRn}$(i)$ and the definition of $C_{y,\kappa}$ it follows that $\lim\limits_{\e\rightarrow0} g(\frack{|z_{i+1}-z_i|\e^{-1}}{r_\e(t_{i+1}-t_i)})g(\frack{\e^{-1}}{r_\e})^{-1}=1$ uniformly for $z\in C_{y,\kappa}$. 
	
\noindent Let $z=(z_1,\dots,z_m)\in C_{y,\kappa}$ and $\tilde z=(\tilde z_1,\dots,\tilde z_m)\in(\bR^n)^m$ be defined by $\tilde z_i:=\x /|\x|^{-2	}\langle \x,z_i\rangle \in [[0, \x]]$. By construction we obtain for $i=1,2,\dots,m$ the inequality $|y_i-\tilde z_i|<|y_i-z_i|$ and hence $\tilde z\in C_{y,\kappa}$. Furthermore the construction yields $|\tilde z_{i+1}-\tilde z_i|\le|z_{i+1}-z_i|$ for $i=0,1,\dots,m$ (with $\tilde z_o:=0$ and $\tilde z_{m+1}:=\x$). Consequently, since each of the summands on the right-hand side of \eqref{466Rn} increases as a function of $|z_{i+1}-z_i|$, the term to be minimized on the right side of \refq{466Rn} becomes only smaller if $z$ is replaced by $\tilde z\in C_{y,\kappa}\cap [[0,\x]]^m$. Note that for $z\in C_{y,\kappa}\cap[[0,\x]]$ we have $\sum_{i=0}^m|z_{i+1}-z_i|=|\x|$, thus we continue \eqref{466Rn} with the help of the asymptotic cancellation property of Lemma~\ref{g_lnRn}$(ii)$ and obtain
\ba\label{467Rn}
&-\limsup_{\e\rightarrow0}S(\e)\ln\P\Big(\max_{i=1,\dots,m}|Y^\e_{t_i}-y_i|\le\kappa\Big)\\
\ge\;&2|\x|\delta+\limsup_{\e\rightarrow0}g(\frack{\e^{-1}}{r_\e})\inf_{z\in C_{y,\kappa}\cap[[0,\x]]^m}\Big(\sum_{i=0}^m|z_{i+1}-z_i|(f'(g(\frack{|z_{i+1}-z_i|\e^{-1}}{r_\e(t_{i+1}-t_i)}))-f'(g(\frack{\e^{-1}}{r_\e})))\\
&\hspace{66mm}-|x|(f'(g(\frack{|x|\e^{-1}}{r_\e T}))-f'(g(\frack{\e^{-1}}{r_\e})))\Big)\\
\ge\;&2|x|\delta+\inf_{z\in C_{y,\kappa}\cap[[0,\x]]^m}\sum_{i=0}^m|z_{i+1}-z_i|\ln\frack{|z_{i+1}-z_i|}{t_{i+1}-t_i}-|x|\ln\frac{|\x|}{T}.
\ea
Recall that the case $m=2$ was presented step-by-step in Subsection~\ref{sss:approach}.II, see formulas \refq{eq:8}-\refq{eq:123}. Sending $\delta\rightarrow0$ we finally obtain the asserted upper bound of Claim 2
\be
\lim_{\kappa\rightarrow0}\limsup_{\e\rightarrow0} -S(\e)\ln\P\Big(\max_{i=1,\dots,m}|Y^\e_{t_i}-y_i|\le\kappa\Big)\spac{8pt}\ge\lim_{\kappa\rightarrow0}\inf_{z\in C_{y,\kappa}\cap[[0,\x]]^m}I_{\x,\tau}(z)\spac{8pt}=I_{\x,\tau}(y).
\ee

\noindent \textbf{Proof of the lower bound:} Let $C_{y}:=\{z\in(\bR^n)^m\mid\max_{i=1,\dots,m}|y_i-z_i|<|\ln\e|^{-2} |\x|\}$. Clearly, we have $\e^{-mn}\lambda_{n,m}(C_y)>1$, and similarly to \refq{465Rn} we obtain
\be\label{469Rn}
\P\Big(\max_{i=1,\dots,m}|Y^\e_{t_i}-y_i|\le\kappa\Big)\spac{8pt}\ge\mu_{r_\e T}(\e^{-1}\x)^{-1}\inf_{z\in C_y}\prod_{i=0}^m\mu_{r_\e(t_{i+1}-t_i)}(\e^{-1}(z_{i+1}-z_i)).
\ee
We apply Theorem~\ref{th:1Rn}$(i)$, to estimate the densities on the right side of \refq{469Rn}. By the construction of $C_y$ together with the robustness result in Remark~\ref{remRn}, item~2. it follows that for every $\delta>0$ we can choose $\e$ sufficiently small, such that the estimate
\ba\label{470Rn}
&-\ln\prod_{i=0}^m\mu_{r_\e(t_{i+1}-t_i)}(\e^{-1}(z_{i+1}-z_i))\\
&\qquad \le \sum_{i=0}^m|y_{i+1}-y_i|\e^{-1}(f'(g(\frack{(y_{i+1}-y_i)\e^{-1}}{(t_{i+1}-t_i)r_\e}))-(1-\delta)g(\frack{(y_{i+1}-y_i)\e^{-1}}{(t_{i+1}-t_i)r_\e})^{-1})
\ea
is valid for all $z\in C_y$. By the definition of setting (B) we have $y\in[[0,x]]$. Therefore the same type of arguments used in \refq{467Rn} can be applied to obtain the lower bound. This completes the proof of~\eqref{eq:74Rn} in setting (B).\\

\noindent As announced at the beginning of the proof of Proposition~\ref{l:43Rn} we continue with the proof of Lemma~\ref{l:44Rn} before concluding the proof of~\eqref{eq:74Rn} in setting (C).\\ 

\noindent \textbf{Proof of Lemma~\ref{l:44Rn} ($S$-exponential tightness)  }
By Theorem~1 and 3 in the appendix of Yu.~V.~Prokhorov \cite{Prokhorov} it is well-known that $(\mathcal D_{[0,T],\bR^n},\mathcal J_1)$ is a complete separable metric space. Thus Lemma~3.3 in Feng and Kurtz \cite{Feng} is applicable in the prove of the $S$-exponential tightness on $(\mathcal D_{[0,T],\bR^n},\mathcal J_1)$. There it is shown that it is sufficient to show the following: For every $M,\kappa>0$ there exists a compact set $\mathcal{K}_{M,\kappa}\subset\mathcal D_{[0,T],\bR^n}$, such that
\be
\limsup_{\e>0}S(\e)\ln\P\Big(\inf_{\varphi\in \mathcal{K}_{M,\kappa}}d(\varphi,Y^\e)>\kappa\Big)\spac{9pt}<-M,\label{471Rn}
\ee
where $d$ denotes the metric that induces the $\mathcal J_1$ topology. Again by Theorem~1 in the appendix of Yu.~V.~Prokhorov~\cite{Prokhorov} it is known that, if the sets $\mathcal{K}_{M,\kappa}$ consist only of continuous functions, then it is sufficient to show
\be
\limsup_{\e>0}S(\e)\ln\P\Big(\inf_{\varphi\in \mathcal{K}_{M,\kappa}}\sup_{t\in[0,T]}|\varphi-Y^\e|>\kappa\Big)\spac{9pt}<-M,\label{471aRn}
\ee instead of \refq{471Rn}. For $c>1$ and $\kappa\in[0,\frac12\wedge T)$  we set
\ba
A_{c,\kappa}\spac{3pt}{:=}\{\varphi\in&\mathcal D_{[0,T],\bR^n}\mid\varphi(0)=0,\\
&\forall s,t\in[0,T],\delta\in(\kappa,\infty):|s-t|<\delta\exp(-c\delta^{-1})\Rightarrow|\varphi(s)-\varphi(t)|<\delta\}.\nonumber
\ea
Note that the function $(\delta\mapsto\delta\exp(-c\delta^{-1}))$ is monotonically increasing on $(0,\infty)$. Therefore the definition of $A_{c,\kappa}$ does not change if we restrict the choice of $\delta$ to $\delta\in(0,\delta_c]$ with $\delta_c$ sufficiently large, such that $\delta_c\exp(-c\delta_c^{-1})\ge T$. For $\kappa=0$ we omit the parameter and set $A_c:=A_{c,0}$. Compactness criteria on $(\mathcal D_{[0,T],\bR^n},\mathcal J_1)$ are given for example in subsection 2.7 of A.V. Skorokhod~\cite{Skorokhod}. By construction those criterias are immediately satisfied by the sets $A_c$ for any $c>0$. Moreover, any function $\varphi\in A_c$ is continuous.\vspace{2mm}

\noindent Choose $c>1$ and $\kappa\in(0,\frac12\wedge T)$ arbitrary, and let $\ell:=\lfloor T\kappa^{-1}\exp(c\kappa^{-1})\rfloor+1$. Fix some $\varphi\in A_{c,\kappa}$ and let $\phi$ be defined as a linear interpolation of $\varphi$ with the following points $\phi(k\frac{T}{\ell})=\varphi(k\frac{T}{\ell})$ for $k=0,1,\dots,\ell$. By the choice of $\kappa < T$ we have $\ell\ge\lfloor T\kappa^{-1}\rfloor+1\ge2$. Thus by construction the distance between those supporting points is given by $\frac{ T}{\ell}\in[\frac12\kappa\exp(-c\kappa^{-1}),\kappa\exp(-c\kappa^{-1})]$. We show the following two statements: 
\begin{align}
&\sup_{t\in[0,T]}|\phi(t)-\varphi(t)|\le2\kappa,\label{(i)}\\
&\phi\in A_{3c}.\label{(ii)}
\end{align}
We note that \eqref{(i)} and \eqref{(ii)} imply $\{\varphi\in\mathcal D_{[0,T],\bR^n}\mid\inf_{\phi\in A_{3c}}||\phi-\varphi||_\infty<2\kappa\}\subseteq A_{3c,2\kappa}$.
Hence showing the follow-up statement 
\begin{equation}\label{(iii)}
\lim_{c\rightarrow\infty}\limsup_{\e\rightarrow0}S(\e)\ln\P(Y^\e\notin A_{c,\kappa})=-\infty \quad \mbox{ for all } \kappa>0
\end{equation}
yields that for every $M,\kappa>0$ there is $c$ sufficiently large, such that \refq{471aRn} is satisfied with the choice $\mathcal{K}_{M,\kappa}=A_c$. Recall that for each $c>0$ the elements of $A_{c}$ are continuous. Therefore, it is sufficient to proof \eqref{(i)}, \eqref{(ii)} and \eqref{(iii)} in order to obtain the desired $S$-exponential tightness of~$(Y^\e)_{\e>0}$.\\

\noindent \textbf{Proof of \eqref{(i)}:} Choose $t\in[0,T]$ arbitrary, and $k\in\{1,2,\dots,\ell\}$ such that \hbox{$t\in[(k-1)\frac{T}{\ell},k\frac{T}{\ell}]$.} Recall $\frac T\ell<\kappa\exp(-\kappa ^{-1})$. Therefore $\varphi\in A_{c,\kappa}$ implies $|\varphi(t_1)-\varphi(t_2)|\le\kappa$ for any $t_1,t_2\in[(k-1)\frac{T}{\ell},k\frac{T}{\ell}]$ and by the construction of $\phi$ we have $|\varphi(k\frac{T}{\ell})-\phi(t)|\le|\varphi(k\frac{T}{\ell})-\varphi((k-1)\frac{T}{\ell})|$. Consequently
\be
|\phi(t)-\varphi(t)|\spac{6pt}\le|\phi(t)-\varphi(k\frack{T}{\ell})|+|\varphi(k\frack{T}{\ell})-\varphi(t)| 
\spac{6pt}\le|\varphi(k\frack{T}{\ell})-\varphi((k-1)\frack{T}{\ell})|+\kappa \spac{6pt}\le 2\kappa.
\ee
This shows \eqref{(i)}.\vspace{2mm} 

\noindent \textbf{Proof of \eqref{(ii)}:} Choose $\delta>0$ arbitrary and $s,t\in[0,T]$ such, that $|t-s|\le\delta\exp(-3c\delta^{-1})$. We show that $|\phi(t)-\phi(s)|\le\delta$.\vspace{1mm}

\noindent We consider first the case $\delta\le2\kappa$: In this case we have 
\ba
\frac{|\phi(t)-\phi(s)|}{|t-s|}\spac{7pt}\le&\max\Big\{\frac{|\varphi(k\frac{T}{\ell})-\varphi((k-1)\frac{T}{\ell})|}{T\ell^{-1}}\Big|k=1,\dots,\ell\Big\}\\
\spac{7pt}\le&\frac{\kappa}{\frac12\kappa\exp(-c\kappa^{-1})}\spac{9pt}=2\exp(c\kappa^{-1}).\label{373}
\ea
By the choices $c\ge1$, $\kappa\le\frac12$ and $\delta\le2\kappa$ we have $3c\delta^{-1}-c\kappa^{-1}\ge\frac12c\kappa^{-1}\ge c> 1$ and the selection of~$s,t$ together with \refq{373} yields 
\be
|\phi(t)-\phi(s)|\spac{5pt}\le2\exp(c\kappa^{-1})|t-s|\spac{5pt}\le2\delta\exp(-3c\delta^{-1}+c\kappa^{-1})\spac{5pt}\le2\delta\exp(-1)\spac{5pt}\le\delta, 
\ee
as required. \\

\noindent We continue with the remaining case $\delta>2\kappa$: Choose $s_1,t_1\in\{k\frac{T}{\ell}|k=0,1,\dots,\ell-1\}$, let $s_2=s_1+\frac{T}{\ell}$, and $t_2=t_1+\frac{T}{\ell}$, such that $s\in[s_1,s_2]$ and $t\in[t_1,t_2]$. For $i,j\in\{1,2\}$ we get
\ba\label{eq:87Rn}
|s_i-t_j|\spac{8pt}\le&|s-t|+2\kappa\exp(-c\kappa^{-1})\spac{8pt}\le\delta\exp(-3c\delta^{-1})+2\kappa\exp(-c\kappa^{-1})\\
\spac{8pt}\le&\delta\exp(-c\delta^{-1})\Big(\exp(-2c\delta^{-1})+2\kappa\delta^{-1}\exp(-c(\kappa^{-1}-\delta^{-1}))\Big)\\
\spac{8pt}\le&\delta\exp(-c\delta^{-1}).
\ea

\noindent To justify the last step in this estimate we use, that $c>1$, $\kappa<\frac12$ and thus
$$\exp(-2c\delta^{-1})+2\kappa\delta^{-1}\exp(-c(\kappa^{-1}-\delta^{-1}))\le\left\{
\begin{array}{ll}
\exp(-1)+\exp(-1)&\mbox{for }\delta\in[2\kappa,4\kappa]\\[2mm]
\exp(-1)+\frac12&\mbox{for }\delta\in[4\kappa,2c]\\[2mm]
(1-c\delta^{-1})+2\kappa\delta^{-1}&\mbox{for }\delta	\ge2c
\end{array}\right\}<1,$$
where in the case $\delta\ge2c$ we have used that $e^{-2x} \le 1-x$ for $0 < x \le \frac{1}{2}$. By~\refq{eq:87Rn} and $\varphi\in A_{c,\kappa}$ with $\kappa<\delta$ we obtain that $|\varphi(s_i)-\varphi(t_j)|\le\delta$ for every choice $i,j\in\{1,2\}$. By construction 
$|\phi(s)-\phi(t)|\le\max_{i, j\in \{1,2\}} |\varphi(s_i)-\varphi(t_j)| < \delta$ as desired. This completes the proof of~\eqref{(ii)}.\\[3mm]
\noindent \textbf{Proof of \eqref{(iii)}:} For $c>1$ and $\kappa\in(0,\frac12)$ set 
\[
\delta_{c}:=\inf\{\delta>0\mid \delta\exp(-c\delta^{-1})>T\}, \quad \mbox{ and }
\quad m_{c,\kappa}:=\lfloor\frack{\delta_{c}}{\kappa}\rfloor+1.
\] 
By the definition it can easily be seen that $\delta_c\sim\frac c{\ln c}$. For our further estimates it will be sufficient to notice that $\delta_c<\frac{2c}{\ln c}$, thus $m_{c,\kappa}\le\frac{3c}{\kappa\ln c}$ for $c$ sufficiently large. This follows by montonicity of $(\delta\mapsto\delta\exp(-c\delta^{-1}))$ and $\delta\exp(-c\delta^{-1})=2\frac{\sqrt c}{\ln c}>T$ for the choice $\delta=\frac{2c}{\ln c}$. Furthermore, for $i=0,1,\dots,m_{c,\kappa}$ set 
\[
s_{c,\kappa,i}:=2(i+1)\kappa\exp\Big(-\frac{c}{(i+1)\kappa}\Big) \quad \mbox{ and }\quad \ell_{c,\kappa,i}:=\lfloor\frac{2T}{s_{c,\kappa,i}}\rfloor+2.
\] 
For $i=0,1,\dots,m_{c,\kappa}$ and $j=0,1,\dots,\ell_{c,\kappa,i}$ set 
\[
t_{i,j}:=\frac{j}{\ell_{c,\kappa,i}}T\quad \mbox{ and }\quad q_{i,j}:=(t_{i,j}+s_{i,c})\wedge T.
\] 
Then we have
\begin{\eq}
&&\hspace{-7mm}(\mathcal D_{[0,T],\bR^n}\setminus A_{c,\kappa})\cap\{\varphi\in\mathcal D_{[0,T],\bR^n}\mid\varphi(0)=0\}\nonumber\\
\hspace{-3mm}&\hspace{-3mm}=&\{\varphi\in\mathcal D_{[0,T],\bR^n}\mid\exists \delta\in(\kappa,\delta_c], s,t\in[0,T]:|t-s|<\delta\exp(-c\delta^{-1}),|\varphi(t)-\varphi(s)|\ge\delta\}\nonumber\\
\hspace{-3mm}&\hspace{-3mm}\subseteq&\bigcup_{i=1}^{m_{c,\kappa}}\{\varphi\in\mathcal D_{[0,T],\bR^n}\mid\exists\delta\in[\kappa,(i+1)\kappa], s,t\in[0,T]:|t-s|<\delta\exp(-c\delta^{-1}),|\varphi(t)-\varphi(s)|\ge \delta\}\nonumber\\
\hspace{-3mm}&\hspace{-3mm}\subseteq&\bigcup_{i=1}^{m_{c,\kappa}}\{\varphi\in\mathcal D_{[0,T],\bR^n}\mid\exists s,t\in[0,T]:|t-s|<(i+1)\kappa\exp(-c((i+1)\kappa)^{-1}),|\varphi(t)-\varphi(s)|\ge i\kappa\}\nonumber\\
\hspace{-3mm}&\hspace{-3mm}\subseteq&\bigcup_{i=1}^{m_{c,\kappa}}\bigcup_{j=1}^{\ell_{c,\kappa,i}}\{\varphi\in\mathcal D_{[0,T],\bR^n}\mid\exists s,t\in[t_{i,j},q_{i,j}]:|\varphi(t)-\varphi(s)|\ge i\kappa\}\nonumber\\
\hspace{-3mm}&\hspace{-3mm}\subseteq&
\bigcup_{i=1}^{m_{c,\kappa}}\bigcup_{j=1}^{\ell_{c,\kappa,i}}\{\varphi\in\mathcal D_{[0,T],\bR}\mid|\varphi(q_{i,j})-\varphi(t_{i,j})|>(i-\frack12)\kappa\}\label{476Rn}\\
\hspace{-3mm}&\hspace{-3mm}&\cup\Big\{\varphi\in\mathcal D_{[0,T],\bR^n}\mid\sup_{0\le s<t\le T}|\varphi(s)|-|\varphi(t)|>\frack\kappa8\Big\}\cup\Big\{\varphi\in\mathcal D_{[0,T],\bR^n}\mid\sup_{0\le t\le T}\inf_{y\in[[0,\x]]}|\varphi(t)-y|>\frack\kappa8\Big\}.\nonumber
\end{\eq}
Note that the last inclusion follows from the following reasoning. For $s,t\in[t_{i,j},q_{i,j}]$ we have 
\ba
|\phi(s)-\phi(t)|\spac{7pt}\le&||\phi(s)|-|\phi(t)||+2\sup_{0\le t\le T}\inf_{y\in[[0,\x]]}|\varphi(t)-y|\hspace{15mm}\mbox{and}\\
||\phi(s)|-|\phi(t)||\spac{7pt}\le&||\phi(q_{i,j})|-|\phi(t_{i,j})||+2\sup_{0\le s<t\le T}(|\varphi(s)|-|\varphi(t)|).
\ea

\noindent By definition we have $\P(Y^\e_o\neq0)=0$. Lemma~\ref{l:45Rn} yields that $\P(\sup_{0\le s<t\le T}|\varphi(s)|-|\varphi(t)|>\frack\kappa8)$ and $\P(\sup_{0\le t\le T}\inf_{y\in[[0,x]]}|\varphi(t)-y|>\frack\kappa8)$ are $S$-exponentially negligibly small. Thus by 
\refq{476Rn} we obtain
\be\label{477Rn}
\limsup_{\e\rightarrow0}S(\e)\ln\P(Y^\e\notin A_{c,\kappa})\spac{8pt}\le\sup_{\substack{t\in[0,T]\\i=1,\dots,m_{c,\kappa}}}\limsup_{\e\rightarrow0}S(\e)\ln\P\Big(|Y^\e_{(t+s_{c,\kappa,i})\wedge T}-Y^\e_t|>(i-\frack12)\kappa\Big).
\ee
Let $0<t_1<t_2<T$. As a direct consequence of  Lemma~\ref{l:45Rn} we have, that $(Y^\e_{t_1},Y^\e_{t_2})_{\e>0}$ is $S$-exponentially tight on $((\bR^n)^2,|\,\cdot\,|)$. Thus $(Y^\e_{t_1},Y^\e_{t_2})_{\e>0}$ satisfies the upper bound of a LDP with speed function $S$ and a rate function
\be\label{380}
I_{\x,(t_1,t_2)}^+(y_1,y_2)\spac{8pt}{:=}-\lim_{\kappa\rightarrow0}\limsup_{\e\rightarrow0}S(\e)\ln\P\Big(\max_{i=1,2}|Y^\e_{t_i}-y_i|<\kappa\Big).
\ee 
For any $k>0$ set $C_k:=\{(y_1,y_2)\in(\bR^n)^2\mid|y_1-y_2|>k\}$. We have already shown that \refq{eq:74Rn} is valid in the settings (A) and (B). Therefore, on $C_k$ the rate function $I^+_{(t_1,t_2)}$ equals the rate function $I_{\x,(t_1,t_2)}$ defined in Proposition~\ref{l:43Rn}. We get
\ba\label{479Rn}
-\limsup_{\e\rightarrow0}&S(\e)\ln\P(|Y^\e_{t_2}-Y^\e_{t_1}|>k)\spac{7pt}=-\limsup_{\e\rightarrow0}S(\e)\ln\P((Y^\e_{t_2}, Y^\e_{t_1})\in C_k)\\
&\spac{7pt}\ge\inf_{y\in C_k}I^+_{\x,(t_1,t_2)}(y_1,y_2) \spac{7pt}=\inf_{y\in C_k}I_{\x,(t_1,t_2)}(y_1,y_2)\\
&\spac{7pt}=\inf_{\substack{y\in C_k\cap[[0,\x]]^2\\|y_1|<|y_2|}}\Big(|y_1|\ln\frack{|y_1|}{t_1}+|y_2-y_1|\ln\frack{|y_2-y_1|}{t_2-t_1}+|\x-y_2|\ln\frack{|\x-y_2|}{T-t_2}-|\x|\ln\frack{|\x|}T\Big)\\
&\spac{7pt}\ge k\ln\frack{k}{t_2-t_1}-\frack{T}{e}-|\x|\ln\frack{|\x|}{T}.
\ea
For the last step we use that for any $t>0$ we have $\inf_{z>0}z\ln\frac zt=-\frac{t}{e}$ to estimate term by term 
\[
|y_1|\ln\frack{|y_1|}{t_1}+|\x-y_2|\ln\frack{|\x-y_2|}{T-t_2}{\spac{7pt}\ge-\frack{t_1+(T-t_2)}{e}}\spac{7pt}>-\frack{T}{e}.
\] 
In order to estimate the limits on the right hand side of \refq{477Rn} it remains to verify \refq{479Rn} for the remaining cases $0=t_1<t_2<T$ and $0<t_1<t_2=T$. By the construction of $Y^\e$, in the case $0=t_1<t_2<T$ we have $\P(|Y^\e_{t_2}-Y^\e_{t_1}|>k)=\P(|Y^\e_{t_2}|>k)$ and in the case $0<t_1<t_2=T$ we have $\P(|Y^\e_{t_2}-Y^\e_{t_1}|>k)=\P(|\x-Y^\e_{t_1}|>k)$. For $t\in(0,T)$ we define the rate functions $I^+_{\x,t}$ on $\bR^n$ analogously to $I^+_{\x,(t_1,t_2)}$ in \refq{380}. Similar reasoning as in \refq{479Rn} with $I^+_{\x,t_1}$ respectively $I^+_{\x,t_2}$ instead of $I^+_{\x,(t_1,t_2)}$ leads to the same upper bound as in the case $0<t_1<t_2<T$.

\noindent Finally we apply \refq{479Rn} to estimate the right side of \refq{477Rn}. By the definition of $s_{c,\kappa,i}$ we obtain for $k = (i- \frac{1}{2}) \kappa$ in each term 
\begin{\eq}\label{382}
&&-\lim_{c\rightarrow\infty}\sup_{\substack{t\in[0,T]\\i=1,\dots,m_{c,\kappa}}}\limsup_{\e\rightarrow0}S(\e)\ln\P\Big(|Y^\e_{(t+s_{c,\kappa,i})\wedge T}-Y^\e_t|>(i-\frack12)\kappa\Big)\nonumber\\
&&\ge\lim_{c\rightarrow\infty}\inf_{i=1,\dots,m_{c,\kappa}}(i-\frack 12)\kappa\ln\frack{(i-\frac12)\kappa}{s_{c,\kappa,i}}-\frack{T}{e}-|\x|\ln\frack{|\x|}{T}\spac{12pt}=\infty.
\end{\eq}
where in the last step we use the estimate
$$(i-\frack 12)\kappa\ln\frack{(i-\frac12)\kappa}{s_{c,\kappa,i}}\spac{7pt}=(i-\frack12)\kappa\Big(\frack c{(i+1)\kappa}+\ln\frack{i-\frac12}{2(i+1)}\Big)\spac{7pt}\ge\frack c4-i\kappa\ln8 $$
with $i\le m_{c,\kappa}<\frac{3c}{\kappa\ln c}$.
This completes the proof of \eqref{(iii)} and hence of Lemma~\ref{l:44Rn}.\\ 

\noindent As announced we finally show the limiting relation \refq{eq:74Rn} of Claim 2 in setting (C) with the help of the preceding Lemma~\ref{l:44Rn}.\\ 

\noindent \textbf{Proof of Claim 2 in setting (C), lower bound:} The lower bound of \refq{eq:74Rn} in Claim 2 is a consequence of the already proven settings (A) and (B) together with  the well-known lower semicontinuity of the function
\be
I_{\x, \tau}^-(y_1,\dots,y_m)\spac{8pt}{:=}-\lim_{\kappa\rightarrow0}\liminf_{\e\rightarrow0}S(\e)\ln\P\Big(\max_{i=1,\dots,m}|Y^\e_{t_i}-y_i|<\kappa\Big).
\ee 
Indeed, let $M_m$ denote the set of vectors $y\in(\bR^n)^m$ which belong to setting (A) or (B). We have already shown that $I^-_{\x, \tau}$ and $I_{\x, \tau}$ coincide on $M_m$. For $y$ from setting (C), it follows from the definition of~$I_{\x, \tau}$ in~\refq{finiteentropy} that $\lim_{\substack{z\to y\\z\in M_m}}I_{\x, \tau}(z)=I_{\x, \tau}(y)$. Together with the semicontinuity of~$I^-_{\x, \tau}$, this implies
\be
I_{\x, \tau}(y)\spac{8pt}=\liminf_{\substack{z\to y\\z\in M_m}}I_{\x, \tau}(z)\spac{8pt}=\liminf_{\substack{z\to y\\z\in M_m}}I^-_{\x, \tau}(z)\spac{8pt}\ge\liminf_{z\to y}I^-_{\x, \tau}(z)\spac{8pt}\ge I^-_{\x, \tau}(y),
\ee
which, together with the definition of $I^-_{\x,\tau}$, yields the lower bound of \refq{eq:74Rn}.\\[2mm]

\noindent \textbf{Proof of the upper bound: } By Lemma~\ref{l:45Rn} we know, that $(Y^\e)_{\e>0}$ is $S$-exponentially tight on $(\mathcal D_{[0,T],\bR^n},\mathcal J_1)$. Let $(\e_i)_{i\in\bN}$, such that $\lim_{i\rightarrow\infty}\e_i=0$. By Puhalskii \cite[Theorem~3.2.8]{Puhalskii} there exists a subsequence $(\e_{i_j})_{j\in\bN}$, such that $(Y^{\e_{i_j}})$ satisfies a LDP on $(\mathcal D_{[0,T],\bR^n},\mathcal J_1)$ with a rate function $J_{\x}$.\\

\noindent For every $m\in\bN$ and $\tau=(t_1,\dots,t_m)$ with \hbox{$0\le t_1<\dots<t_m\le T$} the contraction theorem \cite[Theorem 4.2.1]{Dembo} yields the existence of a LDP for $(Y^{\e_{i_j}}_{t_1},\dots,Y^{\e_{i_j}}_{t_m})$ with a rate function 
\be
J_{\x,\tau}(y_1,\dots,y_m)\spac{8pt}=\inf\{J_{\x}(\varphi)\,|\,\varphi(t_i)=y_i, i=1,\dots,m\}.
\ee
For $m\in\bN$, $\tau=(t_1,\dots,t_m)$ with $0<t_1<\dots<t_m<T$ and every $y\in M_m$ we already know, that $J_{\x,\tau}(y)=I_{\x,\tau}(y)$. In the sequell we use this identity to compare the rate functions $I_{\x}$ and $J_{\x}$.\\[2mm]
\noindent Let $\varphi\in \mathcal D_{[0,T],\bR^n}$ be continuous. Although formula \refq{eq:2} cannot be applied directly to the family $(Y^\e)_{\e>0} \subset \mathcal{D}_{[0,T],\mathbb{R}^n}$, the proof of Theorem \ref{Th:2Rn} outlines an indirect approach for deriving the corresponding rate function using \refq{eq:2}.
It has been shown that for every continuous function $\varphi$, this approach yields the representation of the rate function that appears in the first identity of the subsequent estimate
\ba\label{387a}
J_{\x}(\varphi)\spac{8pt}=&\sup_{0\le t_1<\dots < t_m\le T}J_{\x,(t_1,\dots,t_m)}(\varphi(t_1),\dots,\varphi(t_m))\\
\spac{8pt}\ge&\sup_{\substack{0< t_1<\dots < t_m< T\\(\varphi(t_1),\dots,\varphi(t_m))\in M_m}}J_{\x,(t_1,\dots,t_m)}(\varphi(t_1),\dots,\varphi(t_m))\\
\spac{8pt}=&\sup_{\substack{0<t_1<\dots < t_m< T\\(\varphi(t_1),\dots,\varphi(t_m))\in M_m}}I_{\x,(t_1,\dots,t_m)}(\varphi(t_1),\dots,\varphi(t_m))\\
\spac{8pt}=&\sup_{0<t_1<\dots < t_m< T}I_{\x,(t_1,\dots,t_m)}(\varphi(t_1),\dots,\varphi(t_m))\spac{8pt}=I_{\x}(\varphi).
\ea
The last identity follows similarly to the proof of Theorem \ref{Th:2Rn} by Remark \ref{rem:convI}. To justify the second to last identity for any $\varphi\in D_{\x,T}$ let $A_{\varphi}=\{t\in[0,T]||\varphi|'=0\}$. Then the integral in the definition~\refq{entropyrate} of $I_{\x}$ reads $\int_0^T|\varphi|'\ln|\varphi|'dt=\int_{[0,T]\setminus A_\varphi}|\varphi|'\ln|\varphi|'dt$.\\[2mm]
Let $\varphi\in \mathcal D_{[0,T],\bR^n}$ be discontinuous. We apply Theorem~\ref{th:1Rn}$(ii)$ to obtain $J_{\x}(\varphi)\ge I_{\x}(\varphi)$. In this case there exists an open neighborhood $A_\varphi \subset \mathcal D_{[0,T],\bR^n}$ of $\varphi$ and a constant $\Delta>0$ such that $A_\varphi\subseteq\{\phi\in\mathcal D_{[0,T],\bR^n}|\sup_{t\in [0, T]} |\phi(t)-\phi(t-)|>\Delta\}$. We obtain
\be\nonumber
J_{\x}(\varphi)=-\lim_{\kappa\to0}\lime S(\e)\ln\P(d(\varphi,Y^\e)<\kappa)\ge-\lime S(\e)\ln\P\Big(\sup_{t\in [0, T]} |Y^\e_t - Y^{\e}_{t-}|>\Delta\Big)=\infty\ge I_{\x}(\varphi),
\ee 
where $d$ denotes the metric that induces the $\mathcal J_1$ topology. Let $m\in\bN$, $\tau=(t_1,\dots,t_m)$ mit $0<t_1<\dots<t_m<T$ and $y\in(\bR^n)^m$ satisfy case (C). We get
\ba\label{484Rn}
-\lim_{\kappa\rightarrow0}\lim_{j\rightarrow\infty}S(\e_{i_j})\ln\P\Big(\max_{k=1,\dots,m}|Y^{\e_{i_j}}_{t_k}-y_k|<\kappa\Big)\spac{7pt}=&J_{\x,\tau}(y_1,\dots,y_{m})\\
\spac{7pt}=&\inf\{J_{\x}(\varphi)|\varphi(s_k)=y_k,k=1,\dots,m\}\\
\spac{7pt}\ge&\inf\{I_{\x}(\varphi)|\varphi(s_k)=y_k,k=1,\dots,m\}\\
\spac{7pt}=&I_{\x,\tau}(y_1,\dots,y_{m}). 
\ea
The second identity follows by the contraction theorem and the last identity is a direct consequence of the definitions of $I_\tau$ and $I$. Since the sequence $(\e_i)_{i\in\bN}$ has been chosen arbitrary , we obtain 
\be
-\lim_{\kappa\rightarrow0}\limsup_{\e\to0} S(\e)\ln\P\Big(\max_{k=1,\dots,m}|Y^{\e}_{t_k}-y_k|<\kappa\Big)\spac{8pt}\ge I_{\x,\tau}(y_1,\dots,y_{m})
\ee
This completes the proof of \refq{eq:74Rn} in setting (C) and the proof of Proposition~\ref{l:43Rn}.\\

\bigskip 
\subsection{Proof of the asymptotic empirical path properties of $Y^\e$ (Theorem~\ref{Th:3Rn})}
\label{ss:empirical}

\noindent By the same arguments used in \refq{eq:4.44} during the proof of Theorem~\ref{th:1Rn}$(i)$ setting (b) we may reduce the case $\nu_\eta(\bR^n)\in(0,\infty)$ to $\nu_\eta(\bR^n)=1$. For the remainder of the section we assume $\nu_\eta(\bR^n)=1$.

\subsubsection{Proof of the LDP for the jump frequency $\sqrt{S(\e)} \bar N^{\x, \e}$ (Theorem~\ref{Th:3Rn}$(i)$)}\label{sss:N} 

\noindent \textbf{Proof of Theorem~\ref{Th:3Rn}$(i)$:} By formula \eqref{eq:7} and by $\nu_\eta(\bR^n)=1$ we have
\ba\label{3104}
\P(N^{\x,\e}=m)\spac{7pt}=&\P(N_{r_\e T}=m)\frac{\nu^{*m}(dy)}{dy}_{\Big|y=\x\e^{-1}}\mu_{r_\e T}(\x\e^{-1})^{-1}
\ea
 Let $g_o$ be defined as in \eqref{Def_g0Rn} in the proof of Theorem~\ref{th:1Rn}. 
Let $\tilde m_{\x,\e}:=|\x|\e^{-1}g_o(\frac{|\x|\e^{-1}}{r_\e T})^{-1}$ and $\tilde N^{\x,\e}:=\sqrt{k_\e}(N^{\x,\e}-\tilde m_{\x,\e})$.

The proof of Theorem~\ref{Th:3Rn}$(i)$ consists of two steps: 
First we prove the validity of the asserted LDP for the family $(\sqrt{S(\e)}\tilde N^{\x,\e})_{\e>0}$ instead of  $(\sqrt{S(\e)}\bar N^{\x,\e})_{\e>0}$. Next we show, that the difference $|\tilde N^{\x,\e}-\bar N^{\x,\e}|$ is sufficiently negligible, such that we obtain the asserted LDP for $(\sqrt{S(\e)}\bar N^{\x,\e})_{\e>0}$.\\

\noindent\textbf{Proof of the LDP for $(\sqrt{S(\e)}\tilde N^{\x,\e})_{\e>0}$:} In order to prove the LDP it is sufficient to show, that
\be\label{Claim3}
\lime S(\e)\ln\P(\sqrt{S(\e)}\tilde N^{\x,\e}>M)\spac{7pt}=\lime S(\e)\ln\P(\sqrt{S(\e)}\tilde N^{\x,\e}<-M)\spac{7pt}=-\frack12M^2
\ee
for any $M>0$. Obviously the  limit \refq{Claim3} yields
\ba
\lime S(\e)\ln\P(\sqrt{S(\e)}\tilde N^{\x,\e}\notin[-M,M])\;&=\;-\frack12M^2\hspace{5mm}\mbox{ for any }M>0,\mbox{ and}\\
\lim_{\kappa\to 0}\lime S(\e)\ln\P(|\sqrt{S(\e)}\tilde N^{\x,\e}-M|<\kappa)\;&=\;-J(M) \hspace{5mm}\mbox{ for any }M\in \mathbb{R}. 
\ea
By the first limit we infer the $S$-exponential tightness of $(\sqrt{S(\e)}\tilde N^{\x,\e})_{\e>0}$ and by the second limit combined with the $S$-exponential tightness we obtain the validity of the asserted LDP (see again the combination of Lemma~1.2.18 with Theorem~4.1.11 in Dembo Zeitouni~\cite{Dembo}).\\

\noindent\textbf{Proof of \refq{Claim3}:} By the definition of $r_\e$ we may choose $\gamma<\rho<1$, such that $\e^{-\gamma}<r_\e T<\e^{-\rho}$ for~$\e$ sufficiently small. Analogously to the choice of $s_x$ in the proof of Theorem \ref{th:1Rn} set
\be\label{389}
s_{\x,\e}:=\lfloor2(1-\gamma)^{-1}(1-\rho)^{\frac{\alpha-1}\alpha}\frac\alpha{\alpha-1}g(|\x|\e^{-1})^{-1} |\x|\e^{-1}\rfloor+1.
\ee
First we show that this definition implies the following limits:
\be\label{386}
\lime S(\e)s_{\x,\e}<\infty\mb{5mm}{and}\lime S(\e)\ln\P(N^{\x,\e}>s_{x,\e})=-\infty.
\ee
The first limit of \refq{386} follows directly from the definition of $S(\e)$ and $s_{\x,\e}$. To establish the second limit, similarly to \refq{438Rn} we have
\be
\sum_{m=s_{\x,\e}}^\infty\P(N_{r_\e T}=m)\frac{\nu^{*m}(dx)}{dx}\Big|_{y=\x\e^{-1}}
\le\sum_{m=s_{\x,\e}}^\infty\frac{(|x|^\gamma)^m}{m!}\spac{8pt}
\le 2\frac{(|x|^\gamma)^{s_{\x,\e}}}{s_{\x,\e}!},
\ee
where the factorial $s_{\x,\e}!$ can be estimated by the Stirling formula as in \refq{e:771}. We obtain
\be\label{390}
\lime(f'(g(\frack{|\x|\e^{-1}}{r_\e T}))|\x|\e^{-1})^{-1}\ln\sum_{m=s_{\x,\e}}^\infty\P(N_{r_\e T}=m)\frac{\nu^{*m}(dy)}{dy}\Big|_{y=\x\e^{-1}}\spac{7pt}<-1.
\ee
Finally we apply \refq{390} together with \refq{3104} and the lower bound of \refq{mu.tRn} in Theorem \ref{th:1Rn}, to estimate $\P(N^{\x,\e}\ge s_{x,\e})$ and obtain
\be\label{391}
\lime(f'(g(\frack{|\x|\e^{-1}}{r_\e T}))|\x|\e^{-1})^{-1}\ln\P(N^{\x,\e}\ge s_{x,\e})\spac{7pt}<0.
\ee
By the definition of $S$ we have $\lime f'(g(\frack{|\x|\e^{-1}}{r_\e T}))|\x|\e^{-1}S(\e)=\infty$. Therefore the second limit of~\refq{386} follows directly from \refq{391}. Now, by \refq{386} it follows, that for any family of Borel-measurable sets $(A_\e)_{\e>0}$ we obtain
\be\label{387}
\lime S(\e)\ln\P(N^{\x,\e}\in A_\e)=\lime S(\e)\ln\max_{\substack{m\in A_\e\\m\le s_{\x,\e}}}\P(N^{\x,\e}=m),
\ee
where the probabilities $\P(N^{\x,\e}=m)$ can be calculated by \refq{3104}. Furthermore, the first limit of~\refq{386} allows us to estimate the factor $\P(N_{r_\e T}=m)\frac{\nu^{*m}(dy)}{dy}\Big|_{y=\x\e^{-1}}$ from the right-hand side of \refq{3104} similarly to \refq{434}: For every $\delta>0$ there is $\e$ sufficiently small, such that for all $m=1,2,\dots,s_{\x,\e}$ we obtain
\ba\label{498Rn}
&\big|\ln\P(N_{r_\e T}=m)\frac{\nu^{*m}(dy)}{dy}\bigg|_{y=\x\e^{-1}}+\theta_\e(m)\big|\spac{7pt}\le\delta S(\e)^{-1},
\ea
where $\theta_\e(m):=m(\ln\frack{m}{r_\e T}-1+f(\frack{|\x|\e^{-1}}m)+\frack n2\ln f''(\frack{|\x|\e^{-1}}m)-\frack n2\ln2\pi-\frack{n-1}2\ln(\alpha-1))$.\vspace{2mm}

\noindent By \refq{3104} and \refq{498Rn}  the proof of \refq{Claim3} boils down to the analysis of the function $\theta_\e$ on the intervall $[0,s_{\x,\e}]$.
We continue with estimate of the second derivative of $\theta_\e$. We use an auxiliary function $\sigma\in SR_o$, $\sigma(y)=\ln\frac1y+\frac n2\ln f''(y)$ and obtain
\ba
\theta_\e''(m)&=\frack{d^2}{dm^2}\,m(f(\frack{|\x|\e^{-1}}m)+\sigma(\frack{|\x|\e^{-1}}m))\\
&=\frack{d}{dm}\Big(f(\frack{|\x|\e^{-1}}m)+\sigma(\frack{|\x|\e^{-1}}m)-\frack{|\x|\e^{-1}}m(f'(\frack{|\x|\e^{-1}}m)+\sigma'(\frack{|\x|\e^{-1}}m))\Big)\\
&=\frack{|\x|^2\e^{-2}}{m^3}(f''(\frack{|\x|\e^{-1}}m)+\sigma''(\frack{|\x|\e^{-1}}m)).
\ea
By definition~\refq{389} we have $s_{\x,\e}\ll\e^{-1}$ for $\e\to0$. By Lemma \ref{SR}$(viii)$ combined with Remark~\ref{alpha0} we have $f''\in SR_{\alpha-2}$ and $\sigma''\in SR_{-2}$. Therefore, for every $\delta>0$ and $\e$ sufficiently small we have $|\sigma''(\frack{|\x|\e^{-1}}m)|\le\delta f''(\frack{|\x|\e^{-1}}m)$ for every $m\le s_{\x,\e}$. Furthermore, for every $M>0$, $\delta>0$, $m$ with $|m-\tilde m_{\x,\e}|<(M+1)\sqrt{S(\e)k_\e}^{-1}$ and $\e$ sufficiently small we have 
\be
|\frack{|\x|^2\e^{-2}}{m^3}f''(\frack{|\x|\e^{-1}}m)-\frack{|\x|^2\e^{-2}}{\tilde m_{\x,\e}^3}f''(\frack{|\x|\e^{-1}}{\tilde m_{\x,\e}})|\spac{7pt}\le\delta\frack{|\x|^2\e^{-2}}{\tilde m_{\x,\e}^3} f''(\frack{|\x|\e^{-1}}{\tilde m_{\x,\e}}).
\ee
By $f\in SR_\alpha$ and Lemma \ref{SR}$(ii)$, $(viii)$, combined with the third limit \refq{e:lim3} from Lemma \ref{g_lnRn}$(iii)$ and the definitions of $k_\e$ and $\tilde m_{\x,\e}$ we have that for every $\delta,M>0$ there is some $\e$ sufficiently small, such that  the following estimates are valid for all $|m-\tilde m_{\x,\e}|<(M+1)\sqrt{S(\e)k_\e}^{-1}$:
\ba\label{3.118}
\theta_\e''(m)&\;\le\;(1+\delta)\frack{|\x|^2\e^{-2}}{m^3}f''(\frack{|\x|\e^{-1}}m)\;\le\;(1+2\delta)\tilde m_{\x,\e}^{-1}\Big(\frack{|\x|\e^{-1}}{\tilde m_{\x,\e}}\Big)^2f''(\frack{|\x|\e^{-1}}{\tilde m_{\x,\e}})\\
&\;\le\;(1+3\delta)\tilde m_{\x,\e}^{-1}\alpha(\alpha-1)f(\frack{|\x|\e^{-1}}{\tilde m_{\x,\e}})\;=\;(1+3\delta)\tilde m_{\x,\e}^{-1}\alpha(\alpha-1)f(g_o(\frack{|\x|\e^{-1}}{r_\e}))\\
&\;\le\;(1+4\delta)\tilde m_{\x,\e}^{-1}\alpha|\ln\e|\;=\;(1+4\delta)k_\e.
\ea

The upper bound  corresponding to \refq{3.118} is obtained similarly. Thus, for every $\delta,M>0$ and $\e$ sufficiently small, we obtain
\be\label{397}
(1-\delta)k_\e\spac{7pt}\le\theta_\e''(m)\spac{7pt}\le(1+\delta)k_\e
\ee
for every $|m-\tilde m_{\x,\e}|<(M+1)\sqrt{S(\e)k_\e}^{-1}$. 

\noindent For $M>0$ set \hbox{$\tilde m_{\x,\e,M}:=\lfloor\tilde m_{\x,\e}+M\sqrt{S(\e)k_\e}^{-1}\rfloor+1$.} By \refq{354} and the definition of $g_o$ and $\tilde m_{\x,\e}$ we know, that $m=\tilde m_{\x,\e}$ is a minimizer of $\theta_\e(m)$ on $(0,s_{\x,\e})$ since $\theta_\e$ is monotonically decreasing on $(0,\tilde m_{\x,\e})$ and increasing on $(\tilde m_{\x,\e},s_{\x,\e})$. Furthermore, by Theorem~\ref{th:1Rn} we have that $\lime S(\e)(\theta_\e(\tilde m_{\x,\e})+\ln\mu_{r_\e T}(\x\e^{-1}))=0$. Thus by the definition of $\tilde m_{\x,\e,M}$ together with \refq{387}, \refq{3104}, \refq{498Rn} and \refq{397} we finally obtain
\ba
&\lime S(\e)\ln\P(\sqrt{S(\e)}\tilde N^{\x,\e}>M)\\
&=\lime S(\e)\sup_{\tilde m_{\x,\e,M}\le m\le s_{\x,\e}}\ln\Big(\P(N_{r_\e T}=m)\frac{\nu^{*m}(dy)}{dy}_{\Big|y=\x\e^{-1}}\mu_{r_\e T}(\x\e^{-1})^{-1}\Big)\\
&=-\lime S(\e)(\theta_\e(\tilde m_{\x,\e,M})-\theta_\e(\tilde m_{\x,\e}))\spac{7pt}=-\lime S(\e)\int_{\tilde m_{\x,\e}}^{\tilde m_{\x,\e,M}}\int_{\tilde m_{\x,\e}}^p\theta''(q)\hspace{2pt}dq\hspace{2pt}dp\\
&=-\lime S(\e)\frack12k_\e(\tilde m_{\x,\e,M}-\tilde m_{\x,\e})^2\spac{7pt}=-\frack12M^2.
\ea
The limit $\lime S(\e)\ln\P(\sqrt{S(\e)}\tilde N^{\x,\e}<-M)$ can be calculated equally. 
This completes the proof of \refq{Claim3} and thus of the LDP for $\tilde N^{\x,\e}$.\\

\noindent \textbf{Proof of the LDP for $(\sqrt{S(\e)} \bar N^{\x,\e})_{\e>0}$: } Since the LDP for $(\tilde N^{\x,\e})_{\e>0}$ is already proven, it is sufficient to show that $\lime\sqrt{S(\e)k_\e}\,|m_{\x,\e}-\tilde m_{\x,\e}|=0$.

\noindent Choose $\delta\in(0,\frac16)$ and $K\in(1,\frac43)$. By definition for $y$ sufficiently large, we have $g(K^{-1}y)<g_o(y)<g(Ky)$. 
Furthermore, by Lemma~\ref{g_lnRn}$(i)$ we have $(\ln z)^{\frac1\alpha-\delta}\le g(z)\le(\ln z)^{\frac1\alpha+\delta}$  for $z$ sufficiently large. Combining Lemma~\ref{g_lnRn}$(i)$ with Lemma \ref{SR}$(viii)$ it follows that $\frac d{dz}g(\exp(z))<z^{\frac1\alpha-1+\delta}$ for $z$ sufficiently large. Hence for $y$ sufficiently large we obtain
\ba
|g_o(y)-g(y)|\spac{5pt}\le g(Ky)-g(K^{-1}y)\spac{5pt}=&\int_{\ln y-\ln K}^{\ln y+\ln K}\frack d{dz}g(\exp(z))dz\\
\spac{5pt}\le&2\ln K\sup_{|z-\ln y|\le\ln K}\frack d{dz}g(\exp(z))\spac{5pt}<(\ln y)^{\frac1\alpha-1+\delta}.
\ea
 In particular, we obtain $|g_o(\frac{|\x|\e^{-1}}{r_\e T})-g(\frac{|\x|\e^{-1}}{r_\e T})|\le|\ln\e|^{\frac1\alpha-1+\delta}=o(g(\frac{|\x|\e^{-1}}{r_\e T}))$ for $\e$ sufficiently small. Similarly we obtain	
\ba\label{3100}
\hspace{-5mm}|\tilde m_{\x,\e}-m_{\x,\e}|\,=\,\Big|\frac{|\x|\e^{-1}}{g_o(\frac{|\x|\e^{-1}}{r_\e T})}-\frac{|\x|\e^{-1}}{g(\frac{|\x|\e^{-1}}{r_\e T})}\Big|\,\le\,&|\ln\e|^{\frac1\alpha-1+\delta}\sup_{|y-g(\frac{|\x|\e^{-1}}{r_\e T})|\le|\ln\e|^{\frac1\alpha-1+\delta}}\Big|\frac d{dy}\frac{|\x|\e^{-1}}y\Big|\\
\le2|\ln\e|^{\frac1\alpha-1+\delta}\Big|\frac d{dy}\frac{|\x|\e^{-1}}y\Big|_{y=g(\frac{|\x|\e^{-1}}{r_\e T})}\,\le\,&|\ln\e|^{\frac1\alpha-1+\delta}\frac{2|\x|\e^{-1}}{g(\frac{|\x|\e^{-1}}{r_\e T})^{{}^2}}\,\le\,|\ln\e|^{-1-\frac1\alpha+2\delta}\e^{-1}.
\ea
On the other hand, we have 
\be\label{3101}
\sqrt{S(\e)k_\e}=\sqrt{\alpha|\x|^{-1}\e^2g\big(\frack{|\x|\e^{-1}}{r_\e T}\big)^2|\ln\e|}<|\ln\e|^{\frac12+\frac1\alpha+\delta}\e.
\ee
Combining the estimates \refq{3100} and \refq{3101} with $\delta<\frac16$ yields $\lime\sqrt{S(\e)k_\e}\,|m_{\x,\e}-\tilde m_{\x,\e}|=0$ as required. This completes the proof of the LDP for $\bar N^{\x,\e}$.\\

\bigskip 
\subsubsection{Proof of the asymptotic normality of the increments of $Y^\e$ (Theorem~\ref{Th:3Rn}$(ii)$)}\label{sss:W} 

\bigskip 
\noindent \textbf{Proof of Theorem~\ref{Th:3Rn}$(ii)$:} To lighten notation, we drop the subindex $\eta$ in $\nu_\eta$.

\noindent Due to the rotational invariance, it is sufficient to consider $\x=(\x_1,0,\dots,0)$, $\x_1>0$. Furthermore, we can define functions $F_m:[0,\infty)\mapsto\bR$ for $m=1,2,\dots$ such that $\frac{d\nu^{*m}}{dy}(y)=F_m(|y|)$. We have the disjoint union 
\[
\{\e L_{r_\e T}=\x\}
=\bigcup_{m=1}^\infty\{L_{r_\e T}=\x \e^{-1}, N^{\x, r_\e T}=m\}.
\]
 For any $m\in\bN$ we condition the L\'evy process on the event $\{L_{r_\e T}=\x\e^{-1},N^{\x, r_\e T}=m\}$.  For \hbox{$i=1,2,\dots,m$} let $W^{\x,\e}_{m,i}$ denote the $i$-th jump of this conditioned L\'evy process. 
 By construction $(W^{\x,\e}_{m,1},\dots,W^{\x,\e}_{m,m-1})$ has a joined density that can be given by the conditioned bridge density 
\be\label{4111Rn}
\upsilon_{\x,\e}^m(y_1,\dots,y_{m-1}):=\exp\Big(-\sum_{i=1}^{m-1}f(|y_i|)-f\Big(\Big|\x\e^{-1}-\sum_{i=1}^{m-1}y_i\Big|\Big)\Big)F_m(|\x\e^{-1}|)^{-1},
\ee
where the last jump is uniquely determined by the previous jumps since $W^{\x,\e}_{m,m}=\x\e^{-1}-\sum_{i=1}^{m-1}W^{\x,\e}_{m,i}$. Let $\upsilon_{\x,\e}^{m,1}$ denote the density of the distribution of a single jump $W^{\x,\e}_{m,i}$, which is given as the marginal
\ba\label{3126}
\upsilon_{\x,\e}^{m,1}(y):=\int_{\bR^n}\dots\int_{\bR^n}\upsilon_{\x,\e}^m(y,y_2,\dots,y_{m-1})dy_2\dots dy_{m-1}.
\ea
We see that the law conditioned on the event $\{\e L_{r_\e T}=\x,N^{\x, r_\e T}=m\}$ in \refq{4111Rn} and \refq{3126} only depends on~$\z:=\x\e^{-1}$. Hence for convenience we simplify the notation as follows
\be\label{upsbar}
\bar \upsilon_{\z}^m(y_1,\dots,y_{m-1})=\upsilon_{\x,\e}^m(y_1,\dots,y_{m-1}),\;\mbox{ and }\;\bar \upsilon_{\z}^{m,1}(y):=\upsilon_{\x,\e}^{m,1}(y),
\ee
for which we consider $|\z|\to\infty$. Furthermore it can be read from the structure of the convolution density that $W^{\x,\e}_{m,i}$ is equal in distribution to $W^{\z}_{m,i}:=W^{\z,1}_{m,1}$. Note, by the choice of $\x$ and the definition of $\z$ we have $\z=(\z_1,0,\dots,0)$ with $\z_1>0$.\\

\noindent \textbf{Statements [1.] and [2.]: } We reduce the proof of Theorem~\ref{Th:3Rn}$(ii)$ to the following  complimentary statements [1.] and [2.], where $[1.]$ implies that values $y$ which are far from the mean~$\frac{\z}{m}$ provide an asymptotically negligible contribution to the expectation of $W^{\z,1}_{m,1}$, while statement [2.] gives 
a parabolic upper and lower bounds of $\ln\bar\upsilon^{m,1}_{\z}(y)$ for values $y$ which are close to $\frac{\z}{m}$. \\
\noindent \mbox{$\mathbf{[1.]}$} For $c\in(1-\frac\alpha2,1)$ there are $r,k,m_o>0$ such that  for every $m>m_o$, $\z_1>km$ we have 
\be\label{3.2hoch7}
\int_{\bR^n}\bar \upsilon_{\z}^{m,1}(y)\textbf{1}_{[(\frac{\z_1}{m})^c,\infty)}(|\frack{\z}{m}-y|)dy<\int_{\bR^n}|\frack{\z}{m}-y|\bar \upsilon_{\z}^{m,1}(y)\textbf{1}_{[(\frac{\z_1}{m})^c,\infty)}(|\frack{\z}{m}-y|)dy<\exp(-(\frack{\z_1}{m})^r).
\ee
\mbox{$\mathbf{[2.]}$} For $c\in(1-\frac\alpha2,1)$ let $\bar y_c:=\frac{|\z|}m+2(\frac{|\z|}m)^c$. For every $\delta>0$ and $m, \frac{|\z|}m$ sufficiently large we have for every $|y_1-\frac{\z_1}m|\le(\frac{\z_1}m)^c$ and $|y_2|\le(\frac{\z_1}m)^c$ that
\be\label{3132a}
\Big(\frac{\frac d{dy_1}\bar \upsilon^{m,1}_{\z}(y_1,0,\dots,0)}{\bar \upsilon^{m,1}_{\z}(y_1,0,\dots,0)}-\frac{\frac d{dy_1}\bar\upsilon^{m,1}_{\z}(y_1,0,\dots,0)_{|y_1=\bar y_c}}{\bar\upsilon^{m,1}_{\z}(\bar y_c,0,\dots,0)}\Big)(y_1-\bar y_c)^{-1}
\left\{\begin{array}{l}
\le-(1-\delta)f''(\frack{\z_1}m)\\[2mm]
\ge-(1+\delta)f''(\frack{\z_1}m)
\end{array}\right.
\ee
and
\be\label{3133}
\Big|\frack d{dy_2}\bar\upsilon_{\z}^{m,1}((y_1,y_2,0,\dots,0)+\frack1{\alpha-1}f''(\frack {\z_1}{m})y_2\bar\upsilon_{\z}^{m,1}(y_1,y_2,0,\dots,0)\Big|\le\delta f''(\frack{\z_1}{m})|y_2|\bar\upsilon_{\z}^{m,1}(y_1,y_2,0,\dots,0).
\ee
While for $\frac{\z_1}m$ sufficiently large and $|y_1-\frac{\z_1}m|<(\frac{\z_1}m)^c$ the supremum $\sup_{y_2}\bar\upsilon_{\z}^{m,1}(y_1,y_2,0,\dots,0)$ is obviously obtained for $y_2=0$, the optimizer of $\sup_{y_1}\bar\upsilon_{\z}^{m,1}(y_1,0,\dots,0)$ is unknown. This explains the different structure of estimates \refq{3132a} and \refq{3133}. By the definition we have $\E W^{\z}_{m,i}=\frac{\z}m$. Therefore it seems natural to formulate estimate \refq{3132a} with $\frac{\z_1}m$ instead of $\bar y_c$ on the left-hand side. The seemingly counter-intuitive formulation with $\bar y_c$ is justified by a technical comparison argument which becomes transparent in estimate \refq{3173} below.\\

\noindent\textbf{Claim 3: } The statements $[1.]$ and $[2.]$ imply Theorem~\ref{Th:3Rn}$(ii)$.\\

\noindent\textbf{Proof of Claim 3:} 
\noindent We show that $\lime\P(\mathbb{W}^{\x,\e}\in A)=\P(\mathcal N\in A)$ for every Borel set $A\subset\bR^n$, where $\mathcal N$ is a standard normal vector in $\bR^n$. It is enough to show, for every $\delta>0$ there is $\e$ sufficiently small, such that 
\be\label{3116}
\hspace{-3mm}\upsilon_{\x,\e}(y)
\ge(1-\delta)(\alpha-1)^{\frac{n-1}2}(2\pi f''(\frack{|\x|\e^{-1}}{m_{\x,\e}}))^{-\frac{n}{2}}\exp\Big(-\frack12f''(\frack{|\x|\e^{-1}}{m_{\x,\e}})\Big((y_1-\frack{|x|\e^{-1}}{m_{\x,\e}})^2+\frack1{\alpha-1}\sum_{i=2}^{n}y_i^2\Big)\Big)
\ee
for any $y\in B_{\e,K}:=\{z\in\bR^n\mid|z-\frac{\x\e^{-1}}{m_{x,\e}}|\le Kf''(\frac{|\x|\e^{-1}}{m_{x,\e}})^{{}^{-\frac12}}\}$. The definitions of $B_{\e,K}$, and $\bar W^{\x,\e}$ and $\mathbb{W}^{\x,\e}$, given in \refq{eq:223} and \refq{224} in Theorem \refq{Th:3Rn}$(ii)$, combined with the identity $\frack{|\x|\e^{-1}}{m_{\x,\e}}=g(\frac{|\x|\e^{-1}}{r_\e T})$ yield that \refq{3116} implies implies the following. The density of $\mathbb{W}^{\x,\e}$ lies asymptotically above a standard normal density, uniformly on $\{y\in\bR^n\,\mid\,|(y_1,(\alpha-1)y_2,\dots,(\alpha-1)y_n)|\le K\}$ for any $K>0$, and therefore uniformly on any compact set. 
 We then conclude the desired convergence in distribution as follows: For every $\delta>0$ there is $K$ sufficiently large and $\e$ sufficiently small, such that
\ba\label{3113.a}
\P(\mathbb{W}^{\x,\e}\in A)&\spac{8pt}\ge\int_A(1-\frack\delta2)(2\pi)^{-\frac n2}\exp(-\frack12|y|^2)\textbf{1}_{[0,K]}(|y|)dy\\
&\spac{8pt}\ge(1-\frack\delta2)\P(\mathcal N\in A)-\frack\delta2\spac{8pt}\ge\P(\mathcal N\in A)-\delta
\ea
for every Borell set $A\subseteq\bR^n$. The corresponding upper bound is then obtained by $\P(\mathbb{W}^{\x,\e}\in A)=1-\P(\mathbb{W}^{\x,\e}\in(\bR^n\setminus A))$.\vspace{2mm}

\noindent\textbf{Proof of \refq{3116}:} 
We start with the following estimates of the density $\bar\upsilon_{\z}^{m,1}(\cdot)$.

\noindent Due to the  rotational invariance of $\nu$ and the choice of $\z$ the density $\bar\upsilon_{\z}^{m,1}(\cdot)$ is rotationally symmetric for fixed first component vector in $\bR^n$, that is,
$\bar\upsilon_{\z}^{m,1}(y)=\bar\upsilon_{\z}^{m,1}(y_1,\tilde y_2,0,\dots,0)$ with \hbox{$\tilde y_2=|(0,y_2,\dots,y_n)|$} for any $y=(y_1,\dots,y_n)\in\bR^n$.  Thus by statement $[2.]$ it follows that there is $y_{\z,m}\in \bR$ such that 
\be\label{4115.rn}
\bar \upsilon_{\z}^{m,1}(y)\left\{
\begin{array}{l}
\le\exp\Big(-(1-\delta)((y_1-y_{\z,m})^2+\frack1{\alpha-1}\sum_{i=2}^ny_i^2)\frack12f''(\frac{|\z|}m)\Big)\bar\upsilon_{\z}^{m,1}((y_{\z,m},0,\dots,0))\\
\ge\exp\Big(-(1+\delta)((y_1-y_{\z,m})^2+\frack1{\alpha-1}\sum_{i=2}^ny_i^2)\frack12f''(\frac{|\z|}m)\Big)\bar\upsilon_{\z}^{m,1}((y_{\z,m},0,\dots,0))
\end{array}\right.
\ee
for any $|y-\frac{\z}{m}|\le(\frac{|\z|}{m})^c$. To justify the choice of $\bar\upsilon_{\z}^{m,1}((y_{\z,m},0,\dots,0))$ as standardization factor we need to show that $|y_{\z,m}-\frac{\z_1}m|<(\frac{|\z|}{m})^c$ such that \refq{4115.rn} can be applied for $y=(y_{\z,m},0,\dots,0)$. This follows directly from the estimate \refq{4.114rn} below.  Note that for the derivation of \refq{4.114rn} we use the structure of the exponents in \refq{4115.rn} but not the size of the scaling factor.\\
 We continue to estimate the values of $y_{\z,m}$ and $\bar\upsilon_{\z}^m((y_{\z,m},0,\dots,0))$ in \refq{4115.rn}. Note that by construction we have $\E W^{\z,1}_{m,1}=\frac{\z}{m}$. Let $(W^{\z,1}_{m,1})_1$ denote the first component of $W^{\z,1}_{m,1}$. For any $\tilde\delta>0$ there is a $\delta>0$, such that for $\frac{|\z|}m$ sufficiently large the assumption 
\[
y_{\z,m}>\frack{|\z|}m+\tilde\delta f''\big(\frack{|\z|}{m}\big)^{-\frac12}
\] 
together with \refq{4115.rn} and statement~$[1.]$ leads to the estimate 
\be
\frack{|\z|}m=\E(W^{\z,1}_{m,1})_1\ge y_{\z,m}-\frack{\tilde\delta}3 f''\big(\frack{|\z|}{m}\big)^{-\frack12}-\exp\big(-\big(\frack{|\z|}{m}\big)^r\big)\ge\frack{|\z|}m+\frack{\tilde\delta}3\delta f''\big(\frack{|\z|}{m}\big)^{-\frack12},
\ee
which is absurd. By similar reasoning the assumption $|y_{\z,m}-\frack{|\z|}m|>\delta f''\big(\frack{|\z|}{m}\big)^{-\frac12}$ is shown to be false.  Hence for every $\delta>0$ there is  $k>0$ such that for every and $|\z|>mk$ we obtain 
\be\label{4.114rn}
|y_{\z,m}-\frack{|\z|}m|<\delta f''(\frack{|\z|}m)^{-\frac12}.
\ee
By construction we have $\int_{\bR^n}\bar\upsilon^{m,1}_{\z}(y)dy=1$. Thus by \refq{4115.rn}, \refq{4.114rn} and statement $[1.]$ we obtain the desired estimate for $\bar\upsilon_{\z}^{m,1}(\frac{\z}{m})$:
\be\label{3131}
(1-\delta)(\alpha-1)^{\frac{n-1}2}(2\pi f''(\frack{\z}{m}))^{-\frac n2}\spac{7pt}\le\bar\upsilon_{\z}^{m,1}((y_{\z,m},0,\dots,0))\spac{7pt}\le(1+\delta)(\alpha-1)^{\frac{n-1}2}(2\pi f''(\frack{\z}{m}))^{-\frac{n}{2}}.
\ee
By inserting \refq{4.114rn} and \refq{3131} into \refq{4115.rn} we obtain an upper and lower bound for the density $\bar \upsilon_{\z}^{m,1}(y)$ for any $m$, $\z=(\z_1,0,\dots,0)$ chosen as in statement [2.], and $|y-\frac{\z}{m}|\le(\frac{\z_1}{m})^c$. \vspace{2mm}

\noindent After these preliminary considerations we start to examine the density of $\upsilon_{\x,\e}$. We combine the estimates \refq{4115.rn} - \refq{3131} of $\bar \upsilon_{\z}^{m,1}(y)$ with the estimate of the distribution of the number of jumps given in part~$(i)$ of the theorem. By the definition of $\upsilon_{\x,\e}$ and $\upsilon_{\x,\e}^{m,1}$ the following estimate holds for any $A\subset\bN$:
\be\label{eq:4106}
\upsilon_{\x,\e}(y)\ge\P(N^{\x,\e}\in A) \cdot \inf_{m\in A}\upsilon_{\x,\e}^{m,1}(y)=\P(N^{\x,\e}\in A) \cdot \inf_{m\in A}\bar\upsilon_{\x\e^{-1}}^{m,1}(y),
\ee
 where $N^{\x,\e}$ is defined in Definition \ref{defN}.1. We recall $m_{\x,\e}$ defined in part $(i)$ of the theorem. For small values $\e,\kappa>0$ we set $A=A_{\e,\kappa}:=\{m\in\bN\mid|m-m_{\x,\e}|\le\kappa(S(\e)k_\e)^{-\frac12}\}$. By part $(i)$ of the theorem it follows that $\lime\P(N^{\x,\e}\in A_{\e,\kappa})=1$ for any $\kappa>0$. In order to apply \eqref{4115.rn} to estimate the densities $\inf_{m\in A_\e,\kappa}\upsilon_{\x,\e}^{m,1}(y)$ for $y\in B_{\e,K}$, we continue with the proof that for all $m\in A_{\e,\kappa}$, $y\in B_{\e,K}$ and $\e$ small enough we have $|y-\frack{|x|\e^{-1}}{m}|<(\frack{|x|\e^{-1}}{m_{x,\e}})^c$.\\

\noindent By the definitions of $S$, $k_\e$ and $m_{\x,\e}$ together with Lemma \ref{g_lnRn}$(i)$ we have, $\lime m_{\x,\e}\sqrt{S(\e)k_\e}=\infty$. Thus by the definition of $A_{\e,\kappa}$ we have that for $\e$ sufficiently small the following estimate holds uniformly for all $m\in A_{\e,\kappa}$

\ba\label{3135a}
|\frack{\x\e^{-1}}m-\frack{\x\e^{-1}}{m_{\x,\e}}|&\le|\x|\e^{-1}|m-m_{\x,\e}|\sup_{s\in[m\wedge m_{\x,\e},m\vee m_{\x,\e}]}s^{-2}\\
&\le|\x|\e^{-1}\cdot\kappa\,(S(\e)k_\e)^{-\frac12}\cdot2m_{\x,\e}^{-2}\\
&=2\frack\kappa{\sqrt{\alpha|x|}}g(\frack{|x|\e^{-1}}{r_\e T})|\ln\e|^{-\frac12}.
\ea
By the definition of~$m_{\x,\e}$, using Lemma~\ref{SR}$(ii)$, Lemma~\ref{g_lnRn}$(i)$ and limit~\refq{e:lim2} in Lemma~\ref{g_lnRn}$(iii)$ we obtain
\ba\label{3136b}
\lime g(\frack{|\x|\e^{-1}}{r_\e T})|\ln\e|^{{}^{-\frac12}}f''(\frack{|\x|\e^{-1}}{m_{\x,\e}})^{{}^{\frac12}}&=\sqrt{\lime g\big(\frack{|\x|\e^{-1}}{r_\e T}\big)^{2}|\ln\e|^{-1}f''\big(g(\frack{|\x|\e^{-1}}{r_\e T})\big)}\\
&=\sqrt{(\alpha-1)\lime g\big(\frack{|\x|\e^{-1}}{r_\e T}\big)|\ln\e|^{-1}f'\big(g\big(\frack{|\x|\e^{-1}}{r_\e T}\big)\big)}=\sqrt\alpha.
\ea
By the limit \refq{3136b} inserted into \refq{3135a} we obtain for all $m\in A_{\e,\kappa}$ and $\e$ sufficiently small
\be\label{3126a}
|\frack{|\x|\e^{-1}}m-\frack{|\x|\e^{-1}}{m_{\x,\e}}|\spac{7pt}\le3\frack\kappa{\sqrt{|x|}}f''(\frack{|x|\e^{-1}}{m_{x,\e}})^{{}^{-\frac12}}=o(\frack{|\x|\e^{-1}}{m_{\x,\e}}),
\ee
where in the last step we use $(x\mapsto f''(x)^{-\frac12})\in SR_{1-\frac\alpha2}$. By the definition of $B_{\e,K}$ together with the first estimate in \refq{3126a}  we obtain
\be\label{3127}
|y-\frack{|x|\e^{-1}}{m}|\spac{7pt}\le(K+3\frack\kappa{\sqrt{|x|}})f''(\frack{|x|\e^{-1}}{m_{x,\e}})^{{}^{-\frac12}}\spac{7pt}<(\frack{|x|\e^{-1}}{m_{x,\e}})^c
\ee
for any $y\in B_{\e,K}$ and $m\in A_{\e,\kappa}$, $c>1-\frac\alpha2$, and $\e$ sufficiently small. \\

\noindent Indeed the estimate \refq{4115.rn} together with \refq{4.114rn} and \refq{3131} can be applied to obtain a lower bound for the densities $\bar\upsilon_{\x\e^{-1}}^{m,1}(y)$ in the missing second factor on the right-hand side of \refq{eq:4106}, uniformly for all values $m\in A_{\e,\kappa}$ and $y\in B_{\e,K}$. As a next step we establish estimates for the expressions $f''(\frack{|\x|\e^{-1}}{m})$ and $(y_1-y_{\x\e^{-1},m})^2$ in~\refq{4115.rn} - \refq{3131}. By \refq{3126a} and $f''\in SR_{\alpha-2}$, for any $\delta>0$ there is $\e$ sufficiently small, such that for all $m\in A_{\e,\kappa}$ we have
\be\label{3113}
|f''(\frack{|\x|\e^{-1}}{m})-f''(\frack{|\x|\e^{-1}}{m_{\x,\e}})|<\delta f''(\frack{|\x|\e^{-1}}{m_{\x,\e}}).
\ee
We use the basic estimate $|a^2-b^2|\le2|a||a-b|+|a-b|^2$ to estimate the term $(y_1-y_{\x\e^{-1},m})^2$ in~\refq{4115.rn}. With $a=y_1-\frac{\x_1\e^{-1}}{m_{\x,\e}}$  and  $b=y_1-y_{\x\e^{-1},m}$, we obtain 
\be\label{3134}
\Big|(y_1-y_{\x\e^{-1},m})^2-\Big(y_1-\frac{|\x|\e^{-1}}{m_{\x,\e}}\Big)^2\Big|\le2\Big|y_1-\frac{|\x|\e^{-1}}{m_{\x,\e}}\Big|\Big|y_{\x\e^{-1},m}-\frac{|\x|\e^{-1}}{m_{\x,\e}}\Big|+\Big|y_{\x\e^{-1},m}-\frac{|\x|\e^{-1}}{m_{\x,\e}}\Big|^2,
\ee
where the term $|y_1-\frac{|\x|\e^{-1}}{m_{\x,\e}}|$ can be estimated directly by the first inequality of \refq{3127}. To estimate $|y_{\x\e^{-1},m}-\frac{|\x|\e^{-1}}{m_{\x,\e}}|$ we apply \refq{4.114rn}, \refq{3126a} and \refq{3113}. Consequently, for any $\kappa>0$ there is $\e$ sufficiently small, such that for any $m\in A_{\e,\kappa}$ we have
\be\label{3130}
\Big|y_{\x\e^{-1},m}-\frac{\x_1\e^{-1}}{m_{\x,\e}}\Big|\spac{7pt}\le\Big|y_{\x\e^{-1},m}-\frac{\x_1\e^{-1}}{m}\Big|+\Big|\frac{\x_1\e^{-1}}{m}-\frac{\x_1\e^{-1}}{m_{\x,\e}}\Big|\spac{7pt}\le(3\frack{\kappa}{\sqrt{|\x|}}+\kappa)f''(\frack{|\x|\e^{-1}}{m_{\x,\e}})^{-\frac12}.
\ee
Inserting \refq{3127} and \refq{3130} into \refq{3134} yields, for every $K,\delta>0$ there is $\kappa>0$ sufficiently small, such that 
\ba\label{3131a}
\Big|(y_1-y_{\x\e^{-1},m}&)^2-\Big(y_1-\frac{\x_1\e^{-1}}{m_{\x,\e}}\Big)^2\Big|\\
\spac{7pt}\le&\Big(2\Big(K+\frack{3\kappa}{\sqrt{|\x|}}\Big)\Big(\frack{3\kappa}{\sqrt{|\x|}}+\kappa\Big)+\Big(\frack{3\kappa}{\sqrt{|\x|}}+\kappa\Big)^2\Big)f''(\frack{|\x|\e^{-1}}{m_{\x,\e}})^{-1}
\spac{7pt}\le\delta f''(\frack{|\x|\e^{-1}}{m_{\x,\e}})^{-1}
\ea
for $\e$ sufficiently small. Finally, by inserting \refq{3113} and \refq{3131a} into \refq{4115.rn} we obtain an estimate of the missing second factor of ther right-hand side of \refq{eq:4106}:
\begin{\eq}\label{3.127}
\inf_{m\in A_{\e,\kappa}}\bar\upsilon_{\x\e^{-1}}^{m,1}(y)\ge\exp\Big(-(1+\delta)((y_1-y_{\x\e^{-1},m})^2+\frack1{\alpha-1}\sum_{i=2}^ny_i^2)\frack12f''(\frack{|\x|\e^{-1}}m)\Big)\inf_{m\in A_{\e,\kappa}}\bar\upsilon_{\x\e^{-1}}^{m,1}(\frack{|\x|\e^{-1}}{m})\nonumber\\
\ge\exp\Big(-(1+\delta)\Big(\Big(y_1-\frac{|\x|\e^{-1}}{m_{\x,\e}}\Big)^2+\frack1{\alpha-1}\sum_{i=2}^ny_i^2\Big)\frack12f''(\frack{|\x|\e^{-1}}{m_{\x,\e}})-\delta\Big)\inf_{m\in A_{\e,\kappa}}\bar\upsilon_{\x\e^{-1}}^{m,1}(\frack{|\x|\e^{-1}}{m}).\hspace{5mm}\left.\right.
\end{\eq}
Finally \refq{3131} together with \refq{3113} yields 
\be\label{3.128}
\inf\limits_{m\in A_{\e,\kappa}}\bar\upsilon_{\x\e^{-1}}^{m,1}(\frack{|\x|\e^{-1}}{m})\ge(1-\delta)(\alpha-1)^{\frac{n-1}2}(2\pi f''(\frack{{\x\e^{-1}}}{m_{x,\e}}))^{-\frac n2}
\ee
 for $\e$ small enough. Finally, by the estimates \refq{3.127} and \refq{3.128} with $\tilde\delta<(2+\max\{1,\frac1{\alpha-1}\}K^2)^{-1}\delta$ instead of $\delta$, we obtain the desired bound given in \refq{3116}, uniformly for any $y\in B_{\e,K}$ and $\e$ sufficiently small.

\begin{rem} 
In the light of the upper and lower bound of \refq{4115.rn} it is possible to derive the respective upper bound in \refq{3116}. This implies stronger convergence results than those formulated in Theorem \ref{Th:3Rn}$(ii)$, such as uniform convergence of the density on compact sets. Analogously to~\refq{eq:4106} we may apply the approach 
\be\label{3129}
\upsilon_{\x,\e}(y)\le\sup_{m\in A_{\e,\kappa}}\bar\upsilon_{\x\e^{-1}}^{m,1}(y)+\sum_{m\notin A_{\e,\kappa}}P(N^{\x,\e}=m)\bar\upsilon_{\x\e^{-1}}^{m,1}(y).
\ee
In fact, the first summand on the right-hand side of \refq{3129}, $\sup_{m\in A_{\e,\kappa}}\bar\upsilon_{\x\e^{-1}}^{m,1}(y)$ can be estimated from above with the same arguments as the corresponding infimum from below as given by \refq{3.127} and \refq{3.128}. Moreover, for the second summand, it follows that \hbox{$\lime S(\e)\ln\P(N^{\x,\e}\notin A_{\e,\kappa})<0$} by part $(i)$ of the theorem. However, for $m\notin A_{\e,\kappa}$ it cannot be assumed that  the condition \hbox{$|y-\frac{\x\e^{-1}}{m}|<(\frac{|\x|\e^{-1}}{m})^{{}^c}$} is fulfilled and hence for all $y\in B_{\e,K}$ the upper bound of \refq{4115.rn} can not be used to estimate $\bar\upsilon_{\x\e^{-1}}^{m,1}(y)$. The additional effort to determine an alternative upper bound of $\bar\upsilon_{\x\e^{-1}}^{m,1}(y)$ for the case $m\notin A_{\e,\kappa}$ seems disproportionate to the resulting gain in knowledge and hence omitted at this point.
\end{rem}

\noindent 
In the sequel we continue with the prove of statements $[1.]$ and $[2.]$.\\

\noindent\textbf{Proof of $\mathbf{[1.]}$}: 
We prove a slightly more general result with a polynomial prefactor $h$ under the integral, 
which is used in the proof of $[2.]$. While the first estimate of \refq{3.2hoch7} in statement $\mathbf{[1.]}$ is obviously satisfied, the second estimate follows directly from Lemma \ref{cor:31Rn} with the choice $h_{\z,m}(y)=|y-\frac{\z}m|$ together with the definitions of~$\upsilon$ and $\bar\upsilon$ in \refq{4111Rn} - \refq{upsbar}.

\begin{lemma}\label{cor:31Rn} Let $n$, $\nu$, $f$ and $\alpha$ be defined as in Proposition~\ref{l:2Rn}. Let $\kappa>0$, $c>1-\frac\alpha2$ and let $h_{\z,m}:\bR^n\rightarrow\bR$ satisfy $|h_{\z,m}(y+\frac{\z}{m})|\le|y|^\kappa$ for all $|y|\ge(\frac{|\z|}{m})^c$. Then 
the following estimates hold: For any $\gamma\in(0,2c+\alpha-2)$, $m\in\bN$ and $\frac{|\z|}m$ sufficiently large we have
\ba\nonumber
&\Big|\int_{\bR^n}\dots\int_{\bR^n}h_{\z,m}(y_1)\exp\Big(-\sum_{i=1}^{m-1}f(y_i)-f\Big(\z-\sum_{i=1}^{m-1}y_i\Big)\Big)\textbf{1}_{[(\frac{|\z|}m)^c,\infty)}(|y_1-\frack{\z}{m}|)dy_1\dots dy_{m-1}\Big|\\
\le&\exp(-(\frack{|\z|}m)^\gamma)\int_{\bR^n}\dots\int_{\bR^n}\exp\Big(-\sum_{i=1}^{m-1}f(y_i)-f\Big(\z-\sum_{i=1}^{m-1}y_i\Big)\Big)dy_1\dots dy_{m-1}.
\ea
\end{lemma}
\begin{rem}\label{remark32} $(i)$ By the definition of $W^{\z}_{m,i}$ and its density the statement of Lemma \ref{cor:31Rn} can be interpreted in the sense of $\E(h_{\z,t}(W^{\z}_{m,i})\textbf{1}_{[(\frac {|\z|}m)^c,\infty)}(|(\frack{\z}m-W^{\z}_{m,i}|))<\exp(-(\frack{|\z|}m)^\gamma)$ for $\frac{|\z|}m$ sufficiently large.\\

\noindent$(ii)$ For functions $h_{\z,t}^1$ and $h_{\z,t}^2$ with $h_{\z,t}^1(y)=h_{\z,t}^2(y)$ for $|y-\frac{|\z|}m|<(\frac{|\z|}m)^c$ which satisfy the condition of the lemma, we apply part 1. of the remark for $h_{\z,m}:=|h^1_{\z,m}-h^2_{\z,m}|$. For $\frac{|\z|}m$ sufficiently large, we obtain that the expectations $\E(h^1_{\z,t}(W^{\z}_{m,i}))$ and $\E(h^2_{\z,t}(W^{\z}_{m,i}))$ differ only insignificantly:
$$|\E(h^1_{\z,t}(W^{\z}_{m,i}))-\E(h^2_{\z,t}(W^{\z}_{m,i}))|\spac{7pt}\le\E(h_{\z,t}(W^{\z}_{m,i})\textbf{1}_{[(\frac xm)^c,\infty)}(|(\frack{\z}m-W^{\z}_{m,i}|))\spac{7pt}<\exp(-(\frack{|\z|}m)^\gamma).$$

\end{rem}
\noindent\textbf{Proof of Lemma~\ref{cor:31Rn}: }
Let $\gamma\in(c,\alpha-2+2c)$ and $\rho\in(\gamma-c,\alpha-2+c)$. Recall the definition of the functions $f_{\z,m}$ in equation \refq{deffxm}. We show below, that $\frac{|\z|}m$ can be chosen sufficiently large, such that for the two cases $a=(\frac{|\z|}m)^c$ or $a>(\frac{|\z|}m)^2$ the following estimate is satisfied
\ba\label{4149}
&\int_{\bR^n}\dots\int_{\bR^n}\exp\Big(-\sum_{i=1}^{m-1}f_{\z,m}(y_i)-f_{\z,m}\Big(-\sum_{i=1}^{m-1}y_i\Big)\Big)\textbf{1}_{[a,\infty)}(|y_1|)dy_1\dots dy_{m-1}\\
\spac{8pt}
\le&\exp\Big(-\big(\frack{|\z|}m\big)^{\rho}\cdot a\Big)\int_{\bR^n}\dots\int_{\bR^n}\exp\Big(-\sum_{i=1}^{m-1}f_{\z	,m}(y_i)-f_{\z,m}\Big(-\sum_{i=1}^{m-1}y_i\Big)\Big)dy_1\dots dy_{m-1}.
\ea
First we show, that the statement of the Lemma follows from \refq{4149}. The proof of \refq{4149} is given afterwards. Recall well-known connection between $f$ and $f_{\z,m}$, in particular given by \refq{4.5rn}. Thus, in the case $h_{\z,m}\equiv1$ the assertion of Lemma~\ref{cor:31Rn} follows directly from estimate \refq{4149} by choosing $a=(\frac{|\z|}m)^c$. The same is valid in the case $\kappa\le0$, where we have $\sup_{|y|\ge(\frac{|\z|}m)^c}h_{\z,m}(y+\frac{\z}m)<1$. In the case $\kappa>0$ we apply the integral comparison principle and obtain
\ba\label{3143}
&\int_{\bR^n}\dots\int_{\bR^n}h_{\z,m}(y_1+\frack{\z}m)\exp\Big(-\sum_{i=1}^{m-1}f_{\z,m}(y_i)-f\Big(-\sum_{i=1}^{m-1}y_i\Big)\Big)\textbf{1}_{[(\frac{|\z|}m)^c,\infty)}(|y_1|)dy_1\dots dy_{m-1}\\
\hspace{-8mm}\le\hspace{5pt}&(\frack{|\z|}m)^{2\kappa}\int_{\bR^n}\dots\int_{\bR^n}\exp\Big(-\sum_{i=1}^{m-1}f_{\z,m}(y_i)-f_{\z,m}\Big(-\sum_{i=1}^{m-1}y_i\Big)\Big)\textbf{1}_{[(\frac{|\z|}m)^c,(\frac{|\z|}m)^2)}(|y_1|)dy_1\dots dy_{m-1}\\
&+\int_{\bR^n}|y_1|^\kappa\Big(\int_{\bR^n}\dots\int_{\bR^n}\exp\Big(-\sum_{i=1}^{m-1}f_{\z,m}(y_i)-f_{\z,m}\Big(-\sum_{i=1}^{m-1}y_i\Big)\Big)dy_2\dots dy_{m-1}\Big)\textbf{1}_{[(\frac{|\z|}m)^2,\infty)}(|y_1|)dy_1.
\ea
The first integral on the right-hand side of \refq{3143} can be estimated by \refq{4149} with the choice $a=(\frac{|\z|}m)^c$. In order to estimate the second integral we apply 
\ba\label{3132}
\hspace{-5mm}\int_{\bR^n}|y_1|^\kappa\vartheta(y_1)\textbf{1}_{[k,\infty)}(|y_1|)dy_1\le\sum_{a=k}^\infty\Big((a+1)^\kappa\int_{\bR^n}\vartheta(y_1)\textbf{1}_{[a,\infty)}(|y_1|)dy_1\Big),
\ea
\noindent where in abuse of notation we denote any sum $\sum_{a=k}^\infty q(a)=\sum_{a=0}^\infty q(k+a)$ in the case $k\notin\bN$. We apply \refq{3132} for $k=(\frac{|\z|}m)^2$ and $\vartheta(y_1)=\int_{\bR^n}\dots\int_{\bR^n}\exp(-\sum_{i=1}^{m-1}f_{\z,m}(y_i)-f(-\sum_{i=1}^{m-1}y_i))dy_2\dots dy_{m-1}$. Together with \refq{4149} we obtain
\ba
\int_{\bR^n}\dots\int_{\bR^n}&h_{\z,m}(y_1+\frack{\z}m)\exp\Big(-\sum_{i=1}^{m-1}f_{\z,m}(y_i)-f\Big(-\sum_{i=1}^{m-1}y_i\Big)\Big)\textbf{1}_{[(\frac{|\z|}m)^c,\infty)}(|y_1|)dy_1\dots dy_{m-1}\\
\le&\Big((\frack{|\z|}m)^{2\kappa}\exp(-(\frack{|\z|}m)^{\rho+c})+\sum_{a=\big(\frack{|\z|}m\big)^2}^\infty \big(a+1\big)^\kappa\exp\Big(-(\frack{|\z|}m)^{\rho}a\Big)\Big)\cdot\\
&\hspace{20mm}\cdot\int_{\bR^n}\dots\int_{\bR^n}\exp\Big(-\sum_{i=1}^{m-1}f_{\z	,m}(y_i)-f_{\z,m}\Big(-\sum_{i=1}^{m-1}y_i\Big)\Big)dy_1\dots dy_{m-1}.
\ea
By the choice of $\gamma$ and $\rho$ we have $\rho+c>\gamma$ and therefore
\be
\big(\frack{|\z|}m\big)^{2\kappa}\exp(-(\frack{|\z|}m)^{\rho+c})+\sum_{a=\big(\frack{|\z|}m\big)^2}^\infty (a+1)^\kappa\exp\Big(-(\frack{|\z|}m)^{\rho}\cdot a\Big)\spac{8pt}\le\exp(-(\frack{|\z|}m)^\gamma)
\ee
for $\frack{|\z|}m$ sufficiently large. Indeed, the assertion of the Lemma is a direct consequence of estimate~\refq{4149}.\vspace{2mm}

\noindent \textbf{Proof of \refq{4149}:} We observe the following elementary fact. For any $\ell\in\bN$ and any measurable function $\vartheta: \bR^\ell\to[0,\infty)$, and $r\ge0$ let $A_{r,\vartheta}=\{y\in\bR^\ell\mid\vartheta(y)<r\}$ and $B_{r,\vartheta}=\{y\in\bR^\ell\mid\exp(-\vartheta(y))>r\}$. We have
$$\int_{\bR^\ell}\exp(-\vartheta(y))dy=\int_ o^1\lambda_{\ell}( B_{r,\vartheta})dr=\int_ o^1\lambda_{\ell}( A_{-\ln r,\vartheta})dr=\int_ o^\infty\exp(-z)\lambda_{\ell}( A_{z,\vartheta})dz,$$
where $\lambda_{\ell}$ denotes the Lebesgue measure on $\bR^\ell$.
For the particular choice of $\ell=n(m-1)$ and $\vartheta(y)=\sum_{i=1}^{m-1}f_{\z,m}(y_i)+f_{\z,m}\Big(-\sum_{i=1}^{m-1}y_i\Big)$ and $r\ge0$ we set
\be\label{3145}
A_r=A_{r,\vartheta}=\Big\{y\in(\bR^n)^{m-1}\Big|\sum_{i=1}^{m-1}f_{\z,m}(y_i)+f_{\z,m}\Big(-\sum_{i=1}^{m-1}y_i\Big)<r\Big\}\mb{2mm}{and}\varphi(r)=\lambda_{n(m-1)}(A_r).
\ee
With this notation we have 
\begin{align}
\int_{\bR^n}&\dots\int_{\bR^n}\exp\Big(-\sum_{i=1}^{m-1}f_{\z,m}(y_i)-f_{\z,m}\Big(-\sum_{i=1}^{m-1}y_i\Big)\Big)dy_1\dots dy_{m-1}=\int_o^\infty\exp(-z)\varphi(z)dz.\label{4151}
\end{align}
\begin{rem}\label{rem32}
For later use we note, that the integrals $\int_o^r\exp(-z)\varphi(z)dz$ and $\int_r^\infty\exp(-z)\varphi(z)dz$ are equal to the integral on the left hand side of \refq{4151} over the reduced area of integration $A_r$ and~$(\bR^n)^{m-1}\setminus A_r$ respectively.
\end{rem}
Recall that in the beginning of Subsection \ref{sss:empirical} we assume the function~$f$ to be convex and monotonically non-decreasing. Therefore, by \refq{convexfxm} we have $f_{\z,m}(ky)\ge kf_{\z,m}(y)$ for any $y\neq0$ and $k>1$. 
Conversely, we obtain that the Lebesgue measure of the level sets fulfill the estimate
\be
\varphi(z)\le\left(\frack{z}{z_o}\right)^{n(m-1)}\varphi(z_o)\label{levelsets}
\ee
for $z\ge z_o>0$. We estimate $\exp(-z)\varphi(z)$ for $z\ge2nm$. Therefore we apply \refq{levelsets} with the choices $z_o=nm$ together with the following basic estimate: Let $k:=\frac z{nm}\ge2$. Then we obtain
\ba\nonumber
(\frack z{nm})^{nm}\exp(-(z-nm))&\spac{7pt}=\exp(nm\ln k-(z-nm))\spac{7pt}=\exp(-(1-\frack{nm\ln k}{z-nm})(z-nm))\\
&\spac{7pt}=\exp(-(1-\frack{\ln k}{k-1})(z-nm))\spac{7pt}<\exp(-\frack14(z-nm)),
\ea
thus, combining \refq{levelsets} with the previous estimate we have
\be
\hspace{-10mm}\varphi(z)\exp(-z)\;\leq\;\left(\frack z{nm}\right)^{nm}\varphi(nm)\exp(-z)\spac{7pt}=\,\left(\frack z{nm}\right)^{nm}\exp(-(z-nm))\exp(-nm)\varphi(nm)\nonumber\vspace{-3mm}
\ee
\begin{eqnarray}
&\le&\exp(1-\frack14(z-nm))\exp(-(nm+1))\varphi(nm)\nonumber\\
&\le&\exp(1-\frack14(z-nm))\int_{nm}^{nm+1}\exp(-y)\varphi(y)dy\\
&\leq &\exp(1-\frack14(z-nm))\int_{\bR^n}\dots\int_{\bR^n}\exp\Big(-\sum_{i=1}^{m-1}f_{x,m}(y_i)-f_{x,m}\Big(-\sum_{i=1}^{m-1}y_i\Big)\Big)dy_1\dots dy_{m-1},\nonumber
\end{eqnarray}
\noindent where we use \refq{4151} in the last step. With this estimate at hand, we estimate for any $r\ge2nm$ by \refq{4151} and Remark \ref{rem32}
\ba\label{4155}
&\int_{\bR^n}\dots\int_{\bR^n}\exp\Big(-\sum_{i=1}^{m-1}f_{\z,m}(y_i)-f_{\z,m}\Big(-\sum_{i=1}^{m-1}y_i\Big)\Big)\textbf{1}_{(\bR^{m-1})^n\setminus A_r}(y)dy_1\dots dy_{m-1}\\
=&\int_r^\infty\exp(-z)\varphi(z)dz\\
\le&\int_r^\infty\exp(1-\frack14(z-nm))dz\int_{\bR^n}\dots\int_{\bR^n}\exp\Big(-\sum_{i=1}^{m-1}f_{\z,m}(y_i)-f_{\z,m}\Big(-\sum_{i=1}^{m-1}y_i\Big)\Big)dy_1\dots dy_{m-1}\\
\le&12\exp(-\frack14(r-nm))\int_{\bR}\dots\int_{\bR}\exp\Big(-\sum_{i=1}^{m-1}f_{\z,m}(y_i)-f_{\z,m}\Big(-\sum_{i=1}^{m-1}y_i\Big)\Big)dy_1\dots dy_{m-1}.
\ea
With the choice $r:=\max\{2nm,nm+5(\frac{|\z|}m)^{\rho}\cdot a\}$, for $\frac{|\z|}m$ sufficiently large and $a\ge(\frac{|\z|}m)^c$ we have
\be\label{4156}
12\exp\big(-\frack14(r-nm)\big)\spac{8pt}\le\exp\big(-\big(\frack{|\z|}m\big)^{\rho}\cdot a\big).
\ee
\noindent Inserting \refq{4156} into \refq{4155} yields that the integral on the left-hand side of \refq{4149} is sufficiently small on $(\bR^n)^{m-1}\setminus A_r$.\\

\noindent It remains to estimate the integral on the left-hand side of \refq{4149} on $\{y\in(\bR^n)^{m-1}\mid|y_1|>a\}\cap A_r$.  We consider the cases $r=2nm$ and $r=nm+5(\frac{|\z|}m)^{\rho}a$ separately. \\

\noindent In the case $r=nm+5(\frac{|\z|}m)^{\rho}a$,  the structure of the maximum in the definition of $r$ yields $nm+5(\frac{|\z|}m)^{\rho}a>2nm$, thus $r=nm+5(\frac{|\z|}m)^{\rho}a<10(\frac{|\z|}m)^{\rho}a$. Recall the choice of $\rho<\alpha-2+c$ and let $\delta\in(0,\alpha-2+c-\rho)$. By \refq{convexfxm}, \refq{421} and $f''\in SR_{\alpha-2}$ and by the construction of $A_r$ we obtain for any $a\ge(\frac{|\z|}m)^c$ and $y\in\{z\in(\bR^n)^{m-1}\mid|z_1|>a\}\cap A_r$
$$r\ge f_{x,m}(y_1)\ge|y_1|\big(\frack{|\z|}m\big)^{-c} f_{\z,m}\left(\big(\frack{|\z|}m\big)^c\frack {y_1}{|y_1|}\right)\ge a\big(\frack{|\z|}m\big)^{-c}\cdot\big(\frack{|\z|}m\big)^{\alpha-2-\delta}\big(\big(\frack{|\z|}m\big)^{c}\big)^2>10\big(\frack{|\z|}m\big)^{\rho}a\spac{5pt}>r.$$
\noindent Thus  the set $\{y\in(\bR^n)^{m-1}\mid|y_1|>a\}\cap A_r$ is empty and the integral on the left-hand side of~\refq{4149} over $\{y\in(\bR^n)^{m-1}\mid|y_1|>a\}\cap A_r$ equals 0.\\

\noindent We continue with the second case $r=2nm$. 
\noindent  In this case we infer $nm\ge5\big(\frac{|\z|}m\big)^{\rho}a$. Thus for large values of $\frac{|\z|}m$, the index $m$ is large, too. In order to compare $f_{\z,m}$ with some linear function $\tilde f_{\z,m,a}$ we carry out the following construction. We set 
$$q_a:=\min_{|y|=a}\frac{f_{\z,m}(y)}{|y|}, \quad D_a:=\Big\{y\in\bR^n\,\Big|\,\frac{f_{\z,m}(y)}{|y|}\ge q_a\Big\}
\quad \mbox{ and }\quad k_a:= \min_{y\in D_a}|y|.$$ 
By construction we obtain
\be\label{3152}
\tilde f_{\z,m,a}(y)\spac{7pt}{:=}q_a|y|\;\left\{
\begin{array}{ll}
\ge f_{\z,m}(y)&\mbox{for $y\notin D_a$,}\\
\le f_{\z,m}(y)&\mbox{for $y\in D_a$}.
\end{array}
\right.
\ee
\noindent In order to conclude the proof we need the following lower bound for the quotient $\frac{q_a\cdot k_a}a$. We give this estimate for the cases $a=(\frac{|\z|}m)^c$ and $a\ge(\frac{|\z|}m)^2$ separately. We start with the first case.  For \hbox{$a=(\frac{|\z|}m)^c$} and \hbox{$|y|\le a$}  we apply \refq{421} to approximate $f_{\z,m}(y)$ by the parabolic function $\frac12f''(\frac{|\z|}m)(y_1^2+\frac1{\alpha-1}\sum_{i=2}^{m-2}y_i^2)$. Thus we obtain a lower bound of $k_a$. For any $\delta>0$ and $\frac{|\z|}m$ sufficiently large we have
\be\label{3_155}
k_a\spac{7pt}\in\big(\big(((\alpha-1)\wedge\frack1{\alpha-1})^{\frac12}-\delta\big)a,a\big).
\ee
By the definition of $q_a$, by $f''\in SR_{\alpha-2}$ and estimates \refq{421} and \refq{3_155} we obtain:
For every $0<\delta\le\min\{\alpha-2+2c-\rho,\hspace{3pt}\alpha-1,\hspace{3pt}\frac1{\alpha-1}\}$ we can chose $\frac{|\z|}m$ sufficiently large, such that the following estimates hold:
\ba\label{3154}
\frac{q_a\cdot k_a}a&\ge\big(((\alpha-1)\wedge\frack1{\alpha-1}\big)^{\frac12}-\delta)\min_{|y|=(\frac{|\z|}m)^c}\frac{f_{\z,m}(y)}{|y|}\\
&\ge\big(((\alpha-1)\wedge\frack1{\alpha-1}\big)^{\frac12}-\delta)^2\,\frack12f''\big(\frack{|\z|}m\big)\big(\frack{|\z|}m\big)^c\spac{7pt}\ge\big(\frack{|\z|}m\big)^{\alpha-2+c-\delta}.
\ea
Note, that by the choice of $\delta$ the exponent $\alpha-2+c-\delta$ is positive. 
\noindent For the second case $a\ge(\frac{|\z|}m)^2$ we note that for $|y|\ge\frac12(\frac{|\z|}m)^2$ the right-hand side of the definition~\refq{deffxm} of $f_{\z,m}(y)$ is dominated by its first summand. More precisely, for any $\delta$ choosen as above and $\frac{|\z|}m$ sufficiently large we have
\be\label{3155}
(1-\delta)f(|y|)<f_{\z,m}(y)<(1+\delta)f(|y|).
\ee
 Conversely, by the definition of $k_a$ we obtain
\be\label{3156}
k_a\in((1-\delta)a,a)
\ee
for any $a\ge\big(\frac{|\z|}m\big)^2$ and $\frac{|\z|}m$ sufficiently large. Therefore, by \refq{3156}, \refq{3155}, $f\in SR_\alpha$, the definition of $q_a$, the choices of $a$ and $c<1<\alpha$ we obtain
\be\label{3157}
\frac{q_a\cdot k_a}{a}\spac{3pt}\ge(1-\delta)\min_{|y|=a}\frac{f_{\z,m}(y)}{a}\spac{3pt}\ge(1-\delta)^2f(a)a^{-1}\spac{3pt}\ge a^{\alpha-1-\frac\delta2}\spac{3pt}\ge\big(\frack{|\z|}m\big)^{2\alpha-2-\delta}.
\ee
Note, that by the definion of $c$ we have $2\alpha-2-\delta>\alpha-2+c-\delta$. Thus the right-hand side of~\refq{3154} is valid as an lower bound for $\frac{q_a\cdot k_a}{a}$ in both cases, $a=(\frac{|\z|}m)^c$ and $a\ge(\frac{|\z|}m)^2$.\\

\noindent We continue with the comparison of the integrals appearing in \refq{4149}. To take advantage of the linear function $\tilde f_{\z,m,a}$, we define (similarly to $A_r$) the sets
\[\tilde A_r:=\Big\{y\in(\bR^n)^{m-1}\mid\sum_{i=1}^{m-1}\tilde f_{\z,m,a}(y_i)+\tilde f_{\z,m,a}\Big(\sum_{i=1}^{m-1}y_i\Big)<r\Big\}.\] 
We estimate the quotient of the left- and right-hand side of \refq{4157}, where by \refq{4155} and \refq{4156} in the integral on the left-hand side we restrict the domain of integration $\{y\in(\bR^n)^{m-1}\mid|y_1|>a\}$ to $\{y\in(\bR^n)^{m-1}\mid|y_1|>a\}\cap A_r$. By \refq{4151} and Remark \ref{rem32} we have then\\
\begin{\eq}\label{4157}
&&\frac{\int_{\bR^n}\dots\int_{\bR^n}\exp\Big(-\sum_{i=1}^{m-1}f_{\z,m}(y_i)-f_{\z,m}\Big(-\sum_{i=1}^{m-1}y_i\Big)\Big)\textbf{1}_{A_r}(y)\textbf{1}_{[a,\infty)}(|y_1|)dy_1\dots dy_{m-1}}{\int_{\bR^n}\dots\int_{\bR^n}\exp\Big(-\sum_{i=1}^{m-1}f_{\z,m}(y_i)-f_{\z,m}\Big(-\sum_{i=1}^{m-1}y_i\Big)\Big)dy_1\dots dy_{m-1}}\nonumber\\
&=&\int_o^{r}\exp(-z)\varphi(z)\frac{\lambda_{n(m-1)}(A_z\cap\{y\in(\bR^n)^{m-1}\mid |y_1|\ge a\})}{\lambda_{n(m-1)}(A_z)}dz\Big(\int_o^{\infty}\exp(-z)\varphi(z)dz\Big)^{-1}\nonumber\\
&\le&\int_o^{r}\exp(-z)\varphi(z)\frac{\lambda_{n(m-1)}(A_z\cap\{y\in(\bR^n)^{m-1}\mid y_1\in D_a\})}{\lambda_{n(m-1)}(A_z)}dz\Big(\int_o^{r}\exp(-z)\varphi(z)dz\Big)^{-1}\nonumber\\
&\le&\sup_{z\le r}\frac{\lambda_{n(m-1)}(A_z\cap\{y\in(\bR^n)^{m-1}\mid y_1\in D_a\})}{\lambda_{n(m-1)}(A_z)}\spac{7pt}=\frac{\lambda_{n(m-1)}(A_r\cap\{y\in(\bR^n)^{m-1}\mid y_1\in D_a\})}{\lambda_{n(m-1)}(A_r)}\nonumber\\
&\le&\frac{\lambda_{n(m-1)}(\tilde A_r\cap\{y\in(\bR^n)^{m-1}\mid y_1\in D_a\})}{\lambda_{n(m-1)}(\tilde A_r)},
\end{\eq}
where the last step is justified as follows: For $y\in(\bR^n)^{m-1}$ define the following index function $I_a(y)\subset\{1,2,\dots,m\}$. For $j\in\{1,2,\dots,m-1\}$ we set $j\in I_a(y):\Leftrightarrow y_j\in D_a$. For $j=m$ we set $m\in I_a(y):\Leftrightarrow\sum_{i=1}^{m-1}y_i\in D_a$. Thus by construction we have $|I_a(y)|=\psi_a(y)$ with 
\[
\psi_a(y):=\sum_{i=1}^{m-1}\textbf{1}_{D_a}(y_i)+\textbf{1}_{D_a}\Big(\sum_{i=1}^{m-1}y_i\Big).
\] 
Obviously by the definitions from $A_{r}$ and $\tilde A_{r}$, in particular, the symmetry under permutation of their components, it follows for any choice of $N\subseteq\{1,2,\dots,m\}$ the quotient
\be
\frac{\lambda_{n(m-1)}(\tilde A_r\cap\{y\in(\bR^n)^{m-1}\mid I(y)=N\})}{\lambda_{n(m-1)}(A_r\cap\{y\in(\bR^n)^{m-1}\mid I(y)=N\})}\label{quotient}
\ee
remains constant if $N$ is replaced by any $\tilde N$ with $|\tilde N|=|N|$. Hence it is a function of the cardinality~$|N|$. Moreover, by the definition \refq{3152} of $\tilde f_{\z,m,a}$ we have $\tilde f_{\z,m,a}\le f_{\z,m}$ on $D_a$ and $\tilde f_{\z,m,a}\ge f_{\z,m}$ on $\bR^n\setminus D_a$. This implies that the quotient in \refq{quotient} is non-decreasing as a function of $|N|$ and hence
\be\label{3158}
\frac{\lambda_{n(m-1)}(A_r\cap\{y\in(\bR^n)^{m-1}\mid\psi_a(y)\ge s\})}{\lambda_{n(m-1)}(A_r)}\spac{7pt}\le\frac{\lambda_{n(m-1)}(\tilde A_r\cap\{y\in(\bR^n)^{m-1}\mid\psi_a(y)\ge s\})}{\lambda_{n(m-1)}(\tilde A_r)}
\ee
for every $s=0,1,\dots,m$. Again by symmetry and the definition of $\psi_a$ we obtain 
\be\label{3159}
\frac{\lambda_{n(m-1)}(A_r\cap\{y\in(\bR^n)^{m-1}\mid y_1\in D_a,\psi_a(y)=s\})}{\lambda_{n(m-1)}(A_r\cap\{y\in(\bR^n)^{m-1}\mid \psi_a(y)=s\})}\spac{7pt}=\frac sm.
\ee
By the same argument \refq{3159} remains true with $A_r$ being replaced by $\tilde A_r$. Together with \refq{3159} and \refq{3158} we continue to estimate the right side of \refq{4157}
\begin{\eq}
&&\hspace{-8mm}\frac{\lambda_{n(m-1)}(A_r\cap\{y\in(\bR^n)^{m-1}\mid y_1\in D_a\})}{\lambda_{n(m-1)}(A_r)}\nonumber\\
&=&\sum_{s=1}^{m}\frac sm\frac{\lambda_{n(m-1)}(A_r\cap\{y\in(\bR^n)^{m-1}\mid\psi_a(y)=s\})}{\lambda_{n(m-1)}(A_r)}\nonumber\\
&=&\frac1m\sum_{s=1}^{m}\frac{\lambda_{n(m-1)}(A_r\cap\{y\in(\bR^n)^{m-1}\mid\psi_a(y)\ge s\})}{\lambda_{n(m-1)}(A_r)}\label{3163}\\
&\le&\frac1m\sum_{s=1}^{m}\frac{\lambda_{n(m-1)}(\tilde A_r\cap\{y\in(\bR^n)^{m-1}\mid\psi_a(y)\ge s\})}{\lambda_{n(m-1)}(\tilde A_r)}\nonumber\\
&=&\frac{\lambda_{n(m-1)}(\tilde A_r\cap\{y\in(\bR^n)^{m-1}\mid y_1\in D_a\})}{\lambda_{n(m-1)}(\tilde A_r)}.\nonumber
\end{\eq}
Thus the last step of \refq{4157} is justified. We continue with estimates of the right-hand side of~\refq{4157}.\\[2mm] By construction of $\tilde A_r$ we have, that for any $z>\frac r{q_a}$ the set $\tilde A_r\cap\{y\in\bR^{m-1}\mid|y_1|=z\}$ is empty and due to symmetry, for any $z\le\frac r{q_a}$ the marginal $\lambda_{n(m-2)}(\tilde A_r\cap\{y\in\bR^{m-1}\mid y_1=v\})$ is constant for any $v\in\bR^n$ with $|v|=z$. For the choice $v_z:=(z,0,\dots,0)\in\bR^n$ we obtain
\ba\label{3164}
&\frac{\lambda_{n(m-1)}(\tilde A_r\cap\{y\in(\bR^n)^{m-1}\mid y_1\in D_a\})}{\lambda_{n(m-1)}(\tilde A_r)}\nonumber\\
&\le\frac{\lambda_{n(m-1)}(\tilde A_r\cap\{y\in(\bR^n)^{m-1}\mid|y_1|\ge k_a\})}{\lambda_{n(m-1)}(\tilde A_r)}\\
&=\frac{\int_{k_a}^\infty\lambda_{n(m-1)-1}(\tilde A_r\cap\{y\in(\bR^n)^{m-1}\mid|y_1|=z\})dz}{\int_{o}^\infty\lambda_{n(m-1)-1}(\tilde A_r\cap\{y\in(\bR^n)^{m-1}\mid|y_1|=z\})dz}\noindent\\
&=\frac{\int_{k_a}^{\frac r{q_a}}z^{n-1}\lambda_{n(m-2)}(\tilde A_r\cap\{y\in(\bR^n)^{m-1}\mid y_1=v_z\})dz}{\int_{o}^{\frac r{q_a}}z^{n-1}\lambda_{n(m-2)}(\tilde A_r\cap\{y\in(\bR^n)^{m-1}\mid y_1=v_z\})dz},
\ea
\noindent where we have used that $\lambda_{n-1}(\{y_1\in\bR^n\mid|y_1|=z\})$ is proportional to $z^{n-1}$. 
We continue to estimate the Lebesgue measures for $z\in[0,\frac r{q_a}]$. By the scaling property of the Lebesgue measure we obtain
\ba\label{3166}
&\lambda_{n(m-2)}(\tilde A_r\cap\{y\in(\bR^n)^{m-1}\mid y_1=v_z\})\\
\spac{7pt}=&\lambda_{n(m-2)}\Big(\Big\{y\in(\bR^n)^{m-1}\;|\;\sum_{i=1}^{m-1}q_a|y_i|+q_a\Big|\sum_{i=1}^{m-1}y_i\Big|<r,y_1=v_z\Big\}\Big)\\
\spac{7pt}=&\lambda_{n(m-2)}\Big(\Big\{y\in(\bR^n)^{m-2}\;|\;\sum_{i=1}^{m-2}|y_i|+\Big|v_z+\sum_{i=1}^{m-2}y_i\Big|<\frack r{q_a}-z\Big\}\Big)\\
\spac{7pt}=&\Big(\frac{r-q_az}{r}\Big)^{n(m-2)}\lambda_{n(m-2)}\Big(\Big\{y\in(\bR^n)^{m-2}\;|\;\sum_{i=1}^{m-2}|y_i|+\Big|\frac{r}{r-q_az}v_z+\sum_{i=1}^{m-2}y_i\Big|<\frack r{q_a}\Big\}\Big).
\ea
\noindent Due to symmetry, the map $z\mapsto\lambda_{n(m-2)}(\{y\in(\bR^n)^{m-2}\;|\;\sum_{i=1}^{m-2}|y_i|+|\frac{r}{r-q_az}v_z+\sum_{i=1}^{m-2}y_i|<\frack r{q_a}\})$ is monotonically non-increasing for $z\in[0,\frac r{q_a}]$. We use this monotonicity together with the basic fact that for $\vartheta_1$ and $\vartheta_2$ positive functions on $[0,b]$, $\vartheta_1$ non-increasing and $0<a<b$ we have 
\be\label{3_167}
\frac{\int_a^b\vartheta_1(x)\vartheta_2(x)dx}{\int_0^b\vartheta_1(x)\vartheta_2(x)dx}\le\frac{\int_a^b\vartheta_2(x)dx}{\int_0^b\vartheta_2(x)dx}.
\ee
 Therefore inserting \refq{3166} into \refq{3164} combined with \refq{3_167} we estimate 
\be\label{3_166}
\frac{\int_{k_a}^{\frac r{q_a}}z^{n-1}\lambda_{n(m-2)}(\tilde A_r\cap\{y\in(\bR^n)^{m-1}\mid y_1=v_z\})dz}{\int_{o}^{\frac r{q_a}}z^{n-1}\lambda_{n(m-2)}(\tilde A_r\cap\{y\in(\bR^n)^{m-1}\mid y_1=v_z\})dz}\spac{7pt}\le\frac{\int_{k_a}^{\frac r{q_a}}z^{n-1}(\frac{r-q_az}{r})^{n(m-2)}dz}{\int_{\left.\right._{\hspace{-3pt}o}}^{\frac r{q_a}}z^{n-1}(\frac{r-q_az}{r})^{n(m-2)}dz}.
\ee
\noindent Let $Z$ be a $\beta_{p,q}$-distributed random variable with parameters $p=n$ and $q=n(m-2)+1$, i.e. the distribution of $Z$ exhibits the density function $y\mapsto\frac{\Gamma(p+q)}{\Gamma(p)\Gamma(q)} y^{p-1}(1-y)^{q-1}\textbf{1}_{[0,1]}(y)$ and $\E Z=\frac p{p+q}$. By the choice $p\ge1$ the $\beta_{p,q}$-distribution yields $\P(Z>2b\mid Z>b)\le\P(Z>b)$ for any $b>0$. With this notation the value of the right hand side of \refq{3_166} is equal to $\P(\frac{r}{q_a}Z>k_a)$.  By the estimate \refq{3154} and the choices $a>(\frac{|\z|}m)^c$ and $r=2nm$, the Markov's inequality yields $\P(\frac{r}{q_a}Z>k_a)\le\frac{r}{k_aq_a}\E Z<\frac12$ for~$\frac{|\z|}m$ sufficiently large. Thus $\P(\frac{r}{q_a}Z>k_a)\le2\P(\frac{r}{q_a}Z\in(k_a,2k_a))$. With this estimate at hand we continue to estimate the right-hand side of \refq{3_166}. For~$m$ and $\frac{|\z|}m$ sufficiently large we obtain
\ba\label{3167}
&\frac{\int_{k_a}^{\frac r{q_a}}z^{n-1}(\frac{r-q_az}{r})^{n(m-2)}dz}{\int_{\left.\right._{\hspace{-3pt}o}}^{\frac r{q_a}}z^{n-1}(\frac{r-q_az}{r})^{n(m-2)}dz}\spac{7pt}\le\frac{2\int_{k_a}^{2k_a}z^{n-1}(\frac{r-q_az}{r})^{n(m-2)}dz}{\int_{\frac{k_a}2}^{\frac r{q_a}}z^{n-1}(\frac{r-q_az}{r})^{n(m-2)}dz}\spac{7pt}\le4^{n}\frac{\int_{k_a}^{\frac r{q_a}}\left(\frac{r-q_az}{r}\right)^{n(m-2)}dz}{\int_{\frac{k_a}2}^{\frac r{q_a}}\left(\frac{r-q_az}{r}\right)^{n(m-2)}dz}\\
&\hspace{-0mm}\le\hspace{7pt}4^n\Big(\frac{r-k_aq_a}{r-\frac{k_aq_a}{2}}\Big)^{n(m-2)}\spac{7pt}\le4^n\Big(1-\frac{k_aq_a}{2r}\Big)^{n(m-2)}\spac{7pt}\le4^n\exp\Big(-\frac{k_aq_a}{3r}nm\Big)\\
&\hspace{-0mm}=\hspace{7pt}4^n\exp\Big(-\frac{k_aq_a}{6}\Big)\spac{7pt}\le\exp\left(-\left(\frack{|\z|}m\right)^\rho\cdot a\right).
\ea
Again, in the last step we used estimates \refq{3154}. Combining \refq{4157}, \refq{3163}, \refq{3164}, \refq{3_166} and~\refq{3167} establishes \refq{4149}, and therefore finishes the proofs of Lemma~\ref{cor:31Rn} and of statement~$[1.]$. \\

\noindent\textbf{Proof of $\mathbf{[2.]}$:} Recall the definition of $F_m$ in the first paragraph of Subsection \ref{sss:W}. We start our analysis with an estimate for $F'_m(\cdot)$. Let $h:\bR^n\rightarrow\bR$, $h(y)=\frac d{dy_1}f(|y|)=\frac{y_1}{|y|}f'(|y|)$ for $|y|>0$ and $h(0)=0$. For $x>0$ set $v_x:=(x,0,\dots,0)\in\bR^n$. Due to symmetry  we obtain the following presentation
\ba\label{4.116Rn}
\hspace{-3mm}&F'_m(x)=\frac d{dx}\int_{\bR^n}\dots\int_{\bR^n}\exp\Big(-\sum_{i=1}^{m-1}f(|z_i|)-f\Big(\Big|v_x-\sum_{i=1}^{m-1}z_i\Big|\Big)\Big)dz_1\dots dz_{m-1}\\
&=-\int_{\bR^n}\dots\int_{\bR^n}h\Big(v_x-\sum_{i=1}^{m-1}z_i\Big)\exp\Big(-\sum_{i=1}^{m-1}f(|z_i|)-f\Big(\Big|v_x-\sum_{i=1}^{m-1}z_i\Big|\Big)\Big)dz_1\dots dz_{m-1}\\
&=-F_m(x)\E h\Big(v_x-\sum_{i=1}^{m-1}W^{v_x}_{m,i}\Big)\spac{7pt}=-F_m(x)\E h(W^{v_x}_{m,m})\spac{7pt}=-F_m(x)\E h(W^{v_x}_{m,1}).
\ea
Fix $c\in(1-\frac\alpha2,1)$. With the identity \refq{4.116Rn} we continue to analyze the derivative $\frac d{dy_1}\bar\upsilon_{\z}^{m,1}(y)$ at the position $y=(y_1,0,\dots,0)$ with $|y_1-\frac{|\z|}m|\le(\frac{|\z|}{m})^c$. Recall the choice $\z=(\z_1,0,\dots,0)$,~$\z_1>0$. Therefore, for $\frac{|\z|}m$ sufficiently large, the choice of $y$ yields $|\z-y|=|\z|-|y|$, $\frac{d}{dy_1}|y|=1$ and $\frac{d}{dy_1}|\z-y|=-1$. Consequently we obtain
\be\label{e:4.110}
\frac{\frac{d}{dy_1}\bar\upsilon_{\z}^{m,1}(y)}{\bar\upsilon_{\z}^{m,1}(y)}\spac{7pt}=\frac{\frac{d}{dy_1}(F_1(|y|)F_{m-1}(|\z-y|)F_m(|\z|)^{-1})}{F_1(|y|)F_{m-1}(|\z-y|)F_m(|\z|)^{-1}}\spac{7pt}=\frac{F_1'(y_1)}{F_1(y_1)}-\frac{F_{m-1}'(\z_1-y_1)}{F_{m-1}(\z_1-y_1)}.
\ee
Keep in mind that by the structure of \refq{3132a} we have to estimate the difference  of each of the terms on the right-hand side between the positions $y_1$ and $\bar y_c$. 
Let $\delta\in(0,1)$. By definition of $F_1$ we have $F'_1(y_1)=-f'(y_1)F_1(y_1)$. Due to $f\in SR_\alpha$, the second deviation $f''$ can be considered asymptotically constant in the vicinity of $\frac{|\z|}m$, such that for $\frac{|\z|}m$ sufficiently large we have
\be\label{e:4.111}
-(1+\frack\delta2)f''(\frack{|\z|}m)\spac{8pt}\le\Big(\frac{F'_1(y_1)}{F_1(y_1)}-\frac{F'_1(\bar y_c)}{F_1(\bar y_c)}\Big)(y_1-\bar y_c)^{-1}\spac{8pt}\le-(1-\frack\delta2)f''(\frack{|\z|}m).
\ee
For any $y$ under consideration. It remains to estimate the announced difference of the second term on the right-hand side of equation \refq{e:4.110}
\be\label{3160}
\frac{F'_{m-1}(\z_1-y_1)}{F_{m-1}(\z_1-y_1)}-\frac{F'_{m-1}(\z_1-\bar y_c)}{F_{m-1}(\z_1-\bar y_c)}.
\ee
As a preparation, and in order to estimate the expectation $\E h(W^{v_x}_{m,1})$ on the right-hand side of equation \refq{4.116Rn}, we expand $h$ around the footpoint $h(\frac{\z}m)=f'(\frac{|\z|}m)$. 
 Note that $\frac{d}{dy_2}h(y)=\frac{d}{dy_2}\frac{y_1}{|y|}f'(|y|)=\frac{y_1y_2}{|y|^2}\cdot\big(f''(|y|)-|y|^{-1}f'(|y|)\big)$. Hence with $\tilde y_2=|(0,y_2,\dots,y_n)|$ it follows
\ba\label{3177}
h(y)\spac{7pt}=&f'(\frack{|\z|}m)+(y_1-\frack{\z_1}m)f''(\frack{|\z|}m)+\int_{\frac{|\z|}m}^{y_1}f'''(t)(y_1-t)dt\\
&+\int_o^{\tilde y_2}\frac{y_1t}{\big|(y_1,t)\big|^2}\Big(f''(|(y_1,t)|)-|(y_1,t)|^{-1}f'(|(y_1,t)|)\Big)dt\\
\spac{7pt}=&f'(\frack{|\z|}m)+(y_1-\frack{\z_1}m)f''(\frack{|\z|}m)+o((y_1-\frack{\z_1}m)f''(\frack{|\z|}m))+o(f''(\frack{|\z|}m))\hspace{5mm}\mbox{for $\frac{|\z|}{m}\to\infty$},
\ea
thus 
\ba\label{3164}
\hspace{-46mm}h(\tilde y)-h(y)\spac{7pt}=&(1+o(1))(\tilde y_1-y_1)f''(\frack{|\z|}m)+o(f''(\frack{|\z|}m))
\ea
for any $|\tilde y-\frac{\z}m|\vee|y-\frac{\z}m|\le(\frac{|\z|}m)^c$, where in the last step of \refq{3177} we used $f^{(i)}\in SR_{\alpha-i}$. For~$\z,m$ as before let $h_{\z,m}:\bR^n\to\bR$ such that $h_{\z,m}(y)=h(y)$ for any $|y-\frac{|\z|}m|\le(\frac{|\z|}m)^c$ and \refq{3164} be satisfied with~$h$ being replaced by $h_{\z,m}$ for all $y\in\bR^n$.  Finally we apply \refq{4.116Rn} together with~\refq{3164}, Lemma~\ref{cor:31Rn} and Remark~\ref{remark32}$(ii)$ to bound \refq{3160}:
\ba\label{3173}
&\Big|\frac{F'_{m-1}(\z_1-y_1)}{F_{m-1}(\z_1-y_1)}-\frac{F'_{m-1}(\z_1-\bar y_c)}{F_{m-1}(\z_1-\bar y_c)}\Big|\spac{7pt}=\Big|\E h(W^{u_{\z_1-y_1}}_{m-1,1})-\E h(W^{u_{\z_1-\bar y_c}}_{m-1,1})\Big|\\
\le\hspace{7pt}&\Big|\E\Big[h_{\z,m}(W^{u_{\z_1-y_1}}_{m-1,1})\Big]-\E \Big[h_{\z,m}(W^{u_{\z_1-\bar y_c}}_{m-1,1})\Big]\Big|+\exp\big(-\big(\frack{|\z|}m\big)^\gamma\big)\\
\le\hspace{7pt}&(1+\delta)f''(\frack{|\z|}m)\Big|\E[W^{u_{\z_1-y_1}}_{m-1,1}]-\E[W^{u_{\z_1-\bar y_c}}_{m-1,1}]\Big|+\frack\delta4 f''(\frack{|\z|}m)+2\exp\big(-\big(\frack{|\z|}m\big)^\gamma\big)\\
=\hspace{7pt}&(1+\delta)f''(\frack{|\z|}m)\Big|\frac{\z_1-y_1}{m-1}-\frac{\z_1-\bar y_c}{m-1}\Big|+\frack\delta4 f''(\frack{|\z|}m)+2\exp\big(-\big(\frack{|\z|}m\big)^\gamma\big)\\
\le\hspace{7pt}&\frack\delta3(\bar y_c-y_1)f''(\frack{|\z|}m)
\ea
for $m$ sufficiently large. For the last step recall that by the choice of $\bar y_c$ we have $(\bar y_c-y_1)>(\frac{|\z|}m)^c$ for any $y$ under consideration. Comparing the estimates in \refq{e:4.111} and \refq{3173} we obtain that the contribution of the difference of the second summand in \refq{e:4.110} is negligible in \refq{3132a} as $\frac{|\z|}m\to\infty$. Therefore, inserting \eqref{e:4.111} into \eqref{e:4.110} establishes the first estimate  \refq{3132a} in statement [2.].\vspace{3mm}

 \noindent We continue with the proof of \refq{3133}. To this end we estimate the derivative $\frac d{dy_2}\bar\upsilon_{\z}^{m,1}(y)$ at the position $y=(y_1,y_2,0,\dots,0)$. We choose an approach similar to \refq{e:4.110}. However, from the choice of $y$ and $\z$, here we may not assume $|\z-y|=|\z|-|y|$, and for the corresponding derivatives we have $\frac{d}{dy_2}|y|\not\equiv1$ and $\frac{d}{dy_2}|\z-y|\not\equiv-1$. Then we have
\be\label{eq:4.120}
\frac{\frack d{dy_2}\bar\upsilon_{\z}^{m,1}(y)}{\bar\upsilon_{\z}^{m,1}(y)}\spac{7pt}=(\frack d{d_{y_2}}|y|)\frac{F'_1(|y|)}{F_1(|y|)}+(\frack d{d_{y_2}}|\z-y|)\frac{F'_{m-1}(|\z-y|)}{F_{m-1}(|\z-y|)}.
\ee
Recall $F_1'=-f'\cdot F_1$. In case of $y_2\ge0$ we obtain by the choice of $y_1$, Lemma \ref{SR}$(ii)$ and $f'\in SR_{\alpha-1}$
\be\label{e:4118}
(\frack d{d_{y_2}}|y|)\frac{F'_1(|y|)}{F_1(|y|)}
=-\frack{y_2}{|y|}f'(|y|)\left\{
\begin{array}{l}
\ge-(1+\frac\delta3)y_2\frac m{\z_1}f'(\frac{\z_1}m)\spac{7pt}\ge-(1+\frac\delta2)\frac1{\alpha-1}y_2f''(\frac{|\z|}m)\vspace{2mm}\\
\le-(1-\frac\delta3)y_2\frac m{\z_1}f'(\frac{\z_1}m)\spac{7pt}\le-(1-\frac\delta2)\frac1{\alpha-1}y_2f''(\frac{|\z|}m)
\end{array}\right.
\ee
for $\frac {|\z|}m$ sufficiently large. We remark that for $y_2<0$ all inequalities in \eqref{e:4118} are reversed. In the following we apply \refq{4.116Rn} to estimate $\frac{F'_{m-1}(|\z-y|)}{F_{m-1}(|\z-y|)}$. For any $\delta>0$ there is $m$, $\frac{|\z|}m$ sufficiently large, such that we obtain
\ba\label{3_177}
\Big|\big(\frack{d}{dy_2}|\z-y|\big)\frac{F'_{m-1}(|\z-y|)}{F_{m-1}(|\z-y|)}\Big|
\spac{7pt}=\Big|\frac{y_2}{|\z-y|}\E h(W^{\z-y}_{m-1,1})\Big|\spac{7pt}\le(1+\delta)\z_1^{-1}\E h(W^{\z-y}_{m-1,1})|y_2|
\ea
for any $|y_1-\frac{\z_1}m|\le(\frac{\z_1}m)^c$ and $|y_2|\le(\frac{\z_1}m)^c$. We continue with Lemma~\ref{cor:31Rn} together with Remark \ref{remark32}$(i)$ and obtain
\ba\label{3.166}
\hspace{-7mm}\z^{-1}\E h(W^{\z-y}_{m-1,1})&\spac{7pt}\le\z^{-1}\Big(\E \Big(h(W^{\z-y}_{m-1,1})\textbf{1}_{[0,(\frac{|\z-y|}{m-1})^c]}(|W^{\z-y}_{m-1,1}-\frack{\z-y}{m-1}|)\Big)+\exp\big(-\big(\frack{|\z-y|}{m-1}\big)^\gamma\big)\Big)\\
&\spac{7pt}\le(1+\delta)|\z|^{-1}\sup_{|u-\frac{\z-y}{m-1}|\le(\frac{|\z-y|}{m-1})^{{}^c}}h(u)\spac{9pt}\le(1+2\delta)|\z|^{-1}f'(\frack{|\z|}m)\\
&\spac{7pt}\le(1+3\delta)m^{-1}\frack1{\alpha-1}f''(\frack{|\z|}m).
\ea
Inserting \refq{3.166} into \refq{3_177}, a comparison with \refq{e:4118} yields for $m$, $\frac{|\z|}m$ sufficiently large
\be\label{eq:4.123}
\Big|\big(\frack{d}{dy_2}|\z-y|\big)\frac{F'_{m-1}(|\z-y|)}{F_{m-1}(|\z-y|)}\Big|
\spac{7pt}\le\delta\Big|\big(\frack{d}{dy_2}|y|\big)\frac{F'_{1}(|y|)}{F_{1}(|y|)}\Big|.
\ee
Inserting \refq{e:4118} and \refq{eq:4.123} into \refq{eq:4.120} establishes \refq{3133} in statement [2.]. Thus the proof of~[2.] and therefore of Theorem~\ref{Th:3Rn}$(ii)$ is complete.\vspace{3mm}

\bigskip 

\section*{\textbf{ Acknowledgments}}
\noindent
The research of M.A.H. has been supported by 
the project ``Mean deviation frequencies and the cutoff phenomenon'' (INV-2023-162-2850) 
of the School of Sciences (Facultad de Ciencias) at Universidad de los Andes.

\appendix\markboth{Appendix}{Appendix}
\renewcommand{\thesection}{\Alph{section}}
\numberwithin{equation}{section}

\section{The proofs of the sufficient conditions}\label{a:ProofsSC}\setcounter{equation}{0}

We start by proving Lemma~\ref{l:Orey_Rn} followed by the proof of an auxiliary result in Lemma~\ref{smalljumpRn}. Finally we prove the desired sufficient condition formulated in of Lemma~\ref{hinr_bed2}.\vspace{3mm}

\noindent\textbf{Proof of Lemma~\ref{l:Orey_Rn}:} By the choice of $v_i$ we obtain that \refq{eq:Orey_Rn} is satisfied for any $x\in\bR^n$ with~$|x|=1$ instead of $v_i$. Furthermore, there is $r_o,k>0$, such that for all $r\in(0,r_o)$ and $|x|=1$ we have
\be\label{e:Orey}
\int_{|y|<r}\langle x,y\rangle^2\nu(dy)>kr^\beta.
\ee
Let $\xi^o$ be a L\'evy process
on $\bR^n$ with generating triplet $(0,\nu,\Gamma)$. 
Let $\hat\mu_t^o$ denote the characteristic function of the distribution of $\xi^o_t$. 
By \eqref{e:Orey} there is some $\lambda_o>0$ such that for $\lambda >\lambda_o$ it follows 
\ba\label{440}
-\ln|\hat\mu^o_t(\lambda)|\spac{7pt}=&t\int_{\bR^n}(1-\cos\langle\lambda,y\rangle)\nu(dy)\spac{7pt}\ge t\int_{|y|<|\lambda|^{-1}}\frack13\langle\lambda,y\rangle^2\nu(dy)\\
\spac{7pt}=&t|\lambda|^2\int_{|y|<|\lambda|^{-1}}\frack13\langle\frack{\lambda}{|\lambda|},y\rangle^2\nu(dy)\spac{7pt}\ge\frack13k|\lambda|^{2-\beta}t.
\ea
On the one hand, under the assumption that the distribution of $\xi^o_t$ has a density $\mu^o_{t}$, it can be calculated by Fourier inversion as $\mu^o_{t}(x)=\frac1{(2\pi)^n}\int_{\bR^n}\exp(-i\langle\lambda,x\rangle)\hat\mu^o_{t}(\lambda)d\lambda$. By \refq{440} the following estimate holds 
\ba\label{4_21}
&\hspace{-4mm}\frack1{(2\pi)^n}\sup_{x\in\bR^n}\int_{\bR^n}\exp(-i\langle\lambda,x\rangle)\hat\mu^o_{t}(\lambda)d\lambda\spac{8pt}\le\Big(\int_{[-\lambda_o,\lambda_o]^n}+\int_{\bR^n\setminus[-\lambda_o,\lambda_o]^n}\Big)\left|\hat\mu^o_{1}(\lambda)\right|^td\lambda\\
\spac{7pt}\le&(2\lambda_o)^n+\int_{\bR^n\setminus[-\lambda_o,\lambda_o]^n}\exp(-\frack13k|\lambda|^{2-\beta}t)d\lambda\spac{8pt}\le(2\lambda_o)^n+t^{-\frac n{2-\beta}}\int_{\bR^n}\exp(-\frack23k|\lambda|^{2-\beta})d\lambda\\
\spac{7pt}\le&K(1\vee t^{-\frac n{2-\beta}})
\ea
for $K$ choosen sufficiently large. On the other hand, the finiteness in \refq{4_21} is sufficient to justify the assumption of existence of $\mu_t^o$ and its upper bound 
\be\label{4_21a}
\sup_{x\in\bR^n}\mu^o_{t}(x)\spac{7pt}\le K(1\vee t^{-\frac n{2-\beta}})
\ee
\noindent By construction the distribution of $\xi_t$ equals the convolution of the distribution of~$\xi^o_t$ and the corresponding marginal distribution of the Brownian component of $\xi$. Thus the distribution of $\xi_t$ has a density $\mu_t$ and for any $t>0$ we obtain $\sup_{x\in\bR^n}\mu_{t}(x)\le\sup_{x\in\bR^n}\mu^o_{t}(x)$. Thus \refq{4_21a} implies the statement of Lemma~\ref{l:Orey_Rn} and finishes the proof.\\

\noindent In order to prove Lemma~\ref{hinr_bed2} we need slightly modified versions of two results that have been proven in \cite{Imkeller}. The first result is a multidimensional version of Lemma~3.3 in \cite{Imkeller}.
\begin{lemma}\label{smalljumpRn}
Let $n\in\bN$, $\nu$ be a L\'evy measure on $\bR^n$, $\Lambda>1$, $h_\Lambda\ge1$ and $T_\Lambda>0$, such that $\lim_{\Lambda\rightarrow\infty}\frac{ \Lambda}{h_\Lambda}=\lim_{\Lambda\rightarrow\infty}\frac{\Lambda}{h_\Lambda T_\Lambda}=\infty$. 
Let $\sigma^2\in \bR^{n\times n}$ be symmetric and nonnegative definite, $\Gamma\in\bR^n$ and for given $\Lambda>0$ consider $\xi^\Lambda$ be a L\'evy process defined by the characteristic triplet $(\sigma^2,\nu|_{\{y\in\bR^n\mid~|y|\le h_\Lambda\}},\Gamma)$. Then for every $\delta>0$, $k\ge0$ the parameter $\Lambda$ can be chosen sufficiently large, such that
\be
\sup_{t\le T_\Lambda}(1\vee t^{-k})\P(|\xi^\Lambda_t|>\Lambda)~\le\exp\Big(-(1-\delta)\frac{\Lambda}{h_\Lambda}\ln\frac{\Lambda}{h_\Lambda T_\Lambda}\Big).\label{eq:smallRn}
\ee
\end{lemma}
\noindent \textbf{Proof of Lemma~\ref{smalljumpRn}:} Let $\delta\in(0,1)$,  ${N_\delta=n[1+3\delta^{-1}]}$, $\tilde D_\delta=\{y=(y_1,\dots,y_n)\mid|y_i|\in\{0,\frac 1{N_\delta},\frac 2{N_\delta},\dots,1\}\}$ and \hbox{$D_\delta=\{\frac y{|y|}\mid y\in\tilde D_\delta\}$.} By definition we have $|y|=1$ for every $y\in D_\delta$ and $\sup_{y\in D_\delta}\langle y,z\rangle\ge(1-{\frac\delta3})|z|$ for every $z\in\bR^n$. By the Chebyshev inequality we obtain 

\be\label{414Rn}
\P(|\xi^\Lambda_t|>\Lambda)\spac{8pt}\le\sum_{y\in D_\delta}\P(\langle y,\xi_t^\Lambda\rangle>(1-{\frack\delta3})\Lambda)\spac{8pt}\le\sum_{y\in D_\delta}\exp(-(1-\frack\delta3)s\Lambda)\E\exp(-\langle sy,\xi_t^\Lambda\rangle).
\ee
for every $s>0$. Let $m:=\int_{\bR^n}(1\wedge|y|^2)\nu(dy)$ and $s=s_{\Lambda,t}={\frac1{h_{\Lambda}}\ln\frac{\Lambda}{h_{\Lambda}t}}$. We may estimate $\E\exp(-\langle ky,\xi_t\rangle)$ by the characteristic function of $\xi_t$. We can choose $\Lambda$ sufficiently large, such that the following estimate holds for all $y\in D_\delta$ and $t\le T_\Lambda$:
\ba\label{415Rn}
&\ln\E\exp(-\langle s_{\Lambda,t}y,\xi_t^\Lambda\rangle)\\
&\spac{8pt}=t\Big(\frack12\langle s_{\Lambda,t}y,\sigma^2 s_{\Lambda,t}y\rangle+\langle \Gamma,s_{\Lambda,t}y\rangle+\int_{\bR^n}(\exp(\langle z,s_{\Lambda,t}y\rangle)-1-\langle z,s_{\Lambda,t}y\rangle\textbf{1}_{[0,1]}(|z|))\textbf{1}_{[0,h_\Lambda]}(|z|)\nu(dz)\Big)\\
&\spac{8pt}\le t\Big(s_{\Lambda,t}^3+\int_{\bR^n}\Big(\langle z,s_{\Lambda,t}y\rangle\textbf{1}_{(1,h_\Lambda]}(|z|)+\sum_{r=2}^\infty\frack1{r!}\langle z,s_{\Lambda,t}y\rangle^r\textbf{1}_{[0,h_\Lambda]}(|z|)\Big)\nu(dz)\Big)\\
&\spac{8pt}\le t\Big(s_{\Lambda,t}^3+\int_{\bR^n}\Big({(h_\Lambda s_{\Lambda,t})}(1\wedge|z|^2)+\sum_{r=2}^\infty\frack1{r!}{(h_\Lambda s_{\Lambda,t})}^r(1\wedge|z|^2)\Big)\nu(dz)\Big)\\
&{\spac{8pt}= t\Big(s_{\Lambda,t}^3+\int_{\bR^n}(1\wedge|z|^2)\nu(dz)\sum_{r=1}^\infty\frack1{r!}{(h_\Lambda s_{\Lambda,t})^r\Big)}}\\
&\spac{8pt}\le t(s_{\Lambda,t}^3+m\exp(s_{\Lambda,t}h_\Lambda))\spac{8pt}=ts_{\Lambda,t}^3+{\frac{\Lambda}{h_\Lambda}}.
\ea
{For the next step we remember the choices $h_\Lambda>1$ and $\lim_{\Lambda\to\infty}\frac{\Lambda}{h_\Lambda}=\lim_{\Lambda\to\infty}\frac{\Lambda}{h_\Lambda T_\Lambda}=\infty$. For every $k,\delta>0$ we can chose $\Lambda$ sufficiently large, such that from \refq{415Rn} together with \refq{414Rn} and the definition of $s_{\Lambda,t}$ we obtain
\ba
\P(|\xi_t^\Lambda|>\Lambda)\;&\le\;|D_\delta|\exp(-((1-\frack\delta3)s_{\Lambda,t}\Lambda-ts_{\Lambda,t}^3-m\frack{\Lambda}{h_\Lambda}))\\
&=\;|D_\delta|\exp(-(1-\frack\delta3)(1-\frack t\Lambda s_{\Lambda,t}^2-m(\ln\frack{\Lambda}{h_\Lambda t})^{-1})s_{\Lambda,t}\Lambda)
\spac{7pt}\le\exp(-(1-\frack\delta2)s_{\Lambda,t}\Lambda)\\
&=\spac{7pt}\exp(-(1-\frack\delta2)\frack{\Lambda}{h_\Lambda}\ln\frack{\Lambda}{h_\Lambda t})\spac{7pt}\le(1\wedge T_\Lambda^k)\exp(-(1-\delta)\frack{\Lambda}{h_\Lambda}\ln\frack{\Lambda}{h_\Lambda t})\\
&\le\spac{7pt}(1\wedge t^k)\exp(-(1-\delta)\frack{\Lambda}{h_\Lambda}\ln\frack{\Lambda}{h_\Lambda T_\Lambda})
\ea
uniformly for every $0< t\leq T_\Lambda$. This completes the proof of Lemma~\ref{smalljumpRn}.}\\

\noindent The second result is a slightly modified version of Lemma~5.1 in~\cite{Imkeller}. We stress that this result covers scalar L\'evy processes. 

\begin{lemma}
\label{levytail}
Let $f\in R_\alpha$ be a non decreasing function with $\alpha>1$ and $L = (L_t)_{t\geq 0}$ a L\'evy process with values in $\mathbb{R}$ with jump measure $\nu$. 
Let the jump measure  $\nu$ of $L$ satisfy $\nu([\Lambda,\infty))\le\exp(-f(\Lambda))$ for each $\Lambda\in(0,\infty)$. For $\gamma<1$ denote $d_{\alpha,\gamma}:=\alpha(\frac{\alpha-1}{1-\gamma})^{-(1-\frac1\alpha)}$. Let $q_\Lambda:=\sup\{y\in(0,\infty)\mid f(y)<\ln \Lambda\}$ and $D_\Lambda:=\exp(-d_{\alpha,\gamma}\frac{\ln \Lambda}{q_\Lambda}\Lambda)$. Then for every $\delta\in(0,1)$, $\gamma<1$ the parameter $\Lambda$ can be chosen sufficiently large, such that
\be
\sup_{t<\Lambda^\gamma}\P(L_t>\Lambda)\le D_\Lambda^{1-\delta}.\label{eq:25}
\ee
\end{lemma}

\noindent The first difference between the~present and the original version is, that we consider the supremum $\sup_{t<\Lambda^\gamma}$ for some $\gamma<1$ instead of the supremum $\sup_{t<T}$ for some fixed value $T>0$. Indeed, for the choice $\gamma=0$ the parameter $d_{\alpha,o}$ coincides with the parameter $d_\alpha$ of the original lemma. 
The second slight generalization is the inequality $\nu([\Lambda,\infty))\le\exp(-f(\Lambda))$ instead of an equality. 
A proof of Lemma~\ref{levytail} can be found in the Ph.D. thesis \cite[Lemma 2.25]{Wetzel2}. The necessary modifications compared to the original proof in~\cite{Imkeller} are straightforward.\\

\noindent In the proof of Lemma \ref{hinr_bed2} we will apply estimate \refq{eq:25} together with the following estimate of~$D_\Lambda$. For $f\in R_\alpha$, $\alpha>1$ and any $\delta>0$ there is $y$ sufficiently large, such that $y^{\alpha-\delta}<f(y)<y^{\alpha+\delta}$. Conversely by the definition of $q_\Lambda$ we obtain
\be\label{qlambda}
(\ln\Lambda)^{\frac1\alpha-\delta}\le q_\Lambda\le(\ln\Lambda)^{\frac1\alpha+\delta}\mb{8mm}{and thus}(\ln\Lambda)^{1-\frac1\alpha-\delta}\Lambda\le \ln D_\Lambda\le(\ln\Lambda)^{1-\frac1\alpha+\delta}\Lambda
\ee
any $\delta>0$ and $\Lambda$ sufficiently large.\\

\noindent \textbf{Proof of Lemma~\ref{hinr_bed2}$(i)$}: Let $\aleph_o=\aleph_1\wedge \aleph_2$. Let $\mu_t$, $\mu^1_t$ and $\mu^2_t$ denote the density of the distributions of $\xi_t$, $\xi^1_t$ and $\xi^2_t$ respectively. We obtain $\mu_t$ as a convolution of $\mu^1_t$ and $\mu^2_t$. 
 For any $\Lambda>0$ and $|y|=\Lambda$ we have
\ba\label{423}
\mu_t(y)\;\le&\;\P(|\xi^1_t|\le\frack12\Lambda)\sup_{|z|<\frac12\Lambda}\mu_t^2(y-z)+\P(|y-\xi^2_t|\ge\frack12\Lambda)\sup_{|z|>\frac12\Lambda}\mu_t^1(z)\\
\;\le&\;\sup_{|z|\ge\frac12\Lambda}\mu^2_t(z)+\sup_{|z|\ge\frac12\Lambda}\mu^1_t(z).
\ea
By hypothesis both, $\mu^1_t(x)$ and $\mu^2_t(x)$, satisfy \refq{BEDIIRn} with parameter $\aleph_o$. Therefore, by~\refq{423} it follows, that $\mu_t(x)$ satisfies \refq{BEDIIRn} {for any parameter $\aleph\in(\frac 1\alpha,\aleph_o)$.}\vspace{3mm}

\noindent\textbf{Proof of Lemma~\ref{hinr_bed2}$(ii)$: III.a implies Hypothesis~II: } 
Let the generating triplet $(\sigma^2,\nu_\xi,\Gamma)$ of a L\'evy Process $\xi$ 
satisfy condition~III.a. Consider a L\'evy process $(\xi^o_t)_{t\ge0}$ with generating triplet $(0,\nu_{\xi},\Gamma)$ and let $(W_t)_{t\ge0}$ denote the Gaussian component of $\xi$ with the respective marginal densities $\mu_{\xi^o,t}$ and $\mu_{W,t}$.  By construction the distribution of $\xi_t$ equals the convolution of the distributions of $\xi^o_t$ and $W_t$. For $\Lambda>0$ and $|y|=\Lambda$ we have
\ba\label{424}
\mu_{\xi,t}(y)\spac{7pt}\le&\P(|\xi^o_t-y|\le\frack12\Lambda)\sup_{|z|\le\frac12\Lambda}\mu_{W,t}(z)+\sup_{|z|\ge\frac12\Lambda}\mu_{W,t}(z)\\
\spac{7pt}\le&\P(|\xi^o_t|\ge\frack12\Lambda)\sup_{z\in\bR^n}\mu_{W,t}(z)+\sup_{|z|\ge\frac12\Lambda}\mu_{W,t}(z).
\ea
By the hypothesis $\det \sigma^2>0$, the choice $\gamma<1$ and the well-known shape of $\mu_{W,t}$, we obtain that uniformly for all $t\le\Lambda^\gamma$ the second summand on the right side of \refq{424} is sufficiently small  to satisfy~\eqref{BEDIIRn} for any $\aleph>1-\frac1\alpha$.

\noindent We continue with the first summand. By hypothesis III.a there is some $h>0$ for which $\nu_\xi(\{z\in\bR^n\mid|z|>h\})=0$.  Note that $\xi^o$ is a special case of $\xi^\Lambda$ in Lemma \ref{smalljumpRn}, which is independent of $\Lambda$ for this constant, i.e. $\Lambda$-independent choice $h_\Lambda=h$. Let $T_\Lambda=\Lambda^\gamma$ and $k>\frac n2$. Then Lemma~\ref{smalljumpRn} yields
\be\label{425}
\limsup_{\Lambda\rightarrow\infty}\sup_{t\le\Lambda^\gamma}\frac{\ln(t^{-k}\P(|\xi^o_t|\ge\frac12\Lambda))}{\Lambda\ln\Lambda}\spac{8pt}<0.	
\ee
By the well-known density of the marginal distributions of a Brownian motion together with $\det \sigma^2>0$ and the choice of $k$ we obtain
\be\label{426}
\lim_{t\to0}\,t^k\sup_{y\in\bR^n}\mu_{W,t}(y)=0.
\ee
Thus by \refq{425} and \refq{426} the first summand on the right side of \refq{424} is sufficiently small  to satisfy~\eqref{BEDIIRn} for any $\aleph\in(1-\frac1\alpha,1)$. \\

\noindent\textbf{III.b implies Hypothesis~II: } For every $t>0$ the distribution of $\xi_t$ is equal to the convolution of the distribution of $\xi_{\frac t2}$ with itself.
 Similarly  to \refq{423}, for $\Lambda>0$ and $|y|=\Lambda$ we obtain
\ba\label{4_22}
\mu_{\xi,t}(y)&\le\P(|y-\xi_{\frac t2}|\le\frack12\Lambda)\sup_{|z|\le\frac12\Lambda}\mu_{\xi,\frac t2}(z)+\P(|\xi_{\frac t2}|\ge\frack12\Lambda)\sup_{|z|\le\frac12\Lambda}\mu_{\xi,\frac t2}(y-z)\\
&\le2\P(|\xi_{\frac t2}|>\frack12\Lambda)\sup_{z\in\bR^n}\mu_{\xi,\frac t2}(z).
\ea
Let~$k>\frac{n}{2-\beta}$. Then we estimate $t^{-k}\P(|\xi_{x,\frac t2}|>\frack12\Lambda)$ similarly to \refq{425}. Finally by Lemma~\ref{l:Orey_Rn} we obtain
\be
\lim_{t\to0} t^k\sup_{z\in\bR^n}\mu_{\xi,\frac t2}(z)=0.
\ee
Indeed $\mu_{\xi,t}(y)$ is sufficiently small and satisfies \eqref{BEDIIRn}  for any $\aleph\in(1-\frac1\alpha,1)$.\\

\noindent\textbf{IV implies Hypothesis~II:} Let the generating triplet $(\sigma^2,\nu_\xi,\Gamma)$ of a L\'evy process $\xi$ satisfy condition~IV. Let $\xi^o$ and $\xi^1$ be independent L\'evy processes on $\bR^n$ with generating triplets $(\frac{1}{2}\sigma^2,\nu_{\xi^o},\tilde\Gamma)$ and $(0,\nu_{\xi^1},0)$ respectively, where $\nu_{\xi^o}(\cdot):=\nu_\xi(\cdot\cap[-\Delta,\Delta]^n)$  and $\nu_{\xi^1}(\cdot):=\nu_\xi(\cdot\setminus[-\Delta,\Delta]^n)$ and $\tilde\Gamma$ is defined according to \refq{tildegamma}. Let $W$ denote a Brownian motion with covariance matrix $\frac12\sigma^2$, which is independent of $\xi^o$ and $\xi^1$. 
By construction the L\'evy process $\xi$ is identically distributed to the sum $\xi^o+\xi^1+W$ and  the process $\xi^o$ satisfies condition III.a, and therefore Hypothesis~II. By part~$(i)$ it is therefore sufficient to prove that the process $(\xi^1+W)$ satisfies Hypothesis~II. \\
For $t>0$ denote by $\tilde \mu^1_t$  the density of the distribution of~$\xi^1_t+W_t$ and by $\mu^W_t$~the density of~$W_t$. Let $r_t:=\sup(\{0\}\cup\{|y|\mid\mu^W_t(y)>1\})$. By the convolution of the distributions of $\xi^1_t$ and $W_t$ we have
\ba\label{429}
\tilde \mu^1_t(y)\spac{7pt}\le&\P(|y-\xi^1_t|>\frack12\Lambda)\sup_{|z|>\frack12\Lambda}\mu^W_t(z)+\P(|y-\xi^1_t|\in[r_t,\frack12\Lambda])\sup_{|z|\in[r_t,\frac12\Lambda]}\mu^W_t(z)+\\
&+\P(|y-\xi^1_t|\le r_t)\sup_{|z|\le r_t}\mu^W_t(z).
\ea
Similarly to the proof of Condition III.a we have, that the factor $\sup_{|z|>\frac{1}2\Lambda}\mu^W_t(z)$ and therefore the first summand on the right-hand side of \refq{429} is sufficiently small to satisfy~\eqref{BEDIIRn}  for any~$\aleph>1-\frac1\alpha$.

\noindent Let $\zeta$ be a L\'evy process on $\bR$ with generating triplet $(0,\nu_{\zeta},0)$, where 
$\nu_\zeta(\cdot)=\nu_{\xi^1}(\{z\in\bR^n\mid |z|\in\,\cdot\,\})$.  By construction we have
\be\label{A.1.182}
\P(|y-\xi^1_t|\in[r_t,\frack12\Lambda])\;\le\;\P(|\xi^1_t|\ge\frack12\Lambda)\;\le\;\P\Big(\sum_{s\le t}|\xi^1_{s}-\xi^1_{s-}|\ge\frack12\Lambda\Big)\;=\;\P(\zeta_t\ge\frack12\Lambda).
\ee
By construction the scalar L\'evy process $\zeta$ satisfies $\nu_{\zeta}([\Lambda,\infty))\le\exp(-f_\xi(\Lambda))$. Hence the tail-distribution of $\zeta_t$ can be estimated by Lemma~\ref{levytail}. For any $\aleph\in(1-\frac1\alpha,1-\frac1{\alpha_\xi})$ we apply \refq{eq:25} together with \refq{qlambda} for $\delta\in(0,1-\frac1{\alpha_\xi}-\aleph)$. We obtain
\be\label{427}
\lim_{\Lambda\rightarrow\infty}\sup_{t\le \Lambda^\gamma}\frac{\ln\P(\zeta_t>\Lambda))}{\Lambda(\ln \Lambda)^\aleph}\spac{7pt}=-\infty.
\ee 
By the definition of $r_t$ we have $\sup_{|y|\in[r_t,\frac12\Lambda]}\mu^W_t(y)\le1$. Hence by  \refq{427} it follows that the second summand on the right-hand side of \refq{429} is sufficiently small to satisfy~\eqref{BEDIIRn} for any $\aleph\in(1-\frac1\alpha,1-\frac1{\alpha_\xi})$.

\noindent To estimate the missing third summand in \refq{429}, we use the following estimates. By the well-known marginal distributions of a Brownian motion together with the definition of $r_t$ we obtained
\be
\lim_{t\rightarrow0}\;t^\frac n2\sup_{z\in\bR^n}\mu_t^W(x)<\infty\mb{15mm}{and}\lim_{t\rightarrow0}\frac{r_t}{\sqrt{t|\ln t|}}<\infty.\label{430}
\ee 
Furthermore, let $\bar\nu:=\frac{\nu_{\xi,1}}{\nu_{\xi,1}(\bR^n)}$ and $\bar\nu^{*k}$ denote the $k-th$ convolution of $\bar\nu$ with itself. Note that by~\refq{4_16} for every $k\in\bN$, $k>0$ we have
\be\label{433}
\sup_{z\in\bR^n\setminus\{0\}}\lim_{r\rightarrow0}\;r^{-(n-1)}\bar\nu^{*k}_{\xi,1}(\{y\in\bR^n\mid|y-z|<r\})\spac{7pt}<\infty.
\ee
Let $q\in(1-\frac1\alpha,1-\frac1{\alpha_\xi})$. We consider the supremum in \eqref{BEDIIRn} over $t\in(0,\exp(-\Lambda(\ln\Lambda)^{q}))$ and $t\in(\exp(-\Lambda(\ln\Lambda)^{q}),\Lambda^\gamma)$ separately.\\
We start with the case 
%
 $t>\exp(-\Lambda(\ln\Lambda)^{q})$. Combining~\refq{427} for some $\aleph\in(q,1-\frac1{\alpha_\xi})$  with the first limit of \refq{430} we obtain that the third summand on the right-hand side of \refq{429} is sufficiently small  to satisfy~\eqref{BEDIIRn} for any $\aleph\in(q,1-\frac1{\alpha_\xi})$ and by a monotonicity argument in \refq{429} thus for any $\aleph\in(1-\frac1\alpha,1-\frac1{\alpha_\xi})$.\\
We continue with the case $t<\exp(-\Lambda(\ln\Lambda)^{q})<1$. Then for any $A\subset\bR^n\setminus\{0\}$ the probability  $\P(\xi^1_t\in A)$ can be estimated as
\ba\label{432}
\P(\xi^1_t\in A)&=\sum_{k=1}^\infty\frac{(\nu_{\xi,1}(\bR^n)t)^k}{k!}\bar\nu^{*k}(A)\exp(-(\nu_{\xi,1}(\bR^n)t))\\
&\le t\sum_{k=1}^{[\frac n2]}\frac{(\nu_{\xi,1}(\bR^n))^k}{k!}\bar\nu^{*k}(A)+t^{\frac n2}\sum_{k=[\frac n2]+1}^{\infty}\frac{(\nu_{\xi,1}(\bR^n))^k}{k!}\bar\nu^{*k}(A)\\
&\le t\sum_{k=1}^{[\frac n2]}\frac{(\nu_{\xi,1}(\bR^n))^k}{k!}\bar\nu^{*k}(A)+t^{\frac n2}\sum_{k=1}^{\infty}\frac{(\nu_{\xi,1}(\bR^n))^k}{k!}\bar\nu^{*k}(A)\\
&= t\sum_{k=1}^{[\frac n2]}\frac{(\nu_{\xi,1}(\bR^n))^k}{k!}\bar\nu^{*k}(A)+\exp(\nu_{\xi,1}(\bR^n))\cdot t^{\frac n2}\P(\xi^1_1\in A).
\ea
We estimate $\sup_{|z|\le r_t}\mu^W_t(z)\P(\xi^1_t\in A)$ with $A=\{z\in\bR^n||y-z|<r_t\}$ and start with the first term on the right-hand side of \refq{432}. We combine the estimates \refq{430} and \refq{433} to obtain positive constants $K_1$, $K_2$, such that
\ba\label{1212}
\Big(\sup_{|z|\le r_t}\mu^W_t(z)\Big)\cdot\Big(&t\sum_{k=1}^{[\frac n2]}\frac{(\nu_{\xi,1}(\bR^n)t)^k}{k!}\bar\nu^{*k}(\{z\in\bR^n||y-z|<r_t\})\Big)\\
\le&\;K_1\cdot t^{-\frac n2}\cdot(t\,r_t^{n-1})\;\le\;K_2\cdot t^{\frac12} |\ln t|^{\frac{n-1}2}\;\le\; t^{\frac13}\;\le\;\exp(-\frack13\Lambda(\ln\Lambda)^{q})
\ea
for $\Lambda$ sufficiently large and all $t<\exp(-\Lambda(\ln\Lambda)^{q})$. This shows \eqref{BEDIIRn} with $\aleph\in(1-\frac1\alpha,q)$ for the first summand in \refq{432}. For the second summand in \refq{432} we combine the first limit of \refq{427}  for some $\aleph\in(1-\frac1{\alpha},1-\frac1{\alpha_\xi})$ with  \refq{430} to see that 
\ba\label{1213}
\lim_{\Lambda\rightarrow\infty}\sup_{|y|=\Lambda}\frac{\ln(\sup_{|z|\le r_t}\mu^W_t(z)\hspace{2pt}t^{\frac n2}\hspace{2pt}\P(\xi^1_1\in \{z\in\bR^n||y-z|<r_t\}))}{\Lambda(\ln\Lambda)^{\aleph}}=-\infty.
\ea
Combining estimates \refq{432}, \refq{1212} and \refq{1213} shows that the third summand on the right-hand side of \refq{429} is sufficiently small to satisfy \refq{BEDIIRn} for any $\aleph\in(1-\frac1{\alpha},q)$ in the last missing case $t\in(0,\exp(-\Lambda(\ln\Lambda)^{q}))$. Combining the estimates for these two cases of $t$ and by the choice of $q$ we obtain that the third summand on the right-hand side of \refq{429} is sufficiently small to satisfy~\refq{BEDIIRn} for any $\aleph\in(1-\frac1{\alpha}, 1-\frac1{\alpha_\xi})$.
This finishes the proof of \eqref{BEDIIRn} and of Lemma \ref{hinr_bed2}.

\bigskip 
\section{No LDP for the rescaled L\'evy process $(X^\e)_{\e>0}$}\label{s:LDPnonexistence}\setcounter{equation}{0}

Let $L$ be a L\'evy process with values in $\bR^n$ which satisfies Hypotheses~I and II. Let $r_\e$ be of regular variation with an index in $(-1,\infty)$, $T>0$ and $(X^\e_t)_{t\in[0,\infty)}$ with $X^\e_t=\e L_{r_\e t}$. Let $f$ and $g$ be defined as in Theorem~\ref{th:1Rn} and $S(\e):=f(g(\frac{\e^{-1}}{r_\e}))^{-1}\e$. By applying Theorem~\ref{th:1Rn}(i) it can easily be seen, that for every $t\in(0,\infty)$ the family $(X^\e_t)_{\e>0}$ satisfies a LDP according to the speed function $S$ and the rate function 
\be
I_t(x)=|x|.
\label{450Rn}
\ee
For $t=0$ the initial condition $X^\e_t=0$ implies the trivial LDP with the rate-function $I_o(x)=\infty$ for any $x\neq0$.

 If we assume, that $(X^\e)_{\e>0}$ is $S$-exponentially tight on $(\mathcal D_{[0,T],\bR^n},\mathcal J_1)$, then $(X^\e)_{\e>0}$ satisfies a LDP according to the speed function $S$ and equations \refq{eq:1} and \refq{eq:2} are applicable to calculate the associated rate function $I$. Let $|x|=1$, $T>0$ and $\varphi:[0,\infty)\rightarrow\bR^n$ with $\varphi(t):=x\textbf{1}_{[T,\infty)}(t)$. Then ~\refq{eq:1} and \refq{eq:2} together with \refq{450Rn} yield $I(\varphi)=1$.

\noindent By construction $\varphi$ is discontinuous in $t=T$. Hence there is a $\kappa_1>0$ such that for  $A:=\{\vartheta\in\mathcal D_{[0,\infty),\bR}\mid d(\varphi,\vartheta)<\kappa_1\}$ we have $\kappa_2:=\inf_{\vartheta\in A}\sup_{t\le 2T}|\vartheta(t)-\vartheta(t-)|>0$, where $d$ denotes the metric that induces the $\mathcal J_1$ topology. Consequently 
\ba
\lim_{\kappa\to0}\lime S(\e)\ln\P\big(d(\varphi,X^\e)<\kappa\big)\;\le\;\lime S(\e)\ln\P\big(d(\varphi,X^\e\big)<\kappa_1)&\\
\le\;\lime S(\e)\ln\P\Big(\sup_{t\le2r_\e T}\e|L_t-L_{t-}|>\kappa_2\Big)&\;=\;-\infty\;\neq\;-I(\varphi),
\ea
which contradicts the LDP.
Thus the assumption of $S$-exponential tightness of $(X^\e)_{\e>0}$ on $(\mathcal D_{[0,\infty),\bR},\mathcal J_1)$  is false, and no~LDP can exist for $(X^\e)_{\e>0}$ on $(\mathcal D_{[0,\infty),\bR},\mathcal J_1)$ with the speed function $S$.\\[2mm]
Furthermore, any speed function $S_o$, that would grant $S_o$-exponential tightness of $(X^\e)_{\e>0}$ on $(\mathcal D_{[0,\infty),\bR},\mathcal J_1)$ has to satisfy $\lime S_o(\e)S(\e)^{-1}=\infty$. Therefore, when trying to setup a LDP for $(X^\e_t)$ for any $t\in[0,\infty)$, we obtain
\be
\hspace{-1mm}\lime S_o(\e)\ln\P(|X^\e_t|\ge \Lambda)=-\infty,\mb{5mm}{thus}I_t(x)=-\lim_{\kappa\rightarrow0}\lime S_o(\e)\ln\P(|X^\e_t-x|\ge\kappa)=\infty\nonumber
\ee
for any $\Lambda,|x|>0$.  Applying \refq{eq:1} and \refq{eq:2} the resulting LDP for $(X^\e)_{\e>0}$  on $(\mathcal D_{[0,\infty),\bR},\mathcal J_1)$ has the rate function $I_o(\varphi)=\infty$ for any $\varphi\not\equiv0$. Obviously such a LDP is of no practical value.


\end{document}